\documentclass[11pt,a4paper]{amsart}
\usepackage{amssymb}
\usepackage{amsmath}
\usepackage{amscd}
\usepackage[francais]{babel}
\usepackage{latexsym}
\usepackage[T1]{fontenc}
\usepackage{epsfig}
\newcommand{\SW}{\mathrm{SW}}

\newcommand{\Inv}{\mathrm{Inv}_\epsilon}
\newcommand{\LLc}{\mathrm{L}^2c_1}
\newcommand{\PSL}{\mathrm{PSL}}
\newcommand{\PU}{\mathrm{PU}}
\newcommand{\Unitaire}{\mathrm{U}}
\newcommand{\pq}{{p^q}}
\newcommand{\R}{\mathbb{R}}
\newcommand{\rang}{{\mathrm{rang}}}
\newcommand{\trace}{{\mathrm{tr}}}
\newcommand{\gmod}{{\hat g}}
\newcommand{\gSigma}{{g ^\Sigma}}
\newcommand{\Amod}{{\hat A}}
\newcommand{\gmodt}{{\hat g_t}}
\newcommand{\gHH}{{g ^{\HH}}}

\newcommand{\bparab}{\mathcal{ B}}
\newcommand{\bparabarg}{\bparab_a^{\alpha_1,\alpha_2}}
\newcommand{\gmodj}{{\hat g_j}}
\newcommand{\gmodjj}{{\hat g^j}}
\newcommand{\gkahl}{ {g^\mathrm{K}} }

\newcommand{\degpar}{\mathrm{deg\,par\,}}
\newcommand{\gFS}{g^\mathrm{FS}}
\newcommand{\gDelta}{g^\mathrm{\Delta ^*}}
\newcommand{\Rat}{\mathbb{Q}}
\newcommand{\delbar}{{\overline\partial}}
\newcommand{\barSigma}{{\overline\Sigma}}
\newcommand{\dt}{{\partial_t}}
\newcommand{\dtheta}{{\partial_\theta}}
\newcommand{\dthetad}{{\partial_{\theta_2}}}
\newcommand{\ZZ}{\mathbb{Z}}
\newcommand{\C}{\mathbb{C}}
\newcommand{\Proj}{\mathbb{P}} 
\newcommand{\PP}{\mathbb{P}} 
\newcommand{\CP}{\mathbb{CP}^1} 
\newcommand{\CProj}{\mathbb{CP}} 
\newcommand{\HH}{{\mathbb{H}^2}} 
 
\newcommand{\HDR}{\mathrm{H}_{\mathrm{DR}}} 
 
\newcommand{\HLL}{\mathrm{H}_{L^2}} 
 
\newcommand{\HLLg}[1]{\mathrm{H}_{L^2(#1)}} 
\newcommand{\NLLg}{\|_{L^2(g)}} 
\newcommand{\NLLj}{\|_{L^2(g_j)}} 
\newcommand{\NLLjj}{\|_{L^2(g^j)}} 
 
\newcommand{\vol}{\mathrm{vol}} 
\newcommand{\volg}{\mathrm{vol} ^g}
\newcommand{\volgt}{\mathrm{vol} ^{g_t}}  
\newcommand{\volj}{\mathrm{vol}^{g_j}} 
\newcommand{\voljj}{\mathrm{vol}^{g^j}} 
 
\newcommand{\volss}{\mathrm{vol}^{\mathrm{\CProj^1}}} 
\newcommand{\ffond}{\mathbb{I}}

\newcommand{\Lie}{\mathcal{L}}
\newcommand{\Xtheta}{{\mathcal{X}_\theta}}

\newcommand{\aparabp}{asymptotique au  mod{\`e}le local}
\newcommand{\aparab}{\aparabp{} }
\newcommand{\jauge}{{\mathcal G}}

\newcommand{\Ric}{\mathrm{Ric}}
\newcommand{\Dirac}{\mathrm{D}}

\newcommand{\spinc}{spin^c}

\newcommand{\remarque}{\emph{Remarque}~: }
\newcommand{\remarques}{\emph{Remarques}~: }
\renewcommand{\phi}{\varphi} \renewcommand{\Im}{{\mathrm{Im}}}
\renewcommand{\Re}{{\mathrm{Re}}}

\newtheorem{prop}{Proposition}[section]
\newtheorem{cor}[prop]{Corollaire}
\newtheorem{theo}[prop]{Th{\'e}or{\`e}me}

\newtheorem{theosans}{Th{\'e}or{\`e}me}

\newtheorem{dfn}[prop]{D{\'e}finition}
\newtheorem{lemme}[prop]{Lemme}
\newenvironment{demo}{\emph{D{\'e}monstration}.---}{\hfill $\Box$ \medskip }
\newenvironment{demode}[1]{\emph{{D{\'e}monstration} #1.---}}{\hfill
$\Box$ \medskip} 
\date{}
\title[Surfaces k{\"a}hl{\'e}riennes de volume fini]{Surfaces
k{\"a}hl{\'e}riennes de volume fini et 
{\'e}quations de Seiberg-Witten}
\author{Yann Rollin}
\email{rollin@math.polytechnique.fr}
\begin{document}
{ \maketitle} \bigskip

\begin{center}
\small
\begin{tabular}{p{11cm}}
{\sc Abstract.} Let  $M=\PP(\mathcal E)$ be a ruled surface. We introduce 
metrics of finite volume on $M$ whose singularities are parametrized
by a parabolic structure over $\mathcal E$. Then, we generalise
results of Burns--de Bartolomeis and Le\,Brun, by showing that the
existence of a singular K\"ahler metric of finite volume
and constant non positive scalar curvature on $M$ is equivalent to the
parabolic polystability of $\mathcal E$~; moreover these metrics all
come from finite volume quotients of $\HH\times \CP$.
In order to prove the theorem, we must produce a solution of
Seiberg--Witten equations 
for a singular metric $g$. We use orbifold compactifications
$\overline M$ on which  we  approximate $g$ by a sequence of smooth
metrics~; the desired solution for $g$ is obtained as the limit of a
sequence of Seiberg--Witten solutions for these smooth
metrics.
\end{tabular}
\end{center}\bigskip

\section{Introduction}
L'existence de m{\'e}triques de K{\"a}hler {\`a} courbure scalaire constante
sur les surfaces complexes  est un probl{\`e}me ouvert dans lequel de
r{\'e}centes avanc{\'e}es mettent en {\'e}vidence l'importance de la notion de
stabilit{\'e} (\cite{BB}, \cite{L} et \cite{LS}).

Dans le cas d'une surface g{\'e}om{\'e}triquement r{\'e}gl{\'e}e de la
forme  $M=\Proj(\mathcal E)$, o{\`u} $\mathcal E $ est un fibr{\'e}
holomorphe au dessus d'une surface de Riemann hyperbolique compacte
$\barSigma$, Burns et de Bartolomeis ont d{\'e}montr{\'e} que les seules
m{\'e}triques k{\"a}hl{\'e}riennes {\`a} courbure scalaire $s=0$  sont des
produits locaux~; par le th{\'e}or{\`e}me de Narashiman--Seshadri,
cette condition  est {\'e}quivalente {\`a} la polystabilit{\'e} de $\mathcal E$.
Puis ce r{\'e}sultat {\`a} {\'e}t{\'e} g{\'e}n{\'e}ralis{\'e} par Le\,Brun au cas
$s<0$ en utilisant la th{\'e}orie de Seiberg--Witten.

Mehta et Seshadri {\'e}tendent les r{\'e}sultats
de \cite{NS} en volume  fini gr{\^a}ce {\`a}  la th{\'e}orie des
\emph{fibr{\'e}s paraboliques}  et une condition de 
\emph{polystabilit{\'e} parabolique} (pour une
d\'emonstration <<\`a la Donaldson>>  \cite{D} de leur r\'esultat
cf. \cite{Bq2}). D'apr{\`e}s 
\cite{MS} ces fibr{\'e}s  correspondent, en rang 2, aux cas o{\`u} la
surface r{\'e}gl{\'e}e $M$ est un quotient de volume fini
$\Sigma\times_\rho \CP$, provenant d'une repr\'esentation  $\rho
:\pi_1(\Sigma)\rightarrow \PU(2)$,  o\`u $\Sigma =\overline \Sigma \setminus
\{P_i\}$  est une surface de Riemann  de type hyperbolique obtenue
en enlevant un nombre fini de points appel{\'e}s \emph{points
paraboliques} \`a une surface de Riemann compacte.
En munissant le premier facteur de la
m{\'e}trique hyperbolique {\`a} courbure $-1$ et le second de la
m{\'e}trique de Fubini-Study {\`a} courbure $c>0$, on en d{\'e}duit une
m{\'e}trique k{\"a}hl{\'e}rienne $\gmod$ {\`a} courbure scalaire $s=2(c-1)$ sur $M$.
Dans des coordonn{\'e}es locales adapt{\'e}es $(t,\theta,u)\in \R\times
S^1\times \CP$ sur les bouts
de $M$, o\`u $u$ est une coordonn{\'e}e affine sur $\CP$,
la m{\'e}trique $\gmod$ est donn{\'e}e par 
\begin{equation} 
\label{exprmloc}
\gmod = dt^2+ e ^{-2t}d\theta ^2 +
\frac{4/c}{(1+|u|^2)^2}|du - i\alpha u d\theta|^2,
\end{equation}
o{\`u} $\alpha\in [0,1[$ est le \emph{poids} de cette singularit{\'e} appel{\'e}e \emph{bout parabolique}.
 Nous allons {\'e}tendre les r{\'e}sultats de \cite{BB} et \cite{L}
  {\`a} ce  cadre de 
  volume fini  en supposant que les m{\'e}triques consid{\'e}r{\'e}es 
 sont asymptotiques au sens $C^2$  au mod{\`e}le local d{\'e}fini
  par~(\ref{exprmloc}). 
\begin{theosans}
\label{theoremeA}
Soit $\mathcal E\rightarrow \overline\Sigma $  un fibr{\'e} parabolique
  holomorphe  tel que $\Sigma=\overline\Sigma\setminus\{P_i\}$
  soit hyperbolique. Soit $M=\Proj(\mathcal E)$ la surface complexe
  r{\'e}gl{\'e}e associ{\'e}e  restreinte au dessus de $\Sigma$. Alors
  $M$ admet une m{\'e}trique k{\"a}hl{\'e}rienne $\gkahl$ {\`a}
  courbure scalaire constante $s\leq 0$ 
asymptotique au mod{\`e}le local si et seulement si le fibr{\'e} $\mathcal E$ est paraboliquement
polystable.  Dans ce cas, la m{\'e}trique $\gkahl$  se d{\'e}duit {\`a} un biholomorphisme et une
constante pr{\`e}s du mod{\`e}le $(\Sigma\times_\rho \CProj^1,\gmod) $
o{\`u} $\rho :\pi_1(\Sigma)\rightarrow
\PU (2)$ est une repr{\'e}sentation associ{\'e}e au fibr{\'e}
parabolique polystable $\mathcal E$.
\end{theosans}
\emph{Remarque   : }
si $\gkahl$ admet au moins une
singularit{\'e}, son comportement asymptotique impose que
la constante du th{\'e}or{\`e}me soit {\'e}gale {\`a} $1$.\medskip

On rencontre
la difficult{\'e} majeure de la d{\'e}monstration dans le  cas $s<0$,
o{\`u} pour proc{\'e}der suivant la m{\'e}thode de
Le\,Brun on doit extraire une solution $(A,\psi)$, suffisamment
r{\'e}guli{\`e}re, des {\'e}quations de Seiberg-Witten
 \begin{eqnarray*}
   D_A\psi &= &0\\
F_A ^+ & = & q(\psi) 
 \end{eqnarray*}
pour une m{\'e}trique $g$  \aparab et la structure $\spinc$ canonique
de $M$.
 En l'absence de singularit{\'e}s, Le\;Brun  utilise la
th{\'e}orie des invariants de Seiberg-Witten sur les vari{\'e}t{\'e}s
compactes, ce qui n'est pas le cas de $M$.
Pour traiter ce probl{\`e}me, on montre dans la section~\ref{sec:skrm} que 
lorsque les poids des singularit{\'e}s sont rationnels, $M$ admet une
compactification orbifold $\overline M = M\cup D$, o{\`u} $D$ est une
r{\'e}union de diviseurs de la forme $\CP/\ZZ_q$, sur laquelle on peut
approximer $g$ par une suite de m{\'e}triques lisses $g_j$.
On d{\'e}finit alors les {\'e}quations de Seiberg-Witten perturb{\'e}es sur
$\overline M$ pour chaque m{\'e}trique $g_j$ par
 \begin{eqnarray*}
   D_A\psi &= &0\\
(F_A + 2i\pi \varpi_j)^+ & = & q(\psi),
 \end{eqnarray*}
o{\`u} les perturbation $\varpi_j$ se concentrent vers le courant
d'int{\'e}gration sur $D$.
En calculant l'invariant des {\'e}quations de Seiberg--Witten qui
d{\'e}pend d'une certaine \emph{condition de chambre} pour la m{\'e}trique
$g_j$, nous  obtiendrons une suite de solutions
$(A_j,\psi_j)$ que nous ferons converger vers une solution $(A,\psi)$
des {\'e}quations non perturb{\'e}e pour la m{\'e}trique limite $g$
gr{\^a}ce au 
th{\'e}or{\`e}me suivant~: 
\begin{theosans}
\label{theoremeC}
Soit $g$ une m{\'e}trique sur $M$ asymptotique au mod{\`e}le local
  $\gmod$ avec des poids
  $\alpha(P_i)$ rationnels et soit $g_j$ la suite
  d'approximation de $g$ sur $\overline M$.
  Soit $(A_j,\psi_j)$ une suite de solutions des {\'e}quations de
  Seiberg-Witten perturb{\'e}es associ{\'e}es {\`a} la structure $\spinc$
  canonique induite par la structure complexe et aux
  m{\'e}triques $g_j$ sur $\overline M$. Alors quitte {\`a} faire des changements de jauge
  et {\`a} extraire une sous-suite, $(A_j,\psi_j)$ converge au
  sens $C^\infty$ sur tout compact de $M$ vers une solution des
  {\'e}quations non perturb{\'e}es $(A,\psi)$ pour la m{\'e}trique $g$
  v{\'e}rifiant
  \begin{itemize}
  \item $A = \Amod +a$ avec $a\in L^2_1(g)$,
    \item $\psi\in L^2_1(g)$
      relativement {\`a} $\nabla_{\Amod}$,
  \end{itemize}
 o\`u  $\Amod$ est une
  connexion induite par le mod\`ele local k\"ahl\'erien $\gmod$ sur le
  fibr\'e d\'eterminant  $L=K_M ^{-1}$ de la structure $\spinc$ canonique.
\end{theosans}
\remarques
\begin{itemize}
\item  On pourra se r\'ef\'erer \`a~\cite{Bq3} pour une autre
application des \'equations de Seiberg--Witten \`a des m\'etriques
d'Einstein de volume fini.
\item Il serait naturel que l'existence de la solution $(A,\psi)$
obtenue par  ce
th\'eor\`eme soit assur\'ee par une th\'eorie de Seiberg--Witten
d\'evelopp\'ee directement pour
les m\'etriques asymptotiques au mod\`ele local. Notons  que la
construction de l'espace des modules correspondant et la question de
sa compacit\'e  posent des probl\`emes techniques importants que
notre m\'ethode ne r\'esout pas (cf. par exemple \cite{KM} pour
des m\'etriques asymptotiquement plates).
\item
La perturbation 
$\varpi_j$ est  n{\'e}cessaire pour que la connexion $A$ soit
d{\'e}finie sur le <<bon fibr{\'e}>> et s'explique au niveau des {\'e}quations
par le fait que la m{\'e}trique $g$ appara{\^\i}t comme une limite de
m{\'e}triques sur la vari{\'e}t{\'e} compacte moins une \emph{bulle} de
courbure positive. 
\end{itemize}

La premi{\`e}re {\'e}tape  dans la d{\'e}monstration du
th{\'e}or{\`e}me~\ref{theoremeC} (cf. section \ref{chapitre2}) consiste {\`a} d{\'e}velopper une th{\'e}orie de Hodge pour
les m{\'e}triques $g$ et $g_j$
via des in{\'e}galit{\'e}s de Poincar{\'e} uniformes. {\`A} l'aide d'un lemme
de Poincar{\'e} local pr{\`e}s de $D$ pour la cohomologie $L^2$, on
d{\'e}montre un isomorphisme 
$$\HLL^k(M)\simeq \HDR^k(\overline M),
$$
puis on d{\'e}montre pour $k=1$ ou $2$ que les repr{\'e}sentant
$g_j$-harmoniques d'une 
classe de cohomologie convergent vers le repr{\'e}sentant $g$-harmonique
$L^2$ sur tout compact de $M$. {\`A} partir de l\`a on poss{\`e}de tous les
outils n{\'e}cessaires pour faire converger les connexions, pour
d\'emontrer  le th{\'e}or{\`e}me au
\S\ref{subdemotheoa} puis pour calculer l'invariant de Seiberg--Witten
des m{\'e}triques $g_j$ au \S\ref{subcalcinv}.
Dans le cas de poids irrationnels, on n'a plus de compactification orbifold
ad{\'e}quate  et on doit faire une convergence <<en deux temps>> en
commen{\c c}ant par approximer les poids irrationnels par des  poids
rationnels (cf. \S\ref{secirrat}).
\medskip

\emph{Sur les \'eclatements de
  surfaces r\'egl\'ees.} La structure  parabolique,  en chaque  point
  $P$ \emph{non trivial}   d'un fibr\'e   $\mathcal E$
  (cf. \S\ref{subparab} pour la d\'efinition),  
  d\'etermine une droite complexe de $\mathcal E_P$, donc un point
  $Q$ de
  de la surface compacte $\widehat M = \PP(\mathcal E)_{\barSigma}$ et
  un poids $\alpha = \alpha_2-\alpha_1$.
 Le probl{\`e}me d'existence  de
m{\'e}triques de K{\"a}hler {\`a} courbure scalaire constante sur 
l'\'eclatement $\widetilde M$ de $\widehat M$  aux points $Q_i$
 a {\'e}t{\'e} abord{\'e}  par Le\,Brun et Singer 
lorsque $\widehat M$ poss{\`e}de des champs de vecteurs
holomorphes p{\'e}riodiques et $s=0$~:  suivant \cite{LS}, de telles
  m{\'e}triques existent  si et seulement si $\mathcal E$ est
  paraboliquement polystable. 

La forme de K\"ahler $\omega $ de la m\'etrique $\gmod$,
bien que singuli\`ere sur $\widehat M$, correctement
interpr{\'e}t{\'e}e comme un courant positif (cf. \cite{CG}), se
rel{\`e}ve en un 
courant positif $\widehat \omega$ appel{\'e} \emph{transform{\'e}e
stricte} de $\omega$ 
sur l'{\'e}clatement $\widetilde M$ tel que   
$$ 0< \alpha  = \frac {\widehat \omega\cdot   [ E ] }{
  \widehat \omega\cdot [F] },
$$ 
avec $ E$ le diviseurs exceptionel au point $Q$ et $F$ une fibre g\'en\'erique.
Cette identit\'e  est pr\'ecis\'ement celle  v{\'e}rifi{\'e}e par les
classes de K{\"a}hler 
consid{\'e}r{\'e}es dans \cite{LS} et 
ceci   constitue une indication
 suppl{\'e}mentaire mettant en {\'e}vidence le lien  entre la
notion de stabilit{\'e} et le probl{\`e}me d'existence de
m{\'e}triques k{\"a}hl{\'e}riennes 
{\`a} courbure scalaire constante sur les surfaces complexes.\medskip

Je tiens \`a remercier tout particuli\`erement Olivier Biquard pour
l'ensemble des  discussions que j'ai eues avec lui sur ce sujet.

\section{Surfaces k{\"a}hl{\'e}riennes r{\'e}gl{\'e}es mod{\`e}les}\bigskip
\label{sec:skrm}
\subsection{Un exemple fondamental}
\label{subsub:exemple}
Voici tout d'abord une famille de surfaces complexes
r{\'e}gl{\'e}es construites {\`a} partir de  repr{\'e}sentations 
unitaires du groupe fondamental d'une surface de Riemann de volume fini. Ces exemples sont essentiels dans la th{\'e}orie des fibr{\'e}s
parabolique stables comme nous le verrons au \S\ref{secstable} o{\`u} nous citerons le th{\'e}or{\`e}me important de
Mehta-Seshadri~; en
outre ces exemples  poss{\`e}dent des m{\'e}triques k{\"a}hl{\'e}riennes
<<mod{\`e}les>> {\`a}
courbure scalaire constante avec des singularit{\'e}s que
nous {\'e}tudierons pr{\'e}cis{\'e}ment au \S\ref{submodlocmetr} et que
nous appellerons  \emph{bouts paraboliques}.

\subsubsection{Surfaces de Riemann hyperboliques de volume fini} 
Soit $\Gamma$ un sous groupe discret de $\PSL(2,\R)$ 
agissant librement et avec covolume fini sur le demi plan de Poincar{\'e}
$\HH=\{\xi\in\C \;/\; \Im \xi >0 \} \; ;$
le quotient est une surface de Riemann $\Sigma=\HH/\Gamma$, de groupe 
fondamental $\Gamma$, munie de la 
m{\'e}trique k{\"a}hl{\'e}rienne $\gSigma$ compl{\`e}te de volume fini {\`a}
courbure $-1$   induite par la 
m{\'e}trique de Lobachevsky.
Le groupe $\Gamma$ agit {\'e}galement sur le bord {\`a} l'infini
$\partial_\infty \HH$ du demi-plan de Poincar{\'e}. Puisque $\Gamma$ agit avec covolume
fini,  le stabilisateur  d'un point du bord est soit trivial, soit,
pour un nombre fini de points $P$   
appel{\'e}s  \emph{points paraboliques}, {\'e}gal {\`a} $\langle \tau\rangle\simeq \ZZ$, o{\`u} $\tau$ est un {\'e}l{\'e}ment
parabolique de $\PSL(2,\R)$.
On peut  donc, quitte {\`a}  conjuguer $\tau$ par une  
homographie, supposer qu'il est donn{\'e} par
$\tau : \xi\mapsto\xi+u$  o{\`u} $\xi \in \HH$ et $u\in \R$ ; pour
simplifier  on 
supposera m{\^e}me que $u=2\pi$. Au voisinage de $P$, le quotient
$\HH/\Gamma$ est isomorphe {\`a} $I_a/\langle \tau \rangle $,  o{\`u}
$I_a=\{\xi\in\C \;/\; \Im\xi>a\}$ avec $a>0$  suffisamment 
grand. 
On d{\'e}finit alors un isomorphisme entre  le disque {\'e}point{\'e} $\Delta_a^*=\{z\in \C \; / \;
 0<|z|< e^{-a}\}$ et  $I_a/\langle \tau \rangle $  par 
\begin{eqnarray}
   I_a/\langle \tau \rangle & \longrightarrow & \Delta_a^* \\
  \xi & \longmapsto & z= e^{i\xi} \;; \nonumber
\end{eqnarray}
en utilisant le plongement holomorphe  $\Delta_a^*\subset\Delta_a $, ce qui revient a ajouter le point $P$
 correspondant {\`a} $0$ dans le mod{\`e}le du disque,
 on obtient finalement la compactification
 holomorphe $\overline \Sigma = \Sigma \cup \{P_i\}_{1\leq i \leq k} $.

En partant de la m{\'e}trique  de Lobachevsky donn{\'e}e par~:
$$ g^\HH= \frac{|d\xi| ^2}{|\xi-\bar \xi|^2} ,
$$
on calcule ais{\'e}ment la m{\'e}trique induite sur le disque {\'e}point{\'e}~:
$$ g^{\Delta^*}= \frac {|dz|^2}{|z| ^2\ln^2|z|}.
$$
Cette m{\'e}trique poss{\`e}de une singularit{\'e} en $0$ appel{\'e}e \emph{cusp}.
Nous utiliserons souvent un autre syst{\`e}me de coordonn{\'e}es
locales $(t,\theta)\in \R \times \R/2\pi\ZZ$ sur les bouts de $\Sigma$
d{\'e}fini par 
$z=re^{i\theta}$ et $t =\ln
(- \ln |z|)$~;
 dans ces coordonn{\'e}es 
\begin{equation}
  \label{cusp}
g^{\Delta^*}=dt^2+ e^{-2t}d\theta^2.  
\end{equation}

R{\'e}ciproquement (cf. \cite{Fo}), si $\barSigma$ est une surface de Riemann
compacte avec $k$ points paraboliques marqu{\'e}s $(P_i)$, alors
$\Sigma= \overline \Sigma\setminus \{P_i\}$  est 
 hyperbolique si et seulement si   $2g(\overline
\Sigma)-2+k>0$ et dans ce cas elle admet une unique
m{\'e}trique k{\"a}hl{\'e}rienne $g^\Sigma$ de courbure $-1$ avec des cusps aux
points paraboliques.

\subsubsection{Repr{\'e}sentations et surfaces r{\'e}gl{\'e}es}
Soit $\rho:\Gamma\rightarrow\Unitaire(2)$  une repr{\'e}sentation
unitaire du groupe fondamental de
$\Sigma$. On en
d{\'e}duit une action de $\Gamma$ sur $\HH\times \C^2$ donn{\'e}e par
$ \sigma \cdot (\xi, u )= (\sigma \cdot \xi, \rho(\sigma) \cdot u )$ pour
$\sigma\in\Gamma $ et $(\xi,\mu)\in \HH\times \C^2$. 
Au quotient, on obtient un fibr{\'e} vectoriel holomorphe 
$\mathcal E = \HH \times_\rho \C^2 $ au dessus de $\Sigma$. On en
d{\'e}duit une surface complexe r{\'e}gl{\'e}e $M=\PP(\mathcal{E}) =
\HH\times _\rho \CP$ o{\`u} $\rho$ d{\'e}signe par abus de langage la
repr{\'e}sentation projective induite.
La sph{\`e}re de Riemann $\CProj^1$ est munie d'une m{\'e}trique
k{\"a}hl{\'e}rienne {\`a} courbure 
 sectionnelle $c>0$, appel{\'e}e  m{\'e}trique de 
\emph{Fubini-Study}, d{\'e}finie dans des coordonn{\'e}es homog{\`e}nes
$[u:v]$ par  
$$ \omega^{\mathrm{FS}} = \frac {2i}c \partial\delbar \ln(|u|^2+|v|^2) \; ;
$$
dans la carte $v=1$, on obtient alors
$$ \gFS = \frac {4/c}{(1+|u|^2)^2}|du|^2 .
$$
 La m{\'e}trique produit $\gHH\oplus\gFS$ sur $\HH\times
\CProj^1$ est 
k{\"a}hl{\'e}rienne {\`a} courbure scalaire  constante $s=2(c-1)$, car c'est un produit riemannien de m{\'e}triques k{\"a}hl{\'e}riennes 
{\`a} courbures scalaires constantes.  
Le groupe d'isom{\'e}tries de
$(\CProj^1,\gFS)$ est {\'e}gal {\`a} $\PU(2)$, donc $\Gamma$ agit
isom{\'e}triquement sur $\HH\times \CP$. On en d{\'e}duit par quotient une
m{\'e}trique k{\"a}hl{\'e}rienne {\`a} courbure scalaire constante $\gmod$ sur
$M$ appel{\'e}e \emph{m{\'e}trique mod{\`e}le}.
En outre, $\gmod$ est de volume fini car c'est un produit local donc
$\vol(M)=\vol 
^{{\Sigma}}(\Sigma)\times \vol^{\mathrm {FS}}(\CProj^1)<\infty$. 

\subsection{Fibr{\'e}s paraboliques}
\label{subparab}
Une question naturelle se pose maintenant si nous voulons fabriquer
des m{\'e}triques de K{\"a}hler de volume fini {\`a} courbure scalaire
constante sur les surfaces complexes r{\'e}gl{\'e}es comme dans
l'exemple pr{\'e}c{\'e}dent~:
comment reconna{\^\i}{}tre celles qui sont obtenues {\`a} partir  de
repr{\'e}sentation projective unitaire  du groupe fondamental d'une
surface de Riemann de volume fini~? La r{\'e}ponse {\`a} nous est
donn{\'e}e par la th{\'e}orie des fibr{\'e}s paraboliques stables
d{\'e}velopp{\'e}e par Mehta et Seshadri~; nous rappelons dans cette
partie des
\'el\'ements de la th\'eorie, mais pour de plus amples d\'etails nous
renvoyons le lecteur \`a~\cite{MS}.

\subsubsection{Structures paraboliques}
\label{subsub:structureparab}
Le fibr{\'e} hermitien plat $\mathcal E=\HH\times_\rho\C^2$  
admet un prolongement holomorphe au dessus de
$\barSigma$ sur lequel la m\'etrique hermitienne est
singuli\`ere~; en consid\'erant les morphismes born\'es de fibr\'es
paraboliques, on voit appara\^\i{}tre naturellement la notion de
\emph{structure parabolique}.
\begin{dfn}[Structures paraboliques]
\label{dfnfibreparab}
Une  structure parabolique sur un  fibr{\'e} vectoriel complexe $\mathcal E\rightarrow
\overline\Sigma$ au dessus d'une surface de Riemann compacte
est la donn{\'e}e 
d'un ensemble fini points paraboliques $P_i\in\overline \Sigma$,
de filtrations au dessus des points paraboliques
$$ \mathcal E_{P_i}=\mathcal F^1_{P_i}\supsetneq \dots \supsetneq \mathcal F^{k_i+1}_{P_i}=0, $$
et de poids associ{\'e}s {\`a} la filtration
$ 0 \leq \alpha_1(P_i)< \dots < \alpha_{k_i} (P_i) < 1 \; ;$
on appelle multiplicit{\'e} d'un poids $\alpha_j(P_i)$ l'entier 
$d_j(P_i)=\dim
\mathcal F^j_{P_i} -\dim \mathcal F^{j+1}_{P_i}$. Un fibr{\'e} munit
d'une structure parabolique est appel{\'e} un fibr{\'e} parabolique. \medskip
\end{dfn}
\remarque pour les fibr{\'e}s paraboliques de rang $2$, on a soit $k_i=2$, soit
$k_i=1$. Dans ce dernier cas, la structure parabolique au point $P_i$
est r{\'e}duite {\`a} une filtration triviale et {\`a} un unique poids de
multiplicit{\'e} $2$. Un tel point parabolique est appel{\'e} un  
\emph{point trivial} de la structure parabolique. \medskip

Dans une base orthonorm\'ee $(\epsilon_1,\epsilon_2)$ bien choisie, on a
 \begin{equation}
\label{exprtau}
\rho(\tau)=\left ( \begin{array}{cc} e^{2i\pi\alpha_1} & 0 \\ 0 & e^{2i\pi\alpha_2}
  \end{array}\right )\quad \mbox{ avec $0
\leq\alpha_1\leq\alpha_2< 1$};
\end{equation}
on en d\'eduit une base de sections locales holomorphes $(s_1,s_2)$ de
\mbox{$ \HH\times_\rho\CP$} en $P$ telle
que la m\'etrique hermitienne s'\'ecrive
\begin{equation}
\label{exprmetrherm}
h=\left ( \begin{array}{cc} |z|^{2\alpha_1} & 0 \\ 0 & |z|^{2\alpha_2}
  \end{array}\right ).
\end{equation}
On d{\'e}finit alors une structure parabolique sur $\mathcal E$ comme suit~: en un point parabolique $P$
de $\overline \Sigma$ tel que  $\alpha_1 < \alpha_2$, on pose $\mathcal
F^2_P = \C s_2(P)$, $\alpha_1(P)=\alpha_1$ et
$\alpha_2(P)=\alpha_2$.
Si $\alpha_1=\alpha_2$, on d{\'e}cide que $P$ est un point trivial de
la structure parabolique avec un unique poids $\alpha_1(P)=\alpha_1$.

Avec cette d\'efinition de structure parabolique, on peut voir
  facilement que les isomorphisme de $\mathcal E$ born\'es relativement
  \`a la m\'etrique hermitienne sont exactement ceux 
qui respectent les filtrations au dessus des  points paraboliques.

\subsubsection{Stabilit{\'e} parabolique et th\'eor\`eme de Mehta--Seshadri}
 \label{secstable}
Dans le \S\ref{secstable}, $\mathcal E \rightarrow \barSigma $
d{\'e}signera un fibr{\'e} parabolique holomorphe.
Son \emph{degr{\'e} parabolique}
et sa \emph{pente} sont d{\'e}finis respectivement par 
$$\degpar \mathcal E= \deg \mathcal E +\sum_{i,j} d_j(P_i)\cdot \alpha_j(P_i),\quad
\mu( \mathcal E)= \frac {\degpar
  \mathcal E}{\rang  \;\mathcal E}.$$

Si $\mathcal E=\mathcal L_1\oplus \mathcal  L_2$, est une somme de
deux fibr{\'e}s 
paraboliques en droites complexes,alors $\mathcal E$ poss{\`e}de une
structure parabolique induite.
On dira que $\mathcal E$ est \emph{paraboliquement d{\'e}composable}
s'il admet une d{\'e}composition holomorphe en une somme directe de
fibr{\'e}s 
parabolique  $\mathcal E=\mathcal L_1\oplus
\mathcal L_2$ qui est de plus compatible avec la structure parabolique de $\mathcal
E$.

Un fibr{\'e} parabolique, $\mathcal E\rightarrow \overline \Sigma$, de rang $2$,
est dit \emph{paraboliquement stable}
(resp.  \emph{paraboliquement semi-stable})
si tout sous-fibr{\'e} holomorphe en droites complexes $\mathcal L$ muni de la
structure parabolique induite v\'erifie  $ \mu(\mathcal L)< \mu(\mathcal E)$ (resp. $ \mu(\mathcal L) \leq  \mu(\mathcal E)$).
Le fibr{\'e} $\mathcal E$ sera dit \emph{paraboliquement polystable} s'il est
soit paraboliquement stable, soit {\`a} la fois paraboliquement
semi-stable et paraboliquement d{\'e}composable.

Lorsque $\mathcal E\rightarrow \barSigma$ est un fibr\'e parabolique
tel que $\Sigma=\barSigma \setminus \{P_i\}$ soit hyperbolique, on a 
le th{\'e}or{\`e}me important de Mehta et Seshadri
(cf. \cite{MS} th. 4.1) d'apr\`es lequel~: $\mathcal E$ est polystable
 de pente $\mu = 0$ si et seulement si il est plat et donn\'e \`a
isomorphisme parabolique 
pr\`es par une repr\'esentation $\rho :\pi_1(\Sigma)\rightarrow
\Unitaire(2)$.

Lorsque $\degpar \mathcal E \neq 0$, on peut se ramener au cas du
degr{\'e} parabolique nul de la fa{\c c}on suivante~: soit
$\mathcal L \rightarrow \overline \Sigma$,  un fibr{\'e} en droites complexes de
degr{\'e} tr{\`e}s n{\'e}gatif, de sorte  que $\degpar (\mathcal E\otimes \mathcal L) <
0$. Notons que les surfaces $\Proj(\mathcal E)$ et $ \Proj(\mathcal E\otimes \mathcal
L)$ sont isomorphes. Dans
un deuxi{\`e}me temps, on  ajoute des points paraboliques triviaux $Q_k\in \Sigma$
 {\`a} la structure parabolique de sorte que le nouveau fibr{\'e}
parabolique  $\mathcal E '$ v{\'e}rifie
$\degpar \mathcal E '=\degpar (\mathcal E\otimes\mathcal L)
+\sum_{k}2 \alpha_1(Q_k) =0$.
En choisissant bien les poids, on obtient un fibr{\'e} de
degr{\'e} parabolique~$0$ et 
(cf. \cite{Se}) le fibr{\'e} $\mathcal E '$ est paraboliquement stable
(resp. semi-stable, polystable) si et seulement si le fibr{\'e}
$\mathcal E$ est paraboliquement stable (resp. semi-stable,
polystable). 

\subsubsection{Stabilit{\'e} et m{\'e}triques mod{\`e}les}
Soit $\mathcal E\rightarrow \barSigma$ un fibr{\'e} parabolique
polystable, et $M=\Proj(\mathcal
E)_{|\Sigma}$ la surface r{\'e}gl{\'e}e associ{\'e}e.
  En supposant $\Sigma$ hyperbolique, on sait 
  que $\Sigma$ admet une unique m{\'e}trique k{\"a}hl{\'e}rienne $g^\Sigma$ {\`a}
  courbure sectionnelle $-1$ avec des cusps aux points $P_i$. On peut
quitte {\`a} ajouter des points paraboliques triviaux $Q_i$ supposer que $\degpar \mathcal E=0$ sans changer la
  stabilit{\'e} parabolique du fibr{\'e}. 

Alors, d'apr{\`e}s le th{\'e}or{\`e}me de Mehta et Seshadri, il
existe  
une repr{\'e}sentation projective unitaire
$\rho :\pi_1(\Sigma\setminus\{Q_i\})\rightarrow \Unitaire(2)$ telle que
$\mathcal E \simeq
\HH\times_\rho \C^2$. 
Il est alors facile
de voir que la m{\'e}trique hermitienne, singuli{\`e}re aux points $Q_i$,
induit une m{\'e}trique de Fubini-Study sur la fibre qui est  lisse sur
$M$ tout enti{\`e}re. 
On obtient alors suivant \S\ref{subsub:exemple} une m{\'e}trique k{\"a}hl{\'e}rienne $\gmod$
{\`a} courbure scalaire constante sur $M$.

\remarque la construction que nous venons de faire montre que nous
pouvons facilement nous ramener au cas des fibr{\'e}s paraboliques de
pente  $0$ sans ajouter de singularit{\'e}s {\`a} la m{\'e}trique 
$\gmod$. Par la suite, nous supposerons pour simplifier que tous les fibr{\'e}s
paraboliques sont de pente nulle. De plus, nous ferons comme si les
structure paraboliques {\'e}taient concentr{\'e}es au dessus d'un unique
point parabolique, toutes les constructions s'{\'e}tendant trivialement
au cas de plusieurs points.

R{\'e}ciproquement, il est naturel de se demander si l'existence d'une
m{\'e}trique k{\"a}hl{\'e}rienne {\`a} courbure scalaire constante  sur $M$
avec des singularit{\'e}s de volume fini implique la stabilit{\'e} du
fibr{\'e} parabolique  $\mathcal E$ de d{\'e}part. Nous allons r{\'e}soudre
cette question dans le th{\'e}or{\`e}me $\ref{theoA}$, mais  pour
l'{\'e}noncer, nous devons introduire au pr{\'e}alable une notion de \emph{mod{\`e}le local}.

\subsection{Le mod{\`e}le local}
\label{submodlocmetr}
Nous allons maintenant nous concentrer sur l'{\'e}tude locale de la
m{\'e}trique $\gmod$ au voisinage d'un point parabolique.

La m{\'e}trique $\gmod$ que nous avons d{\'e}finie sur la surface
r{\'e}gl{\'e}e $M$ tout enti{\`e}re lorsque $\mathcal E$ est paraboliquement
polystable, peut {\^e}tre d{\'e}finie localement en g{\'e}n{\'e}ral sur un
fibr{\'e} parabolique. Pour commencer, on choisit une trivialisation
locale de
$\mathcal E$ par
une base  de sections  holomorphes $(s_1,s_2)$ d{\'e}finies au voisinage
d'un point
parabolique $P$ compatible avec la filtration parabolique. Soit $z$
une coordonn{\'e}e 
locale holomorphe telle que $z(P)=0$ et telle que la m{\'e}trique
sur $\Sigma$ s'exprime sous la forme
$$ \gSigma=\frac{|dz|^2}{|z|^2\ln^2|z|}.
$$
On d{\'e}finit alors dans la base $(s_1,s_2)$  une m{\'e}trique hermitienne sur $\mathcal E$
\emph{adapt{\'e}e}  {\`a} la structure parabolique par (\ref{exprmetrherm}).
On en d{\'e}duit une base orthonorm{\'e}e de sections locales
en $P$ donn{\'e}es par
$ e_i =s_i /
|z|^{\alpha_i}$ 
et on calcule la connexion de Chern associ{\'e}e {\`a} $h$ dans la trivialisation $(s_i)$~:
\begin{equation}
\label{exprconnchern}
\nabla ^{h}=d + h^{-1}\partial h= d+\left (\begin{array}{cc}\alpha_1 &0 \\
    0&\alpha_2 \end{array}\right) \frac{dz}z,
\end{equation}
ce qui donne dans les coordonn{\'e}es induites par la base orthonormales $(e_i)$, 
$$ \nabla ^{h} = d+i\left (\begin{array}{cc}\alpha_1 &0 \\ 0
    &\alpha_2 \end{array}\right ) d\theta.$$
Cette connexion est plate, avec une singularit{\'e} logarithmique au
point $P$. On retrouve dans cette expression locale l'holonomie
de la connexion de Chern {\'e}gale {\`a} $diag(\exp(-2i\pi\alpha_1),\exp(-2i\pi\alpha_2))$ autour du point $P$. 
En notant $H_x$ la distribution horizontale de la connexion induite
sur le $\CProj^1$-fibr{\'e} $M\rightarrow\Sigma$, on a une
d{\'e}composition 
$T_x M = H_x \oplus T_x
\Proj(\mathcal{E}_{\pi(x)})$.
Comme la
connexion de Chern respecte par d{\'e}finition la 
m{\'e}trique hermitienne sur $\mathcal E$, on en d{\'e}duit que la
distribution $(H_x)$ est invariante par les isom{\'e}tries de la
m{\'e}trique de Fubini-Study sur les fibres
de $M\rightarrow \Sigma$. 
  
On d{\'e}finit alors localement la m{\'e}trique riemannienne $\gmod$  en prenant la
m{\'e}trique de Fubini-Study sur l'espace tangent aux fibres et
l'image r{\'e}ciproque de  $g^\Sigma$ sur $H_x\simeq T_x\Sigma$. Comme de plus la connexion est
plate, la distribution horizontale $H_x$ est int{\'e}grable, et $\gmod$
appara{\^\i}t comme un produit riemannien local $\HH\times \CProj
^1$. On dit que $\gmod$ est \emph{adapt{\'e}e} {\`a} la structure
parabolique de $\mathcal E$.

\subsubsection{Rev{\^e}tement et mod{\`e}le local}
\label{subsub:revmodloc}
En reprenant les notations de \S \ref{subsub:structureparab},
nous disposons de plusieurs coordonn{\'e}es 
locales d{\'e}finies au voisinage de $P$~: la coordonn{\'e}e $\xi\in \HH$, la
coordonn{\'e}e $z\in \Delta_a$ d{\'e}finie au voisinage de $P$,
les coordonn{\'e}es $(\xi, ( \tilde u, \tilde v))$ sur $\HH\times \C^2$
induites par la base canonique $(\epsilon_i)$ de $\C^2$ et les
coordonn{\'e}es $(z ,(u,v))$ sur $\mathcal E$ induites par  la base
orthonorm{\'e}e de
sections locales  $(e_i)$ de
$\mathcal E$.
 On en
d{\'e}duit des coordonn{\'e}es locales homog{\`e}nes sur $\HH\times \CProj^1$ et
$\Proj(\mathcal E)$ et on a le rev{\^e}tement holomorphe riemannien 
\begin{eqnarray}
\label{exprp}
p: I_a\times \CProj^1 &\rightarrow & \Delta_a^* \times \CProj^1 \\
(\xi,[\tilde u : \tilde v ]) &\mapsto & (e^{i\xi}, [e^{-i\alpha_1 \Re\xi }\tilde
u: e^{-i\alpha_2 \Re \xi}\tilde v]), \nonumber 
\end{eqnarray}
avec
\begin{equation}
\label{expractiontau}
\tau\cdot (\xi,[\tilde u : \tilde v ]) = ( \xi+2\pi, [e^{2i\pi\alpha_1}\tilde
u: e^{2i\pi \alpha_2 }\tilde v])
\end{equation}
qui engendre le groupe d'automorphismes 
$\langle\tau\rangle\simeq \ZZ$ du rev{\^e}tement $p$.

Si on se place par exemple dans les coordonn{\'e}es locales de projection
st{\'e}r{\'e}ographique $\tilde v=1$ et $v=1$, on obtient l'expression de $p$
$$p(\xi,\tilde u) =(e^{i\xi}, e^{i \alpha \Re\xi}\tilde u) \mbox{ o{\`u} 
$\alpha=\alpha_2-\alpha_1$.} $$
 On en d{\'e}duit que $e^{i\alpha \Re\xi}d\tilde u =p^*du - i\alpha u
d\theta$. Or la m{\'e}trique standard sur le rev{\^e}tement s'{\'e}crit 
$$p^*\gmod =g^\HH \oplus \gFS =  \frac{|d\xi|^2}{|\xi-\bar \xi|^2}  + \frac {4/c}{(1+|u|^2)^2 }\left |
  d\tilde u \right |^2.$$
d'o{\`u} 
\begin{equation}\label{g0}
 \gmod=\frac{|dz|^ 2}{ (|z|\ln |z|)^2} + \frac {4/c}{(1+|u|^2)^2 }\left |
  du - i\alpha u d\theta  \right |^2. 
\end{equation}
On peut maintenant d{\'e}gager la notion de \emph{bout parabolique} sur une vari{\'e}t{\'e}.
\begin{dfn}
\label{dfnparab}
Le groupe $\ZZ\simeq \langle \tau \rangle $ a une action libre holomorphe et isom{\'e}trique sur la vari{\'e}t{\'e}
complexe k{\"a}hl{\'e}rienne $\HH\times \CProj^1$ munie de sa m{\'e}trique
standard {\`a} courbure scalaire
$2(c-1)$ d{\'e}finie par  
\begin{equation}
\label{expractionZZ}
\tau\cdot (\xi,[\tilde u : \tilde v])= (\xi,[\tilde u e^{2i\pi\alpha_1 }
: \tilde v e^{2i\pi\alpha_2 }]).
\end{equation}
Le quotient $(I_a \times \CProj^1 )/\ZZ$ est par cons{\'e}quent une vari{\'e}t{\'e}
k{\"a}hl{\'e}rienne munie de la m{\'e}trique de K{\"a}hler quotient
$\gmod$  {\`a} courbure scalaire  $2(c-1)$ et de volume fini.   
Ce quotient sera appel{\'e} le bout
parabolique associ{\'e} aux poids $\alpha_1$,
$\alpha_2$. Nous le noterons $\bparabarg$, ou bien $\bparab$ lorsque
aucune ambigu{\"\i}t{\'e} n'est possible.
\end{dfn}

Par d{\'e}finition, un bout parabolique est muni de \emph{coordonn{\'e}es adapt{\'e}es}
 $(z,[u:v])\in \Delta ^ *_a\times \CProj ^1$ dans lesquelles le
 rev{\^e}tement canonique $p : I_a \times \CProj ^ 1\rightarrow \bparab$
est donn{\'e} par (\ref{exprp}) et la m{\'e}trique $\hat g$ est alors donn{\'e}e localement par
(\ref{g0}). 

Dans des coordonn{\'e}es adapt{\'e}es, la projection suivant le
premier facteur, 
$$\pi:\bparab\simeq
\Delta_a^*\times \CProj^1\rightarrow \Delta_a^*$$
 est holomorphe et donne au bout parabolique une structure de $\CProj ^1$-fibr{\'e} et on a clairement le diagramme
commutatif

\begin{picture}(0,0)(-200,-12)
\put(-50,-32){\vector(2,-1){50}}
\put(-34,-55){\scriptsize $\pi_1$}
\end{picture}
$$
\begin{CD}
I_a\times \CProj^1 @>p>>  \bparab  \\
@.                                      @VV\pi V \\
               @.                  \Delta_a ^*
\end{CD}
$$
o{\`u} $\pi_1$ est induite par le rev{\^e}tement canonique $I_a\rightarrow
\Delta_a ^*$.

\remarques
\label{qqnot}
la forme volume de la m{\'e}trique de Fubini-Study sur $\HH\times \CP$  passe au quotient
et donne une $2$ forme ferm{\'e}e sur le bout parabolique que nous noterons $\volss$.

En rempla{\c c}ant la coordonn{\'e}e $z$ par $(t,\theta)$, la m{\'e}trique mod{\`e}le appara{\^\i}{}t comme un produit tordu
\begin{equation}
\label{g_1}
\gmod= dt^2+e^{-2t}d\theta^2 + \frac{4/c}{(1+|u|^2)^2}\left
|du-i\alpha u d\theta\right |^2  = dt ^2 + \gmod_t,
\end{equation} 
o{\`u} $\gmod_t$ est une m{\'e}trique sur $N\simeq S^1\times \CProj^1$ que
nous appellerons la \emph{tranche} du bout parabolique.

La m{\'e}trique sur $\HH\times \CProj^1$ d{\'e}g{\'e}n{\`e}re dans la direction
$\partial_x= dx^\sharp$ o{\`u} $\xi=x+iy$ et o{\`u} $dx^ \sharp$ d{\'e}signe le
dual suivant la m{\'e}trique de $dx$. Ce champ de vecteurs est invariant
sous l'action de $\tau$. Par cons{\'e}quent $\partial_x=dx^\sharp $ passe
au quotient et nous donne un champ de vecteurs not{\'e} $\Xtheta$ sur le bout parabolique.
Comme $p^*d\theta = dx$, on calcule explicitement en introduisant les
coordonn{\'e}es polaires $u=\rho e ^{i\theta_2}$ 
\begin{equation}
\label{exprXtheta}
\Xtheta =d\theta ^\sharp =\partial_\theta + \alpha\partial_{\theta_2}.
\end{equation}
On v{\'e}rifie ais{\'e}ment que la m{\'e}trique $g^\HH+\gFS $ est invariante
suivant $\partial_x$ et on en
d{\'e}duit que $\gmod$ est invariante suivant $\Xtheta$~:
\begin{equation}
\label{exprLieXtheta}
\Lie_\Xtheta \gmod = 0.
\end{equation}

Nous introduisons maintenant sur une vari{\'e}t{\'e} {\`a} bouts
paraboliques, une classe de m{\'e}triques \emph{asymptotiques
 au mod{\`e}le local} $\gmod$ (au sens $C^2$). Dans le cas o{\`u} $X=\PP(\mathcal E)_{|\Sigma}$ est une
surface r{\'e}gl{\'e}e obtenue {\`a} partir d'un fibr{\'e} parabolique, on choisira
toujours un  mod{\`e}le local
$\gmod$ adapt{\'e} {\`a} la structure parabolique.
\begin{dfn}
 Si une vari{\'e}t{\'e} $X$  poss{\`e}de des bouts paraboliques, on en d{\'e}duit un mod{\`e}le local de
  m{\'e}trique que
  nous noterons $\gmod$. Nous appellerons
  m{\'e}trique <<\aparabp>>, toute m{\'e}trique $g$ sur $X$ 
  telle que $g=\gmod+o(\gmod)$, o{\`u} $o(\gmod)$ tend  uniform{\'e}ment
  vers $0$ en norme $C^2$ lorsqu'on s'approche de 
l'infini  sur bout parabolique.
\end{dfn}
C'est une  condition
technique raisonnable qui va nous permettre de d{\'e}montrer au
\S\ref{subdemotheoa} le
th{\'e}or{\`e}me suivant, dont on d\'eduit imm\'ediatement le
th\'eor\`eme~\ref{theoremeA}.
\begin{theo}
\label{theoA}
  Soit $\mathcal E\rightarrow \overline \Sigma$ un fibr{\'e}
parabolique holomorphe de
  rang $2$ sur une surface de Riemann compacte $\overline \Sigma$.
Soit $(P_i)_{ 1\leq i\leq k }$ l'ensemble 
des  points paraboliques de $\overline \Sigma$. Nous supposons
   en outre $\Sigma=\barSigma\setminus\{P_i\}$ hyperbolique.
Si la  surface complexe r{\'e}gl{\'e}e  associ{\'e}e $M=\PP(\mathcal E)_\Sigma $  admet une m{\'e}trique k{\"a}hl{\'e}rienne $\gkahl$ {\`a}
  courbure scalaire constante $s=2(c-1)\leq 0 $ \aparab $\hat g$,
   alors $\mathcal E$ est paraboliquement polystable.
De plus une telle m{\'e}trique provient
  {\`a} un biholomorphisme pr{\`e}s, du mod{\`e}le
  $(\Sigma\times_\rho\CProj^1 , \gmod)$ o{\`u}
  $\rho:\pi_1(\Sigma)\rightarrow \PU(2)$ est une repr{\'e}sentation
associ{\'e}e au fibr{\'e} parabolique polystable
$\mathcal E$. 
\end{theo} 

\subsubsection{Compactification du bout parabolique rationnel}
Nous faisons maintenant l'hypoth{\`e}se $\alpha_2-\alpha_1=r/q\in\Rat$,
o{\`u} $r$ et $q$ sont premiers entre eux.
Dans ce cas, l'action de $\tau$ d{\'e}finie par
(\ref{expractionZZ}) v{\'e}rifie
$$\tau ^q\cdot (\xi,[\tilde u :\tilde v]) = (\xi+2\pi q,[\tilde u :\tilde v]).
$$
On en d{\'e}duit que $ (I_a\times \CProj ^1)/\langle \tau ^q\rangle
\simeq  (I_a /\langle \tau ^q\rangle) \times \CProj ^1$. 
On vient donc de montrer que le rev{\^e}tement $p: I_a\times \CProj ^1\rightarrow
\bparab$ se factorise par un rev{\^e}tement $\pq$ {\`a} $q$ feuillets dont
le groupe d'automorphismes est {\'e}gal {\`a} $\langle \tau \rangle / \langle \tau ^q \rangle \simeq
\ZZ_q$. On a donc un diagramme commutatif holomorphe

\begin{picture}(0,0)(-190,-12)
\put(-50,-32){\vector(2,-1){60}}
\put(-30,-55){\scriptsize $p$}
\end{picture}
$$
\begin{CD}
I_a\times \CProj^1 @>\tilde p>> I_a/\langle \tau ^q\rangle \times \CProj ^1 \\
@.                                      @VV\pq V \\
               @.                 \bparab
\end{CD}
$$
o{\`u} $\tilde p$ d{\'e}signe la projection canonique sur le
rev{\^e}tement quotient. On note
\begin{eqnarray*}
\pi ^ q : \Delta_{a/q}  \simeq I_a
/\langle \tau^q \rangle\cup\{\infty\} & \rightarrow & \Delta_a \simeq I_a
/\langle \tau \rangle \cup\{\infty\} \\
z &\mapsto & z^q 
\end{eqnarray*}  
le rev{\^e}tement holomorphe ramifi{\'e} {\`a} $q$ feuillets du disque par lui-m{\^e}me.
Le groupe d'automorphisme $\ZZ_q\simeq \langle
\tau\rangle /  \langle
\tau ^ q\rangle$ du  rev{\^e}tement $\pi ^ q $ agit par 
$\tau\cdot \xi=\xi+2\pi$ soit $ \tau\cdot z = \exp({\frac{2i\pi}q})z$.
On d{\'e}duit de $\pi ^q$ une projection holomorphe not{\'e}e {\'e}galement
$$\pi ^q:  I_a /\langle \tau^q  \rangle \cup \{\infty\} \times \CProj ^1  \rightarrow  \Delta_a .$$
Par construction, il est clair que $\pi \circ p=  \pi_1  = \pi ^q
\circ \tilde p$~; on en d{\'e}duit le lemme~:

\begin{lemme}[Structure du bout parabolique rationnel] 
  \label{lemmemodrat}
  Supposons que $\alpha=\alpha_2
  -\alpha_1=r/q \in \Rat$, avec $r $ et $q$ premiers entre eux.  
  Le rev{\^e}tement standard $p:I_a\times\CProj^1 \rightarrow \bparab$
  du bout parabolique $\bparabarg$
  se factorise par un rev{\^e}tement  {\`a} $q$- feuillets $\pq$
  et on a le diagramme holomorphe commutatif

\begin{picture}(0,0)(-160,-12)
\put(-50,-32){\vector(2,-1){60}}
\put(-34,-55){\scriptsize $p$}
\put(50,-32){\vector(2,-1){45}}
\put(83,-43){\scriptsize $\pi ^q$}
\end{picture}
$$
\begin{CD}
I_a\times \CProj^1 @>\tilde p>> I_a/\langle \tau ^q\rangle \times \CProj ^1 \\
@.                                      @VV\pq V \\
               @.                 \bparab @>\pi>> \Delta_a ^*,
\end{CD}
$$
\end{lemme}
L'action de $\ZZ_q$ sur $I_a/\langle \tau ^q \rangle \times \CProj ^1
\simeq \Delta_{a/q} ^*\times \CProj ^1 $ s'{\'e}tend clairement en une
action holomorphe sur $\Delta_{a/q}\times \CProj ^1$~: pour le voir, il
suffit de l'{\'e}crire dans des coordonn{\'e}es adapt{\'e}es $(z,[\tilde u
:\tilde v])$ 
 $$\tau\cdot (z,[\tilde u:\tilde v])=(ze^{\frac{2i\pi}q}, [\tilde u
 e^{2i\pi\alpha_1}:\tilde v  e^{2i\pi\alpha_2}]).$$

Alors, comme nous l'avons d{\'e}j{\`a} remarqu{\'e}, $\pi ^q$ se prolonge de
fa{\c c}on holomorphe en $\pi ^q: \Delta_{a/q}\times \CProj ^1
\rightarrow \Delta_a$~; au quotient, on obtient un orbifibr{\'e}
holomorphe en
droites projectives 
$$\overline \bparab = (\Delta_{a/q}\times \CProj ^1) /
\ZZ_q = \bparab \cup (\CProj^1/\ZZ_q) \stackrel {\pi}\rightarrow
\Delta_a .$$
Reformulons cette discussion dans le corollaire suivant.
\begin{cor}
\label{corcomprat}
   Si $\alpha= \alpha_2 -\alpha_1=r/q \in \Rat$ avec $r$ et $q$
premiers entre eux, le $\CProj ^1$-fibr{\'e} $\bparab\rightarrow
\Delta_a ^*$ admet un prolongement holomorphe en un orbifibr{\'e}
$$\overline \bparab = \bparab \cup (\CProj^1/\ZZ_q) \stackrel
{\pi}\rightarrow \Delta_a.$$ 
Cet orbifibr{\'e} poss{\`e}de un rev{\^e}tement
ramifi{\'e} {\`a} $q$ feuillets $\pq : \Delta_{a/q}\times
\CProj^1\rightarrow \overline \bparab$ tel que $(\pq) ^* \gmod$ est
{\'e}gale {\`a} la m{\'e}trique produit standard sur $\Delta_{a/q}^*\times \CProj ^1$.
La fibre $\pi ^{-1}(0)\simeq \CProj^1/\ZZ_q$
sera appel{\'e}e  <<diviseur {\`a} l'infini du bout parabolique>>.
\end{cor}

En conclusion, une vari{\'e}t{\'e} $M$ avec des bouts paraboliques
v{\'e}rifiant la condition de rationalit{\'e} admet une compactification
orbifold $\overline M$. Si $M$ est de plus une vari{\'e}t{\'e} complexe, la
compactification est holomorphe.

\section{Cohomologie $L^2$.}\medskip

\label{chapitre2}
Dans toute cette partie $M$ est une vari{\'e}t{\'e} $M$ munie de bouts
paraboliques et 
on {\'e}tudie la cohomologie $\HLL^* (M)$  
d{\'e}finie {\`a} l'aide du complexe des formes diff{\'e}rentielles $L^2(\gmod)$  sur $M$.
Notons que pour une m{\'e}trique $g$ asymptotique au mod{\`e}le local
$\gmod$ les formes diff{\'e}rentielles $L^2$ relativement {\`a} $g$ et
$\gmod$ co{\"\i}{}ncident ainsi donc que les espaces de cohomologies $L^2$ associ{\'e}s.

Nous allons d{\'e}montrer un lemme de Poincar{\'e} <<{\`a} l'infini>>
sur le bout 
parabolique de $M$ pour la cohomologie $L^2$, afin de la calculer
dans le cas des surfaces r{\'e}gl{\'e}es.

Lorsque les poids des bouts paraboliques v{\'e}rifient la condition
$\alpha_2-\alpha_1\in \Rat$, nous d{\'e}finirons au \S \ref{seccomp} sur
la compactification orbifold $\overline M$ des m{\'e}triques lisses $g_j$ qui 
approximent la m{\'e}trique singuli{\`e}re $g$. 
On d{\'e}veloppera alors la th{\'e}orie de Hodge en montrant des  r{\'e}sultats
de convergence sur tout compact de $M$ des repr{\'e}sentants $g_j$
harmoniques d'une classe de cohomologie vers le repr{\'e}sentant
$g-$harmonique $L^2$
(cf. prop. \ref{convcompact} et \ref{propconvharm}).

\subsection{Approximation par des  m{\'e}triques lisses}
\label {seccomp}
Dans \S\ref{seccomp} nous supposons que les bout paraboliques de $M$ v{\'e}rifient la
condition de rationalit{\'e} et nous expliquons comment approximer la m{\'e}trique $g$,
singuli{\`e}re sur $\overline M$, par une suite de m{\'e}triques $(g_j)$ 
lisses sur $\overline M$.
\subsubsection{Approximation du mod{\`e}le local $\gmod$}
\label{subappmodloc}
D'apr{\`e}s le lemme \ref{lemmemodrat} et le corollaire \ref{corcomprat},
tout bout parabolique compactifi{\'e} de $\overline M$ admet un
rev{\^e}tement ramifi{\'e} holomorphe {\`a} $q$ feuillets 
$\pq : \Delta_{a/q} \times \CProj^1\rightarrow \overline \bparab$,
 tel que $(\pq) ^* \gmod$ est la m{\'e}trique standard
$\gDelta\oplus \gFS$. On va approximer $\gmod$
en rempla{\c c}ant la m{\'e}trique singuli{\`e}re $\gDelta$ par une m{\'e}trique lisse
sur $\Delta$, invariante sous l'action de $\ZZ_q$. Puisque 
$$\gDelta = \frac {|dz|^2}{(|z|\ln|z|)^2};
$$
on peut par exemple multiplier cette m{\'e}trique par un facteur
conforme singulier de la forme $\lambda (|z|)$, valant $(|z|\ln|z|)^2$
au voisinage de l'origine et $1$ loin de l'origine. Alors
la m{\'e}trique k{\"a}hl{\'e}rienne $g_\lambda=\lambda(|z|)\gDelta$
est lisse sur $\Delta$~; cette op{\'e}ration consiste simplement {\`a} recoller un disque 
plat au bout du cusp et $g_\lambda$  est $\ZZ_q$-invariante car la
fonction $\lambda(|z|)$ est elle m{\^e}me invariante. Alors la m{\'e}trique
$g_\lambda \oplus \gFS$ passe au quotient et d{\'e}finit une m{\'e}trique
k{\"a}hl{\'e}rienne lisse sur $\overline\bparab$.
Si on suit les g{\'e}od{\'e}siques de  $g_\lambda$ normales {\`a} un cercle, on en
d{\'e}duit une nouvelle coordonn{\'e}e locale
$t'$ telle que $t=t'$ l{\`a} o{\`u} $\lambda(|z|)=1$ et $g_\lambda = dt^{'2} +
\phi (t')^2 d\theta ^2$. De fa{\c c}on {\'e}quivalente, on peut donc d{\'e}finir ce
type de m{\'e}trique en se donnant la fonction $\phi(t')$. 

Soit la
m{\'e}trique 
$$g^\Delta_j =dt^2+\phi_j^2 d\theta ^2, $$
o{\`u}  $\phi_j$ est d{\'e}finie comme suit~: pour $t\leq j+1$, on pose $\phi_j=e
^{-t}$. Puis pour $j+1\leq t \leq j+1+\delta_j$, on fait brutalement d{\'e}cro{\^\i}tre $\phi_j''$ de $e ^{-(j+1)}$  jusqu'{\`a} $0$
pour un $\delta_j$ que nous ferons par la suite tendre vers $0$. Enfin
pour $j+1+\delta_j\leq t\leq j+1+\epsilon$, on peut  quitte {\`a} rendre
 $\phi_j''$ tr{\`e}s n{\'e}gative faire  passer $\phi_j'$ dans cet
 intervalle  de $-e
 ^{-(j+1)}$ {\`a} $-1$. 

On a 
 recoll{\'e} la m{\'e}trique la m{\'e}trique de cusp (pour $t\leq j+1$)  
{\`a} la m{\'e}trique plate. En effet si $\phi_j'=-1$, alors pour un $T_j$
l{\'e}g{\`e}rement sup{\'e}rieur {\`a} $j+1+\epsilon$ bien choisit, on a
$\phi_j=T_j-t$.  Comme le montre le changement
de variable $t'=T_j-t$, la m{\'e}trique $g ^\Delta_j$ se recolle {\`a} la
m{\'e}trique plate sur le disque~; il en r{\'e}sulte que $g_j^\Delta $ est en
r{\'e}alit{\'e} lisse sur le disque $\Delta_{a/q}$ tout entier.
On v{\'e}rifie que quitte {\`a} choisir $\delta_j$ suffisamment petit, on
 peut supposer que le
rapport $-  \phi_j' /\phi_j$  est croissant, et reste donc minor{\'e} par $1$.
Nous r{\'e}sumons cette construction dans le  dessin ci-dessous.
\bigskip

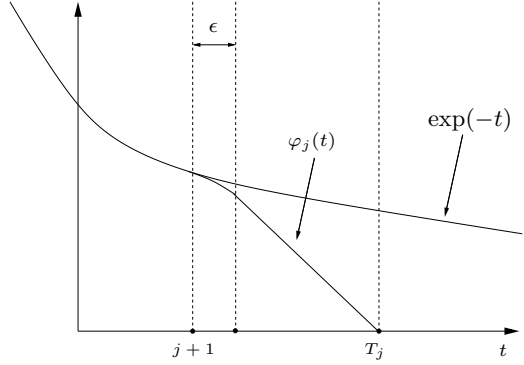
\begin{center}       
\begin{picture}(0,0)%
\epsfig{file=phij.pstex}%
\end{picture}%
\setlength{\unitlength}{1184sp}%
\begingroup\makeatletter\ifx\SetFigFont\undefined%
\gdef\SetFigFont#1#2#3#4#5{%
  \reset@font\fontsize{#1}{#2pt}%
  \fontfamily{#3}\fontseries{#4}\fontshape{#5}%
  \selectfont}%
\fi\endgroup%
\begin{picture}(10749,7443)(1564,-7192)
\put(9001,-7111){\makebox(0,0)[lb]{\smash{\SetFigFont{6}{7.2}{\familydefault}{\mddefault}{\updefault}$T_j$}}}
\put(11821,-7156){\makebox(0,0)[lb]{\smash{\SetFigFont{7}{8.4}{\familydefault}{\mddefault}{\updefault}$t$}}}
\put(10336,-2341){\makebox(0,0)[lb]{\smash{\SetFigFont{9}{10.8}{\familydefault}{\mddefault}{\updefault}$\exp(-t)$}}}
\put(7396,-2761){\makebox(0,0)[lb]{\smash{\SetFigFont{7}{8.4}{\familydefault}{\mddefault}{\updefault}$\phi_j(t)$}}}
\put(4981,-7111){\makebox(0,0)[lb]{\smash{\SetFigFont{6}{7.2}{\familydefault}{\mddefault}{\updefault}$j+1$}}}
\put(5746,-376){\makebox(0,0)[lb]{\smash{\SetFigFont{8}{9.6}{\rmdefault}{\mddefault}{\updefault}$\epsilon$}}}
\end{picture}

\end{center}  
\medskip

Finalement la m{\'e}trique $g^\Delta_j \oplus \gFS$ est lisse et $\ZZ_q$-invariante
sur $\Delta_{a/q}\times\CProj^1$. On en d{\'e}duit une m{\'e}trique $\gmodj$
lisse sur $\overline \bparab$
qui s'{\'e}crit dans les  coordonn{\'e}es $(t,\theta,u)$ 
\begin{equation}
  \label{g1}
\gmodj=dt^2 + \phi_j ^2 d\theta ^2  + \frac {4/c}{(1+|u|^2)^2 }\left |
  du-i\alpha u d\theta \right |^2.  
\end{equation}

\subsubsection{Approximation d'une m{\'e}trique asymptotique au mod{\`e}le}
\label{subapproxmetraparab}
{\'E}tant donn{\'e}e une m{\'e}trique riemannienne $g$ sur $M$, \aparab sur le
bout parabolique de $M$, on
approxime la m{\'e}trique $g$ par des m{\'e}triques  $g_j$ lisses sur la
compactification $\overline M$  comme suit~:
soit $\chi(t)$ une fonction lisse  telle que 
\begin{itemize}
\item $\chi$ est croissante,
\item $\chi(t)=0$ pour $t \leq 0$,
\item $\chi(t)=1$ pour $t\geq \frac 12$.
\end{itemize}
Notons $\chi_j (t)= \chi(t-j)$, puis posons
\begin{equation}
g_j=(1-\chi_j) g +\chi_j \gmodj.
\end{equation}
Alors puisque $g$ est \aparabp, les m{\'e}triques $g_j$ et $\gmod_j$
sont {\'e}galement tr{\`e}s proche au sens $C^2$ sur la zone de recollement
$j\leq t \leq j+1$, et on en d{\'e}duit en particulier que 
\begin{equation}
\label{exprgjaparab}
\sup_{t\geq j } \left \{ |\gmodj-g_j|+ |\nabla ^\gmodj - \nabla ^{g_j}| + |R ^\gmodj - R
^{g_j}| \right \} \rightarrow 0 \mbox{ lorsque } j\rightarrow \infty.
\end{equation}

Le mod{\`e}le local de  m{\'e}trique $\gmod$ est  {\`a} courbure  born{\'e}e 
puisque c'est un produit local model{\'e} sur  $\HH\times
\CProj^1$. 
Donc toute m{\'e}trique \aparab  est {\`a} courbure born{\'e}e. 

Lorsqu'on approxime la m{\'e}trique $\gmod$ par des m{\'e}trique $\gmodj$
lisses sur la compactification $\overline M$, on cr{\'e}e une <<bulle>> de
courbure tr{\`e}s positive. Pour calculer la courbure de $\gmodj$, il nous suffit de calculer celle du facteur modifi{\'e}~:
$$g^\Delta_j =d t^2+\phi_j^2 d\theta ^2. $$
La connexion associ{\'e}e est explicite et on a~:
\begin{alignat*}{2}
&\nabla_\dt\dt =  0 , \qquad &   &\nabla_\dt (\phi_j^{-1}\dtheta) = 0, \\
&\nabla_{(\phi_j^{-1}\dtheta)} \dt = \frac {\dt \phi_j}{\phi_j} \phi_j
^{-1} \dtheta, \qquad &
& \nabla_{(\phi_j^{-1}\dtheta)} \phi_j ^{-1} \dtheta = - \frac {\dt
  \phi_j}{\phi_j}\dt \; ;
\end{alignat*}
on en d{\'e}duit que la courbure sectionnelle de $g^\Delta_j$ est {\'e}gale
{\`a}  
\begin{equation*}
  K^{g^\Delta _j}(\dt,\phi_j^{-1} \dtheta ) = -\frac{\dt^2\phi_j}{\phi_j}.
\end{equation*}
Cette courbure par d{\'e}finition de $\phi_j$ sup{\'e}rieure {\`a} $-1$,
tr{\`e}s positive sur l'anneau $t\in [j+1, j+1+\epsilon ]$ et
uniform{\'e}ment born{\'e}e en dehors.  
Il en r{\'e}sulte  que la courbure
scalaire de $\gmodj= \pq_*(g^\Delta_j+\gFS)$ est donn{\'e}e par 
$$   s^\gmodj =2(c - \frac{\dt^2\phi_j}{\phi_j}) , $$
et que $   s_b^\gmodj := s^\gmodj + 2 \chi_j\frac{\dt^2\phi_j}{\phi_j} $
 est uniform{\'e}ment born{\'e}. 
Puisque la m{\'e}trique $g$ est \aparab $\gmod$, on en d{\'e}duit d'apr{\`e}s
(\ref{exprgjaparab}) un r{\'e}sultat analogue
pour les m{\'e}triques $g_j$~:
\begin{lemme}
\label{courburescalaire}
  La courbure scalaire de la m{\'e}trique $g_j$ s'{\'e}crit 
  \begin{equation}
  s^{g_j}=s^{g_j}_b-2\chi_j\frac{\dt^2\phi_j}{\phi_j},
  \end{equation}
o{\`u} $s^{g_j}_b$ est uniform{\'e}ment born{\'e}.
\end{lemme}
Par construction des m{\'e}trique $\gmodj$, nous avons contract{\'e} le volume initial de $\gmod$. On en d{\'e}duit le lemme technique suivant~:
\begin{lemme}[Contr{\^o}le du volume]
  \label{volume}
Quel que soit $\epsilon>0 $, il existe un compact $K$ de M, tel que
$\volj(M\setminus K)\leq\epsilon$ et $\vol^g(M\setminus
K)\leq\epsilon$. Le volume des m{\'e}triques $g_j$ est donc uniform{\'e}ment born{\'e}.
\end{lemme}

\subsection{Lemme de Poincar{\'e} {\`a} l'infini}
\label{seccohompoinc}
Dans tout le \S\ref{seccohompoinc}, on ne suppose plus en g{\'e}n{\'e}ral que la
condition de rationalit{\'e} des poids est v{\'e}rifi{\'e}e. Nous commen{\c c}ons par
d{\'e}montrer des estim{\'e}es sur les fonctions   ne
d{\'e}pendant que de 
la courbure moyenne des m{\'e}triques. 

\subsubsection{Courbure moyenne}
\label{subsub:courbmoy}
\begin{dfn}
  Soit  $g$ une m{\'e}trique riemannienne sur $[t_1,t_2]\times N$ de la
  forme $g=dt^2+g_t$ o{\`u} $g_t$ est une m{\'e}trique sur $N$ d{\'e}pendant de
  $t$. On d{\'e}finit la courbure  moyenne de la m{\'e}trique $g$ comme la fonction
$$H= \frac 13 \trace_g \ffond = -  \frac 13 \frac{\dt \vol^{g_t}}{\vol^{g_t}},
$$
o{\`u} $\ffond$ d{\'e}signe la deuxi{\`e}me forme fondamentale de la tranche
$N$ et o{\`u} $\vol^{g_t}$ est la forme volume associ{\'e}e {\`a} la m{\'e}trique
$g_t$. Par commodit{\'e}, on d{\'e}finit {\'e}galement la fonction $h=\frac 32 H.$
\end{dfn}
D'apr{\`e}s  (\ref{g0}) et (\ref{g1}) les m{\'e}triques  $\gmod$ et
$\gmodj$ sont des produits tordus v{\'e}rifiant respectivement $h=1/2$
et 
$$ h= -\frac 12 \phi_j^{-1}\dt \phi_j,$$
qui est par construction de $\phi_j$ minor{\'e} par $1/2$. Pr{\`e}s de
$t=T_j$, $\phi_j= T_j - t$ et $h= \frac 12  (T_j-t)^{-1}$ donc $h$
tend vers $+\infty$ quand $t\rightarrow T_j$.

D{\`e}s que l'on a un contr{\^o}le du type $h\geq h_0 >0$,
nous obtenons le lemme suivant qui
s'applique en particulier aux m{\'e}triques $\gmod$ et $\gmodj$~:
\begin{lemme}
\label{courburemoyenne}
Soit $g=dt^2+g_t$ une m{\'e}trique riemannienne sur $[t_1,t_2]\times
N$. Supposons qu'il existe $\delta,h_0\in \R$ telles que
$g$ v{\'e}rifie $h \geq h_0 > \delta$   ou $\delta > h_0 \geq h$.
Alors pour toute fonction $f$, on a 
\begin{equation}
\begin{split}
  \int_{[t_1,t_2]\times N}  \!\!\!\! \!\!\!\! |e^{\delta t}\dt f|^2
 \vol^g 
 \geq (h_0-\delta)^2 & \int_{[t_1,t_2]\times N}   \!\!\!\! |fe^{\delta
 t}|^2\vol^g \\ 
+(h_0-\delta ) \Big  \{ \int_{t=t_2 }  \!\!\!\!\! & |fe^{\delta t}|^2\vol^g
- \int_{t=t_1}\!\!\!\!\!  |fe^{\delta t}|^2\vol^g \Big \}.
\end{split}
\end{equation}
\end{lemme}
\begin{demo}
Voir \cite{Bq3} lemme 4.1, pour la d\'emonstration dans le cas $h>h_0>0$. En g\'en\'eral, la
d\'emontration est identique en int\'egrant $fe^{\delta t}$ par
parties cette fois-ci.
\end{demo}
\subsubsection{Identification de la cohomologie $L^2$}
\label{subPoincare}
Le \S \ref{subsub:courbmoy} fournit des outils suffisants pour
obtenir le lemme de Poincar{\'e} {\`a} l'infini pour la cohomologie $L^2$ suivant~:
\begin{lemme}[Lemme de Poincar{\'e} $L^2$ {\`a} l'infini]
\label{lemmePoincare}
Soit $\gamma$, une $k$-forme ferm{\'e}e $L^2$ sur la vari{\'e}t{\'e} {\`a} bouts
paraboliques $M$. 
\begin{itemize}
\item Pour $k=1,3$ ou $4$, il existe une $k-1$-forme $\beta$, $L^2$ d{\'e}finie sur le bout
  parabolique de $M$ telle que $\gamma=d\beta$ sur le bout parabolique.
\item   Dans le cas o{\`u} $k=2$, il existe $\lambda\in \R$ et une
$1$-forme $\beta$ $L^2$ d{\'e}finie sur le bout parabolique de $M$ tels que
$ \gamma = \lambda \volss + d\beta \mbox{ sur le bout parabolique.}$
\end{itemize}
\end{lemme}
Rappelons que $\volss$ est la $2$-forme ferm{\'e}e induite sur les bouts
parabolique par la m{\'e}trique de Fubini-Study (cf. \ref{subsub:revmodloc}).
Admettons pour l'instant ce lemme qui
va nous permettre de calculer la cohomologie $L^2$ de $M$.

Dans le cas o{\`u} $\alpha_2 - \alpha_1 \in \Rat$, on dispose {\'e}galement d'un lemme de
Poincar{\'e} analogue pour la cohomologie de de Rham orbifold sur le bout parabolique de $\overline M$. 
Rappelons qu'une forme diff{\'e}rentielle lisse au sens orbifold sur le bout
parabolique compactifi{\'e} $\overline \bparab$ de $\overline M$ est
donn{\'e}e par une forme diff{\'e}rentielle lisse et $\ZZ_q-$invariante sur le
rev{\^e}tement ramifi{\'e}.
On remarque alors que si $\gamma = d\beta $ est $\ZZ_q$-invariante sur
$\Delta_a\times \CProj^1$ avec $\beta$ a priori non invariante, alors
il existe une forme $\ZZ_q$-invariante $\beta'$ obtenue en moyennant $\beta$ sous l'action du groupe, telle que $\gamma=d\beta'$.
 On en d{\'e}duit 
le lemme de Poincar{\'e} suivant pour la cohomologie de
De Rham de l'orbifold $\overline M $~:
\begin{lemme}
\label{lemmepoincareorbifold}
Soit $\gamma$, une $k$-forme lisse ferm{\'e}e  sur $\overline M$. 
\begin{itemize}
\item Pour $k=1,3$ ou $4$, il existe une $k-1$-forme lisse $\beta$ 
telle que $\gamma=d\beta$ sur le bout parabolique de $\overline M$.
\item   Dans le cas o{\`u} $k=2$, il existe un r{\'e}el $\lambda$ et une
$1$-forme lisse $\beta$  sur  $\overline M$ tels que
$ \gamma = \lambda \volss + d\beta$  sur le bout parabolique.
\end{itemize}
\end{lemme}

Soit $\overline M'$ la compactification de $M$ obtenue en changeant
la valeurs des poids $\alpha_1(P_i),\alpha_2(P_i)$ de la structure
parabolique en $\alpha_1'(P_i)$ et $\alpha_2'(P_i)$
tels que $\alpha'(P_i)=\alpha_1'(P_i)-\alpha_2'(P_i)\in \Rat$.
En utilisant le lemme de Poincar{\'e} local,  on d{\'e}duit le
corollaire suivant.
\begin{cor}
\label{cohomologie}
On a un isomorphisme entre les espaces de cohomologie
$$ \HLLg{g}^*(M)\simeq \HDR^*(\overline M').$$
\end{cor}
\begin{demo}
D{\'e}finissons  un isomorphisme pr{\'e}servant le cup-produit
$F:\HLL^*(M)\rightarrow
   \HDR^*(\overline M')$ comme suit~:
soit $b \in \HLL^k(M)$. Pour $k=0$, on a clairement $\HLL ^0(M)\simeq \R$ puisque la m{\'e}trique est de volume fini. Par ailleurs $\HDR ^0(M)\simeq \R$ et on d{\'e}finit $F$ par cet isomorphisme.
Si $k=1,3$ ou $4$, on peut d'apr{\`e}s le lemme de Poincar{\'e}
 \ref{lemmePoincare}, choisir un repr{\'e}sentant $\gamma$ {\`a} support compact
 dans $M$ de  $b$ et on pose $F(b)=[\gamma']$.

Pour $k=2$, on peut choisir un repr{\'e}sentant $\gamma$
de $b$
tel que $\gamma =\lambda \volss$ sur le bout
parabolique en utilisant le lemme de Poincar{\'e}. La forme $\volss$ n'est pas a priori lisse sur
$\overline M'$. 
   Notons $\vol ^{S^2}$ la forme  lisse associ{\'e}e {\`a} $\overline
M'$. En int{\'e}grant sur une fibre o{\`u} 
 o{\`u} ces formes   sont toutes deux d{\'e}finies, on a par construction
$\int_{\CProj^1} \volss =\int_{\CProj^1}\vol^{S^2}$,
donc $\volss = \vol^{S^2}+ d\beta_2$, o{\`u} $\beta_2$ est  une
$1$-forme lisse sur $\bparab '$ sans aucune hypoth{\`e}se de
r{\'e}gularit{\'e} {\`a} l'infini.    
N{\'e}anmoins, en utilisant une fonction $\chi_j(t)$ d{\'e}finie au
\S\ref{subapproxmetraparab}) pour $j$  
suffisamment grand,   $\gamma' =
\gamma-\lambda d(\chi_j\beta )$ est lisse sur $\overline M'$, car
elle est {\'e}gale {\`a} $\lambda \vol^{S^2}$ pr{\`e}s de l'infini~;
on pose alors $F(b)=[\gamma']$.

On d{\'e}finit de m{\^e}me la r{\'e}ciproque de $F$ en utilisant cette
fois-ci le lemme de Poincar{\'e} local pour l'orbifold $\overline M'$.
On en d{\'e}duit le corollaire.
\end{demo}

\subsubsection{Cohomologie $L^2$ des surfaces r{\'e}gl{\'e}es}
Soit $M=\PP(\mathcal E)_{|\Sigma}$ une surface r{\'e}gl{\'e}e
obtenue {\`a} partir d'un fibr{\'e} parabolique $\mathcal E\rightarrow
\barSigma$ munie d'un mod{\`e}le  local de m{\'e}trique adapt{\'e} $\gmod$.
On consid{\`e}re maintenant la compactification 
 $\widehat M = \Proj (\mathcal E )_{|\overline \Sigma }$.

Pour un choix de poids
$\alpha'_1(P_i)=\alpha_2'(P_i)=0$, on a $\widehat M=\overline M'$, 
d'o{\`u} par le corollaire \ref{cohomologie} un isomorphisme
\begin{equation}
\label{isomcohom}
\HDR^*(\widehat M)\simeq \HLL^*(M).
\end{equation}
Or la cohomologie de la surface  r{\'e}gl{\'e}e $\widehat M$ est bien connue
(cf. \cite{Bv} prop. III.18)~; en particulier la $2$-cohomologie v{\'e}rifie~:
$$H^2(\widehat M,\ZZ)= \ZZ h\oplus \ZZ F,
$$
o{\`u} $F$ est le dual de Poincar{\'e} d'une fibre, et $h$ la classe du
fibr{\'e} tautologique $\mathcal O_{\widehat M}(1)$ sur $\widehat
M$. Alors $h^2= \deg \mathcal E$, $F^2=0$ et $h\cdot F=1$
; on en d{\'e}duit via l'isomorphisme (\ref{isomcohom}) la proposition
suivante~:  
 \begin{prop}
On a des isomorphismes
$$\HLL^0(M)\simeq \HLL^4(M)\simeq \R,\quad \HLL^1(M)\simeq  \HLL^3(M)\simeq H^1(\overline \Sigma,\R).
$$
   Il existe des classes de cohomologie $F,h\in \HLL^2(M)$ telles
   que $F^2=0$, $h^2=\deg \mathcal E$, $h\cdot F =1$ et 
$\HLL^2(M)= \R h\oplus \R F.$
En cons{\'e}quence $b_2=2$ et la forme d'intersection est de signature $(1,1)$.
\end{prop}

\subsubsection{D{\'e}monstration du lemme de Poincar{\'e} \ref{lemmePoincare}}
\medskip

\subparagraph{\emph{Cas des $1$-formes.}}
Sur le bout parabolique,  $\gamma $ se d{\'e}compose en
$\gamma = f dt + \gamma_2,$
o{\`u} $\gamma_2$ est une $1$-forme telle que $i_\dt\gamma_2=0$. 
Dans le cas o\`u $f$ est {\`a} support compact, d{\'e}finissons une
<<primitive>> de $f$ par
$\tilde f = \int_0^tf dt$ nulle en $t=0$.
D'apr{\`e}s le lemme \ref{courburemoyenne},
$$\int_{t=0} ^T |f|^2\vol^\gmod =\int_{t=0} ^T |\dt \tilde f|^2\vol^\gmod  \geq h_0^2 \int_{t=0}^T | \tilde
 f|^2 \vol^\gmod, $$
car le seul terme de bord non nul est positif.
On d{\'e}duit de cette derni{\`e}re in{\'e}galit{\'e} que pour $f\in L^2$, la
fonction $\tilde f$ est
{\'e}galement d{\'e}finie dans $ L^2$.
Or  $\gamma -d\tilde f$ v\'erifie par d{\'e}finition de $\tilde f$ 
$$\Lie_\dt
 (\gamma - d \tilde f)=(d\circ i_\dt +i_\dt\circ d)(  \gamma - d \tilde
 f) =0.$$
 Par cons{\'e}quent  $(\gamma - d \tilde f)$  est ind{\'e}pendante de $t$,
 et repr{\'e}sente une classe de cohomologie de $\HDR^1(S^1\times \CProj^1)=\R[d\theta]$~; il existe donc une fonction $f_2$ ind{\'e}pendante de
 $t$ (donc $L^2$), telle que $\gamma -d\tilde f -df_2 = \mu d\theta$,
avec  $\mu$ constante.

Soit $\Xtheta$ le champ de vecteurs d{\'e}fini par (\ref{exprXtheta}). Posons
$$V(t)= \int_{\{t\}\times N}\!\!\!\!\! \left(i_{e^t\Xtheta} \gamma\right)\vol^{\gmod_t}.$$
Comme $|e^t\Xtheta|=|e^{-t}d\theta|=1$ et $\gamma, \vol^{\gmod}\in L^2$, on en d{\'e}duit que   $\left(i_{e^t\Xtheta}
\gamma\right)\vol^{\gmod}$ est int{\'e}grable et on peut donc choisir des
tranches $t_k\rightarrow\infty$ telles que $V(t_k)\rightarrow 0$.
 
Par ailleurs 
\begin{equation*}
\begin{split}
 \int_N \Big(i_{e^t\Xtheta} &  d(\tilde f+f_2)
\Big )\vol^{\gmod_t} 
= \int_N \Big ( \Lie_{ e^t\Xtheta } (\tilde f+f_2)\Big ) \vol^{\gmod_t} \\
&= \int_N  \Lie_{e^t\Xtheta} \Big ( (\tilde f+f_2) \vol^{\gmod_t} \Big ) -
\int_N   (\tilde f+f_2) \Lie_{e^t\Xtheta} \vol^{\gmod_t}   \\
&= \int_N  d\circ i_{e^t\Xtheta} \Big ( (\tilde f+f_2) \vol^{\gmod_t} \Big ) -
\int_N   (\tilde f+f_2) \Lie_{e^t\Xtheta} \vol^{\gmod_t} =  0 ,
\end{split}
\end{equation*}
car la m{\'e}trique est invariante suivant $\Xtheta$ (cf.
(\ref{exprLieXtheta})). 
On en d{\'e}duit que 
$$V(t)=\int_N i_{e^t\Xtheta} \mu d\theta \vol^{\gmod_t} =\mu \int_N
e^t\vol^{\gmod_t}. $$ 
Or   $e^t\vol^{\gmod_t} =  d\theta\wedge \volss$ est ind{\'e}pendante de
$t$~; puisque $V(t_k)\rightarrow 0$, on a n{\'e}cessairement $\mu=0$. Finalement, on a d{\'e}montr{\'e} que 
$\gamma = d(\tilde f + f_2 ),$
o{\`u} $\tilde f$ et $f_2$ sont $L^2$ ce qui ach{\`e}ve la
d{\'e}monstration pour les $1$-formes. $\hfill \Box$
\medskip

\emph{Cas des $2$-formes.}
Soit $\gamma$ une $2$-forme ferm{\'e}e et $L^2$.
On la d{\'e}compose en $\gamma=dt\wedge \beta + \gamma_2$,
o{\`u} $\beta$ et $\gamma_2$ sont des formes diff{\'e}rentielles telles
que $i_\dt\gamma_2=0$ et  $i_\dt\beta=0$. 
Comme dans le cas des fonctions, nous cherchons une $1$-forme $\tilde
\beta$ dans $L^2$, telle que $i_\dt\tilde \beta = 0$ et
$\Lie_\dt\tilde \beta=\beta$. On a un la d{\'e}composition
$\beta=\mu d\theta + \beta_2,$
 o{\`u} $\mu$ est une fonction, et $\beta_2$  une $1$-forme 
telle que $i_\dt \beta _2  = i_\dtheta \beta_2= 0$. Comme
$\gamma$ est $L^2$ par hypoth{\`e}se, on en 
d{\'e}duit que $\beta_2$ et $\mu e ^t$ sont $L^2$.
Pour la partie  $\beta_2$, la technique est semblable au cas des
fonctions~; on  pose
$\tilde \beta_2 = \int_0^t \beta_2 dt$,
alors d'apr{\`e}s le lemme \ref{courburemoyenne}, toujours en
utilisant le fait  que le terme de bord est nul en $t=0$, on obtient
pour $\beta_2$ {\`a} support compact
$$  \int_0^\infty |\tilde \beta_2|^2\vol^\gmod  \leq \int_0^\infty |\dt\tilde
\beta_2|^2\vol=  \int_0^\infty
|\beta_2|^2\vol^ \gmod .$$
On en d{\'e}duit que $\tilde\beta_2$ est d{\'e}finie dans $L^2$ alors
$\beta_2\in L^2$. 

Dans le cas o{\`u} $\mu$ est {\`a} support compact, on d{\'e}finit $\tilde \mu =
\int_t^\infty \mu dt$~;
d'apr{\`e}s le lemme \ref{courburemoyenne}
\begin{equation*}
\int_0^\infty |\tilde \mu d\theta |^2\vol ^\gmod  = \int_0^\infty  |\tilde
\mu e^t |^2\vol^\gmod   
 \leq   \int_0^\infty | (\dt \tilde \mu) e^t |^2\vol ^\gmod  
= \int_0^\infty | \mu d\theta |^2\vol ^\gmod .
\end{equation*}
On en d{\'e}duit que
$\tilde \mu d\theta\in L^2$ si  $\mu d\theta\in L ^2$.
Posons maintenant $\tilde \beta = \tilde \mu d\theta +\tilde \beta_2
\in L^2$ 
et  $\gamma_3= \gamma - d\tilde\beta$. Alors  $i_\dt \gamma_3 = 0$
et $d\gamma_3 = 0$, d'o{\`u} $\Lie_\dt \gamma_3= 0$, ce qui implique que $\gamma_3$ est une
$2$-forme sur la tranche ind{\'e}pendante de $t$. On {\'e}crit
alors 
$\gamma_3=\beta_3\wedge d\theta + \gamma_4$,
avec $\beta_3$ et $\gamma_4$ des formes sur la tranche, ind{\'e}pendantes
de $t$ telles que $i_\Xtheta \beta_3 =0$ et $i_\Xtheta \gamma_4 = 0 $.

 On dira que $X$ est un  \emph{ champ de vecteurs horizontal} sur 
 $\bparab$ s'il est ind{\'e}pendant de $t$ et si $p^*X$ (o{\`u} $p : I_a \times
 \CProj^ 1\rightarrow \bparab$ est  le rev{\^e}tement universel
 riemannien associ{\'e} au bout parabolique) est  un champ de vecteurs
 tangent au facteur $\CProj^1 $  invariant sous l'action de $\tau$.   
 Si $X$ est un champ de vecteurs horizontal, 
$|X|_\gmod$ et  $|\Lie_X \volss|_\gmod$ sont uniform{\'e}ment
 born{\'e}s~; le volume {\'e}tant fini, ceci implique que $\Lie_X
\volss \in L^ 2$. 

Posons
$$ U(t)=\int_{\{t \}\times N } (i_X \gamma) \wedge\volss ,$$  
et 
$$V(t)=\int_{\{t\}\times N} (i_X d\tilde \beta)\wedge
\volss =\int _N \Lie_X\tilde \beta\wedge\volss =\int_N \tilde\beta \wedge
\Lie_X\volss .$$
Comme $\gamma \in L^2$ on a  $i_X\gamma \in L^2$~; de plus
$\Lie_X\volss$ et $\tilde\beta$ sont dans $L^2$. On peut donc choisir
des tranches $t_k\rightarrow \infty $  de sorte que $U(t_k)
\rightarrow 0$ et $V(t_k)\rightarrow 0$.  Mais  
$$U(t_k)- V(t_k) = \int_{\{t_k\}\times N}
i_X \gamma_3 \wedge \volss = \int_N
(i_X\beta_3) d\theta \wedge \volss $$ 
est ind{\'e}pendant de $k$. D{\`e}s que
$\beta _3 \neq 0$, on peut choisir $X$
tel que cette int{\'e}grale soit non nulle. N{\'e}cessairement on a alors
$\beta_3 = 0$. On en d{\'e}duit que  $\gamma_3=\gamma_4$ est une
$2$-forme sur la 
tranche ind{\'e}pendante de $t$ telle que $i_\Xtheta  \gamma_3 = 0$, $d\gamma_3=0$,
d' o{\`u} $\Lie_\Xtheta \gamma_3= 0$. 
On en d{\'e}duit
en passant au rev{\^e}tement que $p^* \gamma_3$ est une $2$-forme
$\ZZ$-invariante sur le facteur 
$\CProj ^ 1$. 

En utilisant la cohomologie de De Rham
usuelle, on peut {\'e}crire $\gamma_3 = \lambda \volss + d\beta_4$, avec $\beta_4$
une $1$-forme sur $\CProj ^1$.
Si $\alpha_2-\alpha_1\in \Rat$, l'action de $\ZZ$ induite sur $\CProj^1$ se
factorise en une action de $\ZZ_q$.  La $1$-forme $\beta_4$ n'est pas 
n{\'e}cessairement invariante, mais quitte {\`a} moyenner
sous l'action de $\ZZ_q$, on peut la supposer invariante de sorte que
$\beta_4 $ 
passe au quotient et $\beta_4 \in L^2$. On a donc d{\'e}montr{\'e} que $\gamma =
\lambda \volss + d(\tilde \beta + \beta_4)$ avec $ \tilde \beta$ et
$\beta_4 \in L^2 $.

Dans le cas o{\`u} $\alpha_2-\alpha_1\notin \Rat$, l'action de $\ZZ$ sur $\CProj
^1$ ne se factorise pas. Cette action est  engendr{\'e}e par
une rotation $\mathcal R_\alpha$ donn{\'e}e par 
$ \tau \cdot [\tilde u : \tilde v] =  [\tilde u e^{2i\pi\alpha_1} : \tilde v
e^{2i\pi\alpha_2}]$,
d'angle $-2\pi\alpha$. Comme $\alpha\notin\Rat$, le groupe engendr{\'e} par
$\mathcal R_\alpha$ est dense dans le groupe des rotations de m{\^e}me
axe. Par densit{\'e}, $\gamma_4$ est invariante par toute rotation
$\mathcal R_\phi$ de ce type.
Puisque la forme volume $\volss$ est invariante par rotation, la
forme exacte $\gamma_3 - \lambda \volss $ est {\'e}galement invariante
par rotation. On en d{\'e}duit comme dans le cas rationnel, en moyennant
sur l'action de $\mathcal R_\phi$ pour $\phi\in S^1$, une
$1-$forme $\beta_4$ invariante par rotation telle que $\gamma_3 -
\lambda \volss =d\beta_4$. En particulier, $\beta_4$ est invariante
par $\tau$ et passe au quotient. $\hfill \Box$
\medskip

\emph{Cas des $3$-formes.}
Une  $3$-forme sur le bout parabolique  se d{\'e}compose en
$ \Omega  = dt\wedge (\beta\wedge d\theta +\lambda \volss ) + \Omega_2$,
o{\`u}  $\beta$ et $\Omega_2$ v{\'e}rifient
$i_\dt \beta =i_\dtheta \beta = 0$ et $i_\dt \Omega_2 = 0$.
On cherche des <<primitives>> de $\beta$ et $\lambda$ en posant
$\tilde \beta  = \int_t^\infty \beta dt$ et $\tilde \lambda =
\int_0^t \lambda dt$ dans le cas de formes {\`a} support compact.
D'apr{\`e}s le lemme \ref{courburemoyenne}, 
$$
\int_0^\infty |\tilde \beta\wedge  d\theta |^2\vol ^\gmod  
 \leq  \frac 1{h_0^2} \int_0^\infty | (\dt \tilde \beta ) e^t |^2\vol ^\gmod  
= \frac 1{h_0^2}\int_0^\infty | \beta\wedge  d\theta |^2\vol ^\gmod ,
$$
$$  \int_0^\infty |\tilde \lambda \volss |^2\vol^\gmod=  \int_0^\infty |\tilde
\lambda  |^2 \vol^\gmod  \leq \frac 1{h_0^2} \int_0^\infty |\dt\tilde
\lambda |^2\vol=   \frac 1{h_0^2} \int_0^\infty
|\lambda |^2\vol^ \gmod ,$$
ce qui donne un sens  {\`a} $\tilde\beta\wedge d\theta$ et $\tilde
\lambda$ dans $L^2$
pour $\lambda$ et $\beta\wedge d\theta \in L^2$.
On en d{\'e}duit une $2$-forme $L^2$ d{\'e}finie par $\omega = \tilde \beta
\wedge d \theta + \tilde \lambda\volss $,
telle que $\Omega_3=\Omega -d\omega$ v{\'e}rifie $i_\dt\Omega_3 =0$, par d{\'e}finition de
$\omega$. Alors 
$\Lie _\dt \Omega_3 = (d\circ i_\dt + i_\dt \circ d) \Omega_3 =0$,
 donc $\Omega_3$ est une $3$-forme  ind{\'e}pendante de $t$ sur la
tranche $N\simeq S^1\times \CProj ^ 1$ du bout parabolique. On
peut donc l'{\'e}crire 
$\Omega_3 = \mu d \theta \wedge \volss,$ o{\`u} $\mu $ est une
fonction ind{\'e}pendante de $t$. 
 
S'il existe une fonction $f$ constante suivant $\dt$ et $\Xtheta$, telle que 
$\int_N f \mu d\theta\wedge \volss \neq 0$,
on a alors une contradiction~: $f$, $|df|$ sont born{\'e}es car
$f$ ne d{\'e}pend que du facteur $\CProj^1$. On en d{\'e}duit que $f$ et $df$ sont $L^ 2$
puisque nous sommes en volume fini. Maintenant 
$$
\int_{\{t\}\times N} f\Omega = \int_{\{t\}\times N} f\Omega_3 +
\int_{\{t\}\times N} fd\omega 
 = \int_{\{t\}\times N} f\Omega_3 -
\int_{\{t\}\times N} df\wedge\omega,
$$
et nous pouvons  choisir
des tranches $t_k$ telles que  $\int_{t_k} f\Omega$ et $\int_{t_k} df\wedge
\omega$  tendent vers $0$.
Par ailleurs  $\int_{\{t\}\times N} f\Omega_3 $ est ind{\'e}pendante de
$t$ d'o{\`u} une contradiction.
Par cons{\'e}quent $\mu$ est orthogonale aux fonctions  constantes
suivant $\Xtheta$ ce qui entra{\^\i}ne l'existence de $\tilde
\mu$ telle que $\Xtheta\cdot\tilde \mu=
\mu$. Comme $\tilde \mu \volss $ est born{\'e}e, donc $ L^2$, on en
d{\'e}duit que  $\Omega = d(\omega+ \tilde\mu \volss)$
est exacte au sens de la cohomologie $L^2$ sur le bout
parabolique. $\hfill \Box$\medskip      

\emph{Cas des $4$-formes.} Une $4$-forme $L^2$ sur le bout
parabolique s'{\'e}crit 
$\gamma = \mu dt\wedge d\theta \wedge \volss$,
avec $\mu d\theta \in L^2$. Nous avons d{\'e}j{\`a} d{\'e}montr{\'e}
(cf. cas des $2$-formes) qu'il existe une
fonction $\tilde \mu$ d{\'e}finie sur le bout parabolique telle que
$\tilde \mu d\theta \in L^ 2$  et $\dt\tilde \mu = \mu$. On en d{\'e}duit
$d\tilde \gamma =\gamma$ avec $\tilde \gamma = \tilde \mu d
\theta\wedge \volss \in L^ 2$, d'o{\`u} le r{\'e}sultat. $\hfill\Box$

\subsection{Convergence des formes harmoniques}
\label{sub:convharm}
Dans tout le \S\ref{sub:convharm},
on suppose que les poids des bouts paraboliques de $M$ v{\'e}rifient
la condition  
$\alpha=\alpha_2-\alpha_1 \in\Rat$ par souci de clart{\'e}. Puis  nous  reviendrons plus tard sur le cas irrationnel
(cf. \S\ref{secirrat})  pour
lequel on a des r{\'e}sultats analogues.

\subsubsection{In{\'e}galit{\'e}s de Poincar{\'e} pour les fonctions}
Les lemmes~\ref{volume} et \ref{courburemoyenne}  nous
permettent  facilement 
d'obtenir une borne inf\'erieure uniforme sur la premi{\`e}re
valeur propre du laplacien des m\'etriques $g_j$ pour les fonctions~(cf. \cite{Bq3} corollaire 4.3).
\begin{lemme}[In{\'e}galit{\'e} de Poincar{\'e}]
\label{laplacien}Supposons que $M$ poss{\`e}de des bouts paraboliques
v{\'e}rifiant la condition $\alpha=\alpha_2-\alpha_1\in \Rat$. Soit
$g$ une m{\'e}trique asymptotique au mod{\`e}le local $\gmod$  et soit 
$g_j$  la suite d'approximations de $g$ lisses sur
$\overline M$.
Il existe une constante $c>0$, ind{\'e}pendante de j, telle que pour
  toute fonction $f$ sur $\overline M$ v{\'e}rifiant $\int_M f\volj =0$, on ait
$$
    \int_M|df|^2 \volj\geq c\int_M|f|^2\volj,
$$
les normes {\'e}tant prises par rapport {\`a} la m{\'e}trique
$g_j$. 
\end{lemme}
On en d\'eduit l'in\'egalit\'e de Poincar\'e \`a la limite~:
\begin{cor}[In{\'e}galit{\'e} de Poincar{\'e}]
Soit $g$ une m{\'e}trique \aparab sur $M$. Il existe une constante $c>0$
telle que pour toute fonction $f\in L^2_1(g)$ v{\'e}rifiant $\int_M f\vol^g
= 0$, on ait
$$\int_M|df|^2\vol^g \geq c \int_M |f|^2\vol^g.
$$  
\end{cor}\medskip
\begin{demo}
  On applique le
corollaire \ref{laplacien}  {\`a} $h_j=(1-\chi_j)(f+c_j)$, o{\`u} $c_j$
est une constante choisie de sorte que $\int h_j\volj = 0$. On en d{\'e}duit le r{\'e}sultat en passant {\`a} la limite.
Comme nous le verrons au \S\ref{secirrat}, cette d{\'e}monstration,
dans laquelle nous avons fait implicitement l'hypoth{\`e}se
$\alpha_2-\alpha_1\in \Rat$, s'{\'e}tend m{\^e}me si cette condition n'est
pas v{\'e}rifi{\'e}e.
\end{demo}

\remarque
Cette in{\'e}galit{\'e}  de Poincar{\'e} permet classiquement de d{\'e}velopper  la
th{\'e}orie de Hodge des $1$-formes $g$-harmoniques $L^2$. On en d{\'e}duit que
toute classe de cohomologie de $\HLL^1(M)$ admet un unique
repr{\'e}sentant $g$-harmonique $L^2$.

\subsubsection{Convergence des $1$-formes harmoniques}
L'in{\'e}galit{\'e} de Poincar{\'e} du lemme~\ref{laplacien} sur les fonctions
permet maintenant de d{\'e}montrer le r{\'e}sultat suivant.
\begin{prop}
\label{convcompact}
Soit une classe de cohomologie $b\in \HDR^1(\overline M)$. La
 suite $\gamma_j$ 
 de ses   repr{\'e}sentants $g_j$-harmoniques  sur $\overline
M$  converge au sens 
$C^\infty$ sur tout compact de $M$ vers le repr{\'e}sentant 
$g$-harmonique  $L^2$ de $b$. 

Plus pr{\'e}cis{\'e}ment, il existe une constante $c>0$, telle que 
si l'on choisit un repr{\'e}sentant {\`a} support compact $\beta$ de $b$, et
qu'on {\'e}crit $\gamma_j = \beta + df_j$, o{\`u} $f_j$ est une fonction sur $\overline
M$ v{\'e}rifiant $\int_M f_j\volj = 0$, alors $f_j$ converge  sur tout compact de $M$ vers une
fonction $f\in L^2_1(g)$ telle que $\gamma= \beta + df $ soit $g$-harmonique, 
$$\int_Mf\vol^g = 0,\mbox{ et } \| f\|_{L^2_1} \leq c\|\beta \NLLg. $$
\end{prop}
\begin{demo}
Le repr{\'e}sentant $g_j$-harmonique de $b$ s'{\'e}crit sous la forme 
$\gamma_j=\beta + d f_j$, 
 o{\`u} $f_j$ v{\'e}rifie
$$\int_M f_j\volj = 0,  \quad \Delta^{g_j}f_j= -d^{*_{g_j}}\beta ;
$$  
d'apr{\`e}s le lemme \ref{laplacien},
$$c\int |f_j|^2\volj \leq \int |df_j|^2\volj= -\int f_j d^{*_j}\beta \volj =
-\int \langle df_j, \beta\rangle \volj, $$
et d'apr{\`e}s l'in{\'e}galit{\'e} de Cauchy-Schwarz, on en d{\'e}duit un contr{\^o}le 
$$\sqrt c \, \|f_j\NLLj
\leq \| df_j \NLLj \leq \|\beta \NLLj.$$
Comme $\beta$ est {\`a} support compact, il est clair que $\|\beta \NLLj =
\|\beta \NLLg$ pour $j$ suffisamment grand. 
Sur un compact $K$ de $M$
la norme $L^2_1$ de $f_j$ est donc uniform{\'e}ment
born{\'e}e.
Quitte {\`a} extraire une sous suite, on peut supposer  que $f_j$ converge sur tout compact de $M$ vers une limite
faible $f\in L^2_1(g)$ qui v{\'e}rifie 
$$ \sqrt c \, \|f \NLLg \leq \|\beta\NLLg ,\; \| df \NLLg \leq
\|\beta \NLLg, 
\; \Delta^{g}f= -d^{*_{g}}\beta , \; \int _M f\vol^g= 0.$$
Seule la derni{\`e}re identit{\'e} n'est pas {\'e}vidente~: on sait que
$f\in L^2 $ et que $\|f_j\NLLj$ uniform{\'e}ment born{\'e}e~; d'apr{\`e}s le lemme
\ref{volume},  on peut donc choisir un compact
$K$ de $M$ suffisamment grand, tel que 
$$ \int_{M\setminus K} |f| \vol^g  \leq \vol^g(M\setminus K)^{\frac 12
}\left
(\int_{M\setminus K} |f|^2\vol^g \right )^{\frac 12} \leq\epsilon,$$
$$ \left | \int_{M\setminus K} f_j \volj \right |  \leq \volj
(M\setminus K)^{\frac 12 } \; \| f_j \NLLj \leq \epsilon \mbox{ pour
tout $j$.}$$
Par compacit{\'e} de l'inclusion $L^2_1\hookrightarrow L ^2$,  $(f_j)$
converge vers $f$ au sens $L^2$-fort sur $K$.
Alors 
$$\int_K f\vol^g  = \lim \int_K f_j \volj = -  \lim \int_{M\setminus K} f_j \volj,
$$
d'o{\`u}  $ |\int_K f\vol^g| \leq \epsilon $. 
Maintenant, 
$$\left | \int_M f\vol^g \right | \leq   
\left | \int_K f\vol^g \right | +
\int_{M\setminus K} |f| \vol^g  \leq 2\epsilon ,
$$
d'o{\`u} $\int_M f\vol^g = 0$.

Par r{\'e}gularit{\'e} elliptique la
convergence est en r{\'e}alit{\'e} $C^\infty$ sur tout compact de $M$.
Finalement $\gamma = \beta + df$ est une $1$-forme $g$-harmonique
$L^2$ qui repr{\'e}sente la classe de cohomologie d{\'e}finie par $b$ en
cohomologie $L^2$ d'o{\`u} le r{\'e}sultat.
\end{demo}

\emph{Comportement {\`a} l'infini.}
Nous avons {\'e}galement un contr{\^o}le pr{\`e}s de l'infini sur les $\gamma_j$
de la proposition~\ref{convcompact}
{\`a} l'aide du lemme suivant~:
\begin{lemme}
\label{lemmecompinfty}
En reprenant les notations de la proposition \ref{convcompact},  pour tout $\epsilon>0$, il
existe un $T$ suffisamment grand tel que
$$ \int_{t\geq T}
 |\gamma|^2 \vol^g \leq \epsilon \quad \mbox{ et } \quad \int_{t\geq T}
 |\gamma_j|^2 \vol^{g_j} \leq \epsilon .$$
\end{lemme}
\begin{demo}
  Reprenons  la d{\'e}monstration de la
  proposition \ref{convcompact} ainsi que ses notations.
Par hypoth{\`e}se, $\beta$ est {\`a} support compact, donc pour $T$ assez grand, on a
$\gamma_j=\beta + df_j =df_j$ d'o{\`u}
$$\int_{t\geq T} |\gamma_j |^2\volj = \int_{t\geq T} |df_j |^ 2 \volj .$$
 En int{\'e}grant par
parties, il vient 
\begin{eqnarray*}
\int_{t\geq T} |df_j|^2 \volj  & = & \int_{t\geq T} df_j\wedge *_j df_j 
 = \int_{t\geq T} d(f_j  *_j df_j) -  \int_{t\geq T} f_j d *_j df_j \\
&= & \int_{t= T} f_j *_j df_j + \int_{t\geq T} \langle \Delta
^{g_j}f_j,f_j\rangle \volj  .
\end{eqnarray*}
Comme nous nous sommes plac{\'e} hors du support de $\beta$, 
$ \Delta ^{g_j}f_j = \Delta
^{g_j} f_j+ d^{*_j} \beta  =\Delta
^{g_j}\gamma_j =0$,
d'o{\`u}
$$\int_{t\geq T} |df_j|^2 \volj = \int_{t= T} f_j *_j df_j .$$
Le m{\^e}me calcul reste valable pour la limite $\gamma = \beta +df$ car $f\in
L_1^2$. 
Choisissons un $T$ suffisamment grand, de sorte que $\int _{t\geq
  T}|df|^2\vol^g  \leq \epsilon /2 $. Alors comme $f_j$
tend vers  $f$ au sens $C^\infty$ sur tout compact on a 
$$  \int_{t= T} f_j *_j df_j \rightarrow   \int_{t= T} f * df \mbox{ pour } j
\rightarrow\infty $$
soit
$$\int_{t\geq T} |df_j|^2 \volj \rightarrow \int_{t\geq T} |df|^2 \volg \mbox{ pour } j
\rightarrow\infty  $$
dont on d{\'e}duit le lemme.
\end{demo}

\subsubsection{In{\'e}galit{\'e} de Poincar{\'e} pour les $1$-formes}
Afin de d{\'e}montrer un r{\'e}sultat de convergence des $2$-formes
harmoniques (cf. prop. \ref{propconvharm}) analogue {\`a} celui des
$1$-formes, nous commen{\c c}ons par d{\'e}montrer la proposition suivante~:
\begin{prop}
\label{laplacien2}
Supposons que $M$ soit une vari{\'e}t{\'e} avec des bouts paraboliques
  v{\'e}rifiant la condition $\alpha\in \Rat$. Soit $g$ une m{\'e}trique
  \aparab et $g_j$ la suite d'approximations associ{\'e}e.
  Il existe une constante $c>0$ telle que pour toutes les m{\'e}triques
  $g_j$, et toute $1$-forme orthogonale aux formes $g_j$-harmoniques, on ait~:
$$\int_M(|d\beta|^2+|d^{*_j}\beta|^2)\volj \geq c\int_M|\beta|^2\volj,
$$
les normes {\'e}tant prises par rapport {\`a} la m{\'e}trique
$g_j$. L'in{\'e}galit{\'e}  est {\'e}galement v{\'e}rifi{\'e}e pour la m{\'e}trique $g$ et
toute $1$-forme $\beta \in L^2_1(g)$ orthogonale aux formes
$g$-harmoniques $L^2$.
\end{prop}
\remarque Cette in{\'e}galit{\'e} permet de d{\'e}velopper la th{\'e}orie de
Hodge $L^2$ pour la m{\'e}trique $g$ dans le cas des $2$-formes. Ainsi,
toute classe de $\HLL^2(M)$ admet un unique repr{\'e}sentant
$g$-harmonique  $L^2$.

Nous avons besoin d'un lemme pr{\'e}liminaire avant de commence la
d{\'e}monstration.
\begin{lemme}
\label{courburemoyenne2}
Soit $\gmod$  un mod{\`e}le local de m{\'e}trique dont les bouts
    v{\'e}rifient la condition de rationalit{\'e}.
    Il existe une constante $T$ suffisamment grande telle que pour toute  m{\'e}trique
  $g=\gmod$ ou $\gmodj$, et toute $1$-forme $\beta$ sur
  $[t_1,t_2]\times N$ avec $t_1,t_2\geq T$, on ait~:
$$ h_0 ^ 2 \int_{[t_1,t_2]\times 
  N}\!\!\!\!|\beta|^2  \;\vol + h_0 \int_{\partial [t_1,t_2]\times N}
\!\!\!\!\!\!\!\!\!\!\!\!\!\!\!\!\!\!|\beta |^2  \;\vol ^N \leq  \int_{[t_1,t_2]\times
N}\left (|\nabla \beta|^2 +  \Ric(\beta,\beta) \right )\;\vol,
$$
les normes {\'e}tant prises par rapport {\`a} la m{\'e}trique $g$. 
\end{lemme}
\begin{demo}
Soit $g$ une m{\'e}trique  de la forme $g=dt  ^ 2 + g_t$, sur
$[t_1,t_2]\times N$. On d{\'e}compose les $1$-formes en $\beta= f dt +
\mu$, avec $i_\dt \mu =0$. Alors la connexion de Levi-Civita se
d{\'e}compose en 
$\nabla \beta = dt\otimes \nabla _\dt \beta + \nabla_{|N} \beta,$
o{\`u} $\nabla_{|N}\beta $ est la restriction de $\nabla\beta $ {\`a} une
section de $T^*N\otimes \Omega ^ 1([t_1,t_2]\times N)$. Calculons ce terme plus
pr{\'e}cis{\'e}ment~: notons $\nabla ^N$ la connexion de Levi-Civita induite
sur $N$ par $g$ et $d^N$ la diff{\'e}rentielle suivant la tranche, alors
$$\nabla_{|N} \beta = d ^Nf\otimes dt + f \nabla_{|N}dt + \nabla ^N\mu
+ \ffond (\mu)\otimes dt,
$$
or $ \nabla_{|N} dt = \langle ., \nabla_{|N}\dt\rangle = - \langle
\nabla_{|N} .,\dt \rangle= -\ffond$,
donc 
$$\nabla \beta = dt\otimes \nabla _\dt \beta + ( d ^Nf+ \ffond (\mu)
)\otimes dt  - f \ffond  + \nabla ^N\mu .$$
Appliquons ce calcul dans le cas o{\`u} $g$ est {\'e}gale au mod{\`e}le local
$\gmod$, o{\`u} {\`a} $\gmodj$.
La m{\'e}trique  $g_t$ sur la tranche $N\simeq S^ 1\times \CProj
^1$ s'exprime  en prenant par exemple la coordonn{\'e}e $u$ de projection
st{\'e}r{\'e}ographique d{\'e}finie par $v=1$,
\begin{equation}
 g_t= \phi ^2(t)  d\theta^2 + \frac{4/c}{(1+|u|^2)^2}\left
|du-i\alpha u d\theta\right |^2 ~;
\end{equation} 
la deuxi{\`e}me forme fondamentale est donn{\'e}e par 
$$\ffond = -\frac 12 g_t = - \frac{\dt \phi}\phi (\phi d\theta) ^2
\quad\mbox { et }\quad  h= \frac 12 \trace_g \; \ffond = - \frac 12 \frac
{\dt\vol^{g_t}}{\vol^{g_t}} = - \frac 12 \frac{\dt \phi}\phi.
$$
En notant $\mu = f_2 \phi d\theta  + \eta $ avec
$i_\Xtheta\eta=0$ et $d^{\CP} f$ la partie de $df$ orthogonale {\`a}
$d\theta$ et $dt$, la connexion de Levi-Civita sur la tranche est
{\'e}gale {\`a} 
$$ \nabla ^N \mu = (\phi ^ {-1} \Xtheta \cdot f_2) \phi d\theta \otimes
\phi d\theta + (d^{\CProj^1}f_2)  \otimes \phi d\theta + \phi
d\theta\otimes (\Lie_{\phi ^ {-1} \Xtheta} \eta) + \nabla
^{\CProj ^1} \eta .$$
Finalement  on peut {\'e}crire 
\begin{multline*}
\nabla \beta = dt\otimes \nabla _\dt \beta +  (\phi ^{-1}\Xtheta\cdot f + f_2 )
\phi d\theta
\otimes dt +(\phi ^{-1}\Xtheta \cdot f_2 - f) \phi d\theta\otimes \phi d\theta \\
+ (d^{\CProj ^1}f) \otimes dt+ (d^{\CProj ^1}f_2) \otimes \phi d\theta
+ \phi d\theta\otimes (\Lie_{\phi^{-1}\Xtheta} \eta)+ \nabla ^{\CProj
  ^1}\eta , 
\end{multline*}
or 
$\Ric^g(\beta,\beta)= \Ric ^{\CProj ^1}(\eta, \eta) -  ({\dt ^ 2
  \phi}/\phi) ( |f|^2 + |f_2|^2),$
 o{\`u}  $\Ric^{\CProj^1}$ est la
  courbure  provenant de la m{\'e}trique de Fubini-Study~; on en
d{\'e}duit l'in{\'e}galit{\'e} suivante en faisant une int{\'e}gration par
parties sur la tranche  avec le terme
  en $|\nabla ^{\CP}\eta|^2$ et en en utilisant le fait que par
construction  $-  {\phi ^{-1}\partial_t ^ 2
  \phi} \geq  -1$. 
\begin{multline}
\label{exprintermpoinc1}
\int_N |\nabla \beta |^2 + \Ric^g(\beta,\beta) \geq \int_N |\nabla_\dt \beta |^ 2 +
|\phi ^{-1} \Xtheta \cdot  f|^2 +  |\phi ^{-1}\Xtheta \cdot  f_2|^2 
 \\
 +2 \langle \phi ^{-1} \Xtheta
\cdot  f,f_2  \rangle - 2 \langle f , \phi ^{-1} \Xtheta
\cdot  f_2 \rangle + |d^{\CProj ^1}f|^2+ |d^{\CProj ^1}f_2|^2.
\end{multline}

Nous allons  maintenant montrer que les termes mixtes de
(\ref{exprintermpoinc1}) sont contr{\^o}l{\'e}s.
En int{\'e}grant par parties sur $N$ il vient
$$  \int_N \langle  \phi ^{-1}\Xtheta
\cdot  f,f_2  \rangle \vol ^g = - \int_N \langle f , \phi ^{-1} \Xtheta
\cdot f_2 \rangle\vol ^g\; ; $$
d{\`e}s que $f$ ou $f_2$ est
constante suivant $\Xtheta$, ces int{\'e}grales sont nulles et on conclut  {\`a}
l'aide du lemme \ref{courburemoyenne}.
On supposera {\`a} partir de maintenant que $f$ et $f_2$ sont
orthogonales aux fonctions constantes suivant $\Xtheta$.
Dans le cas d'un bout parabolique
dont les poids v{\'e}rifient la condition de rationnalit{\'e}, le
champ de vecteurs $\Xtheta$ est p{\'e}riodique et donne une structure
de $S^1$-fibr{\'e}. Dans le rev{\^e}tement {\`a} $q$
feuillets $\pq :\Delta ^*_{a/q}\times \CProj ^1\rightarrow
\bparab$, les orbites de $\Xtheta$ correspondent aux cercles
centr{\'e}s en $0$ de
$\Delta ^*_{a/q}$. 

En utilisant la d{\'e}composition en s{\'e}rie
de Fourier sur ce rev{\^e}tement local, $f=\sum_k c_k e ^{k i \theta/q}$, on voit que la condition
d'{\^e}tre orthogonale aux constantes se traduit par $c_0=0$. Alors 
$$ \int_{N} |\Xtheta \cdot  f | ^2 =   \int _{N} \sum_k |\frac kq
c_k|^2 \geq  \frac 1 {q ^ 2}\int_{N}|f|^2, $$
d'o{\`u}
$$  \int_N |\phi ^ {-1}\Xtheta \cdot f| ^2 \geq  \frac 1
{(q \phi ) ^ 2} \int_N | f | ^2 , \quad  \int_N |\phi ^ {-1}\Xtheta
\cdot f_2| ^2 \geq  \frac 1
{(q \phi ) ^ 2} \int_N | f_2 | ^2 . $$
Puisque $\phi(t)\leq e ^{-t}\rightarrow 0$, on en d{\'e}duit que pour
$t$
suffisamment grand, les termes mixtes de (\ref{exprintermpoinc1}) sont
tr{\`e}s bien contr{\^o}l{\'e}s par les termes en $|\phi
^{-1}\Xtheta\cdot  f| ^ 2$
et  $|\phi ^{-1}\Xtheta\cdot f_2|^ 2$, d'o{\`u}
$$\int_N |\nabla \beta |^2 + \Ric^g(\beta,\beta) \geq \int_N
|\nabla_\dt \beta |^ 2 
$$
et le  lemme \ref{courburemoyenne} nous permet de conclure.
\end{demo}

\remarque  les hypoth{\`e}ses que nous devons faire dans ce lemme  sur la m{\'e}trique
sont assez faibles~; par exemple, la courbure positive du facteur $\CP$
n'intervient pas. 
On peut ainsi ais{\'e}ment corriger au passage la d{\'e}monstration de
l'in{\'e}galit{\'e} de Poincar{\'e} sur les $1$-formes de~\cite{Bq3}, o{\`u}
une erreur s'est gliss{\'e}e  dans les termes mixtes de
la formule de Bochner du lemme 4.5.

Ce lemme reste vrai (en modifiant l{\'e}g{\`e}rement les
constantes) si l'on remplace $\gmod$ et $\gmodj$
par une m{\'e}trique \aparab $g$ et ses approximations $g_j$.
\begin{cor}
\label{corcontrpoinc1}
Soit $g$ une m{\'e}trique \aparab $\gmod$.
    Il existe des constantes $T,c>0$ telles que pour toute
    m{\'e}trique
  $g$ ou $g_j$, et toute $1$-forme $\beta$ sur
  $[t_1,t_2]\times N$ avec $t_1,t_2\geq T$, on ait~:
$$ h_0 ^ 2 \int_{[t_1,t_2]\times 
  N}\!\!\!\!\!\!\!\!|\beta|^2  \;\vol  + h_0 \int_{\partial [t_1,t_2]\times N}
\!\!\!\!\!\!\!\!|\beta |^2  \;\vol ^N  \leq  c \int_{[t_1,t_2]\times
N}\left (|\nabla\beta|^2 + \Ric(\beta,\beta)\right ) \;\vol .
$$
\end{cor}
\begin{demo}
Comme les m\'etriques $\gmodj$ et $g_j$ sont \'egales pour $t\geq
j+\frac 12$ il nous suffit de v\'erifier le corollaire pour la
m\'etrique $g$. Il faut voir que l'in{\'e}galit{\'e} est
perturb{\'e}e lorsque on passe de $\gmod$ \`a $g$ par des termes  en
  $\epsilon(T)O(\|\beta\|^2_{L^2_1})$, avec $\epsilon(T)\rightarrow 0$
  quand $T\rightarrow\infty$. Ce n'est pas compl{\`e}tement {\'e}vident pour
 les termes de bord~; pour la m{\'e}trique $\gmod$
$$\int_{t'} |\beta| ^2 \vol^\gmodt = \int_{t\leq t'} \langle \beta,
\nabla_\dt \beta \rangle -2h |\beta| ^2 \; \vol ^\gmod
$$
et on en d\'eduit le corollaire.
\end{demo}

Nous avons maintenant fait un premier pas dans la d{\'e}monstration de la
proposition \ref{laplacien2} gr{\^a}ce au corollaire suivant~:
\begin{cor}
   \label{corcontrloc2}
 Il existe une constante $c>0$ et $T$ suffisamment grand, tels que pour toute
 $1$-forme lisse $\beta$ sur le bout parabolique compactifi{\'e} $\overline
 \bparab$ nulle en $t=T$ on ait
 \begin{equation}
   \int_{t\geq T} \left ( |d\beta|^2+|d^{*_j}\beta|^2 \right ) \volj \geq c
   \int_{t\geq T}|\beta |^2 \volj ,
 \end{equation}
les normes {\'e}tant prises par rapport {\`a} la m{\'e}trique $g_j$. Cette
in{\'e}galit{\'e} est {\'e}galement valable pour la m{\'e}trique $g$, et pour  une
$1$-forme $\beta\in L^2_1(g)$ d{\'e}finie sur $\bparab$ et nulle en $t=T$.
\end{cor}
\begin{demo}
  On applique le corollaire \ref{corcontrpoinc1} sur
$[T,T_j]\times N$ pour la m{\'e}trique $g_j$. Le terme de bord est nul par hypoth{\`e}se pour
$t=T$. Il est nul {\'e}galement pour $t=T_j$ car le
volume de cette tranche est {\'e}gal {\`a} $0$. 
Dans le cas de la m{\'e}trique
$g$ on fait tendre $T_j$ vers
l'infini et le terme de bord  correspondant tend vers $0$ car
$\beta\in L^2_1$. On en d\'eduit le corollaire par la formule de
Bochner $\Delta = \nabla ^*\nabla +\Ric$ en int\'egrant par parties.
\end{demo}

\begin{demode}{de la proposition \ref{laplacien2} }
  supposons la proposition fausse~; on en d{\'e}duit une suite de
  $1$-formes $\beta_j $ est orthogonales aux
  formes $g_j$-har\-mo\-ni\-ques telles que 
$\|\beta_j\NLLj =1$, $\|d\beta_j\NLLj \rightarrow 0$ et
  $\|d^{*_j} \beta_j\NLLj \rightarrow 0$.
Puisque $\|\beta_j\NLLj$ est born{\'e}e on peut en  extraire  
une limite faible  $\beta\in L_1^2(g)$ sur tout compact  de $M$.
La limite v{\'e}rifie n{\'e}cessairement $\|d\beta_j\NLLg = \|d^{*_g}
  \beta \NLLg = 0$ et $\|\beta\NLLg = 1$~; $\beta$ est donc $g$-harmonique. 

Montrons maintenant que $\beta$ est orthogonale aux formes
$g$-harmoniques $L^2$.
Soit $\gamma$ une $1$-forme $g$-harmonique
$L^2$ et $\gamma_j$ les repr{\'e}sentant $g_j$ harmoniques de $[\gamma]$.
On d{\'e}compose l'int{\'e}grale 
\begin{equation}
\label{exprinterm1}
0= \int_M\langle \beta_j, \gamma_j\rangle\volj =  \int_{M\setminus \{t\geq
  T\} } \langle\beta_j,\gamma_j\rangle\volj + \int_{t\geq
  T}\langle\beta_j,\gamma_j\rangle\volj
\end{equation}
et on majore  
$$ \left | \int_{t\geq T}\langle\beta_j,\gamma_j\rangle\volj \right |^2\leq
\|\beta_j\NLLj ^2  \int_{t\geq T} |\gamma_j|^2\volj  =   \int_{t\geq T} |\gamma_j|^2\volj .
$$
De m{\^e}me,
$$ \left | \int_{t\geq T}\langle\beta,\gamma\rangle\vol^g \right |^2\leq
\|\beta\NLLg ^2  \int_{t\geq T} |\gamma |^2\vol^g  =   \int_{t\geq T} |\gamma |^2\vol^g .
$$
D'apr{\`e}s le lemme
 \ref{lemmecompinfty} choisir un $T$ suffisamment
 grand tel que 
$$ \left | \int_{t\geq T}\langle\beta_j,\gamma_j\rangle\volj \right
|\leq \epsilon ,
\quad \left | \int_{t\geq T}\langle\beta,\gamma\rangle\vol^g \right
| \leq {\epsilon}. 
$$
D'apr{\`e}s (\ref{exprinterm1}), on a 
$\left | \int_{M\setminus \{t\geq T\}} \langle\beta_j,\gamma_j\rangle\volj
\right | \leq \epsilon$.
Sur un compact $\gamma_j \rightarrow \gamma $ et $\beta_j
\rightarrow \beta $ au 
sens $L^2$-fort par compacit{\'e} de l'inclusion $L^2_1\hookrightarrow
L_2$~; en  faisant tendre $j$ vers l'infini
$$\int_{M\setminus \{t\geq
  T\} }\langle\beta_j,\gamma_j\rangle\volj  \longrightarrow 
\int_{M\setminus \{t\geq
  T\} }\langle\beta,\mu\rangle\vol^g .$$
et il est facile d'en d{\'e}duire que $\beta $ est orthogonale {\`a} $\gamma$.
Finalement  $\beta$ est {\`a} la fois harmonique et orthogonale aux formes
harmoniques donc $\beta=0$.

Nous allons voir maintenant que nos hypoth{\`e}ses  impliquent  une
contradiction. Choisissons un $T>0$ suffisamment grand afin
d'{\^e}tre dans le cadre du corollaire \ref{corcontrloc2}. Soit $\chi(t)$
une fonction cut-off valant $1$ sur le compact $M\setminus \{t\geq
T\}$ et $0$ hors d'un compact. Alors on d{\'e}coupe les $\beta_j$ en
$$\|\beta_j\NLLj^2=1\leq 2\|\chi\beta_j \NLLj^2   +2\|(1-\chi)\beta_j\NLLj^2,$$
afin d'appliquer \ref{corcontrloc2} {\`a} $(1-\chi)\beta_j$~:
\begin{equation*}
  \begin{split}
c \|(1-\chi)\beta_j\NLLj^2  &\leq \|d (1-\chi)\beta_j\NLLj^2
+\|d^{*_j}(1-\chi)\beta_j\NLLj^2 \\
 &\leq \|d\beta_j\NLLj^2+ \|d^{*_j} \beta_j\NLLj^2 + 2\|(\dt\chi)
\beta_j\NLLj^2.     
  \end{split}
\end{equation*}
Puisque, $\beta_j$ converge vers $\beta=0$ au sens $L^2$-fort  sur
tout compact de $M$ et que par hypoth{\`e}se $\|d\beta_j\NLLj$ et $\|d^{*_j} \beta_j\NLLj$ tendent
vers $0$, il en r{\'e}sulte que  $\|\beta_j\NLLj$ tend vers $0$, ce qui contredit
l'hypoth{\`e}se $\|\beta_j\NLLj =1$.
\end{demode}

\subsubsection{Convergence des $2$-formes harmoniques}
Comme dans le cas des $1$-formes, on va utiliser l'in{\'e}galit{\'e} de
 Poincar{\'e} du  lemme \ref{courburemoyenne2} afin 
 d'{\'e}tudier la convergence des 
 repr{\'e}sentants  $g_j$-harmoniques d'une classe de cohomologie
 fix{\'e}e. 

\begin{prop} 
\label{propconvharm}
Soit une classe de cohomologie $b\in \HDR^2(\overline M)$. La suite
$\gamma_j$ de ses  repr{\'e}sentants $g_j$-harmoniques sur $\overline M$
converge sur tout compact de 
$M$ vers un repr{\'e}sentant $g$-harmonique $L^2$ de $b$.

Plus pr{\'e}cis{\'e}ment~:  soit un repr{\'e}sentant  $w$ de $b$  {\'e}gal {\`a} $\lambda
\volss$ pr{\`e}s de l'infini  sur le bout parabolique de $M$  (avec $\lambda\in\R$). On
{\'e}crit alors $\gamma_j =  w + d\beta_j$, o{\`u} $\beta_j$ est une $1$-forme lisse
 sur $\overline
M$  orthogonale aux formes $g_j$-harmoniques et telle que $d^{*_j}\beta_j=0$. 
Alors $\beta_j$ converge sur tout compact de $M$ vers une
$1$-forme $\beta$ dans $L^2_1(g)$ orthogonale aux formes
$g$-harmoniques  $L^2$, telle que $\gamma=w + d\beta $ soit
$g$-harmonique,   
$d^{*_g}\beta=0$ et $\quad \| \beta \|_{L^2_1} \leq c\| w
\NLLg$, o{\`u} $c>0$ est un constante fix{\'e}e.
\end{prop}
\begin{demo}
L'{\'e}criture de $\gamma_j$ sous la forme $\omega+ d\beta_j$ rel{\`e}ve de
la th{\'e}orie de Hodge {\'e}l{\'e}mentaire sur $\overline M$ et du lemme de
Poincar{\'e}~\ref{lemmePoincare}.
On remarque dans un premier temps que la forme $\volss$ est harmonique relativement aux m{\'e}triques $\gmod$ 
et $\gmodj$. On en d{\'e}duit le lemme suivant pour le m{\'e}triques $g$ et $g_j$~:
\begin{lemme}
  \label{lemmecontrolevolss}
Quel que soit $\epsilon >0$, il existe $T$ suffisamment grand tel que 
$$\int_{t\geq T} |d^{*_g}\volss |^2 \vol^g \leq \epsilon \quad et \quad
\int_{t\geq T}|d^{*_j}\volss|^2 \volj \leq \epsilon. 
$$
\end{lemme}
Par hypoth{\`e}se,  $\gamma_j = w + d\beta_j$ est $g_j$-harmonique, donc 
$d^{*_j}d\beta_j = -d^{*_j} w$.
En int{\'e}grant contre $\beta_j$, puis par in{\'e}galit{\'e} de
Cauchy-Schwarz, il vient 
$$ \|d\beta_j\NLLj^2 = -\int \langle
\beta_j, d^{*_j}\gamma \rangle\volj\leq \|\beta_j\NLLj
\|d^{*_j} w \NLLj~;
$$
par  l'in{\'e}galit{\'e} de Poincar{\'e} de la proposition
\ref{laplacien2}, on en 
d{\'e}duit  que
$$ \| \beta_j\NLLj^2 \leq c\|d\beta_j\NLLj^2 \leq c \|\beta_j\NLLj
\|d^{*_j} w \NLLj.
$$
D'apr{\`e}s le lemme
\ref{lemmecontrolevolss}, $ \|d^{*_j} w \NLLj$ est uniform{\'e}ment born{\'e}e,
ce qui entra{\^\i}ne que 
 $\|\beta_j\NLLj$  puis que  $\|d\beta_j\NLLj$
 sont uniform{\'e}ment born{\'e}es. Comme 
de plus $d^{*_j}\beta_j = 0$, on en d{\'e}duit que $\| \beta_j
\|_{L^2_1(g_j)}$ est uniform{\'e}ment born{\'e}e.
{\`A} partir de l{\`a}, on peut extraire  une limite
faible $\beta\in L^2_1 $ de $\beta_j$ sur tout 
compact de $M$ et terminer la d{\'e}monstration comme  dans le
lemme~\ref{convcompact}.

Seul le fait que la limite $\beta$ soit orthogonale aux formes
$g$-harmoniques $L^2$ n'est pas compl{\`e}tement {\'e}vident. Fixons une
classe de cohomologie $b\in \HLL(M)$. Soit $\omega_j$ la suite de ses
repr{\'e}sentants $g_j-$harmoniques et $\omega$ le repr{\'e}sentant
$g-$harmonique $L^2$. Nous allons voir que 
\begin{equation}
\label{exprcontroorth}\int \langle \omega_j,\gamma_j\rangle \volj \rightarrow \int \langle
\omega ,\gamma \rangle \volg \; ;\end{equation}
puisque la suite de gauche est identiquement nulle par hypoth{\`e}se,
on en d{\'e}duira que $\gamma$ est orthogonale {\`a} toute forme
$g-$harmonique $L^2$. 
{\`A} ce point, nous avons besoin d'un nouveau lemme technique~:
\begin{lemme}
\label{lemmecontrinfty2}
  Pour tout $\epsilon> 0$, et toute classe de cohomologie $b\in \HLL^
   2(M)$, il existe un $T$ suffisamment grand tel que avec les
   notations de \ref{propconvharm}, on ait 
$$\int_{t\geq T} |\gamma|^2\vol^g \leq \epsilon,\quad et \quad \int_{t\geq T} |\gamma_j|^2\volj \leq \epsilon.
$$
\end{lemme}
\begin{demode}{du lemme }
En utilisant le lemme~\ref{lemmecontrolevolss}, on en d{\'e}duit pour
tout $\epsilon>0$, que pour $T$ suffisamment grand, on a
$$\|\beta_j\NLLj \left ( \int_{t\geq T} |d^{*_j} w|^2 \volj\right
)^{1/2} \leq  \frac \epsilon 2 ,\quad\mbox{pour tout $j$.}
$$
En int{\'e}grant  par parties, il vient~: 
\begin{equation}
\int_{t\geq T} |d\beta_j|^2\volj 
 =  - \int_{t\geq T} \langle
d^{*_j} w ,\beta_j\rangle \volj +
\int_{t =  T} \beta_j\wedge *_j d\beta_j ,\label{expr1}
\end{equation}
d'o{\`u} la majoration
\begin{multline*}
  \int_{t\geq T} |d\beta_j|^2\volj   \leq 
\left |  \int_{t\geq T} \langle
d^{*_j}d\beta_j,\beta_j\rangle \volj \right |  + \left | 
\int_{t =  T} \beta_j\wedge *_jd\beta_j \right |   \\
\leq   \|\beta_j\NLLj \left ( \int_{t\geq T} \!\!\!\!\!|d^{*_j} w|^2 \volj\right
) ^{1/2}+ \left |  \int_{t =  T} \!\!\!\!\!\beta_j\wedge *_j
d\beta_j \right | 
\leq   \frac \epsilon 2 + \left |  \int_{t =  T}\!\!\!\!\! \beta_j\wedge *_j
d\beta_j \right |.
\end{multline*}
Puisque $\beta_j$ converge vers $\beta$ au sens $C^\infty$ sur tout
compact de $M$, 
le terme $\int_{t =  T} \beta_j\wedge *_jd\beta_j$ converge vers
$\int_{t =  T} \beta \wedge *d\beta$ lorsque $j$ tend vers l'infini~;
comme $\beta \in L^2_1(g)$, en utilisant la
formule~(\ref{expr1}) pour la m{\'e}trique $g$ et le
lemme~\ref{lemmecontrolevolss}, on peut
choisir une tranche $t=T$ pr{\`e}s de 
l'infini telle que $|\int_{t =  T} \beta \wedge *d\beta|\leq \epsilon
/ 4 $. Donc pour $j$ suffisamment grand  
$ |  \int_{t \geq   T}  \beta_j \wedge *d\beta_j  | \leq
 \epsilon / 2$, d'o{\`u} le lemme.
\end{demode}

\emph{Fin de la d{\'e}monstration de la proposition \ref{propconvharm}~:}
comme $\|\beta_j\NLLj $ est uniform{\'e}ment born{\'e}e, on peut choisir 
d'apr{\`e}s le lemme~\ref{lemmecontrinfty2} un compact $K$ tel que 
$$ \left | \int_{M\setminus K}\!\!\!\!\! \langle \gamma_j ,\beta_j\rangle \volj
\right | \leq
\|\beta_j\NLLj \left | \int _{M\setminus K} \!\!\!\! |\gamma_j| ^2\volj \right |^{1/2}\leq
 \frac \epsilon 2 $$
avec de plus 
$ \left | \int_{M\setminus K}  \langle \gamma_j ,\beta_j\rangle \volj
\right | \leq \frac \epsilon 2$. On d{\'e}coupe les int{\'e}grales
de~(\ref{exprcontroorth}) en deux morceaux suivant $M\setminus K$  o{\`u} on {\`a} le
contr{\^o}le ci-dessus et suivant  $K$ o{\`u} la convergence de
$\omega_j$ et de $\gamma_j$ est $C^\infty$~; on en d{\'e}duit 
facilement~(\ref{exprcontroorth}). 
\end{demo}

\emph{Convergence {\`a} l'infini.}
Les repr{\'e}sentants
$g_j$-harmoniques  $\gamma_j$ et $\omega_j$ 
de deux classes cohomologie fix{\'e}es  de $\HLL^2(\overline M)$ 
convergent sur tout compact de $M$ vers les
repr{\'e}sentants 
$g$-harmoniques $L^2$, $\gamma$ et $\omega$. 
Une autre cons{\'e}quence int{\'e}ressante du
lemme~\ref{lemmecontrinfty2} est la suivante~:
\begin{equation}
\label{exprconvcup}
\int_M \langle \gamma_j,\omega_j \rangle \volj \rightarrow \int_M
  \langle \gamma,\omega\rangle \vol^g \quad \mbox{ lorsque } j\rightarrow
  \infty .
\end{equation}
L'{\'e}toile de Hodge associ{\'e}e {\`a} $g$ 
induit un d{\'e}composition des $2$-formes $g$-harmoniques $L^2$ 
en leur partie autoduale et anti-autoduale. On en d{\'e}duit  
la  d{\'e}composition correspondante de $\HLL^2(M)$ en
\begin{equation}
\label{eqpol}
  \HLL^2(M) = H^{+_g}\oplus H^{-_g}.
\end{equation}
Par ailleurs nous avons   les d{\'e}compositions usuelles
pour les m{\'e}triques $g_j$ sur la vari{\'e}t{\'e} compactifi{\'e}e
$\overline M$ 
$$  \HDR^2(\overline M) = H^{+_j}\oplus H^{-_j}.$$
De (\ref{exprconvcup}), on d{\'e}duit alors le corollaire suivant~:
\begin{cor}
\label{convpol}
  Soit une classe de cohomologie $b\in \HLL^2(M)$, alors  
$$ b^{+_j} \rightarrow b^+\quad \mbox{ et }  \quad b^{-_j} \rightarrow b^-.$$
\end{cor}
\begin{demo}
Soit  $(\gamma_j)$ les repr{\'e}sentants
  $g_j$-harmoniques de $b$  et   $\gamma$ le
  repr{\'e}sentant $g$-harmonique $L^2$. Soit $(\gamma ^+)_j$ les
  repr{\'e}sentants $g_j$-harmoniques de $b^+$. 
Par la proposition~\ref{propconvharm},  la suite $(\gamma ^+)_j$ converge
  sur tout compact vers son repr{\'e}sentant $g$-harmonique qui est par hypoth{\`e}se  $\gamma ^
  + $. Parce que les m{\'e}trique $g_j$ convergent vers $g$ sur
  tout compact  on en d{\'e}duit que  $\gamma_j^{+_j} \rightarrow
  \gamma ^+$ sur tout compact. Alors
 pour toute classe de cohomologie $a\in \HLL^2(M)$, en notant
  $\omega_j$ la suite 
  de ses repr{\'e}sentants $g_j$-harmoniques, on a 
$$ (b^+ - b^{+_j})\cdot a = [(\gamma ^ +)_j - \gamma_j^{+_j} ]\cdot[\omega_j]
\leq \|(\gamma ^ +)_j - \gamma_j^{+_j}\NLLj \|\omega_j \NLLj.
$$
D'apr{\`e}s la proposition~\ref{propconvharm}, $\|\omega_j \NLLj$ est born{\'e}~;
  on en d{\'e}duit  que $(b ^+ - b^{+_j})\cdot a \rightarrow 0$
  d'apr{\`e}s (\ref{exprconvcup}). 
 Puisque le cup produit est une forme 
quadratique non d{\'e}g{\'e}n{\'e}r{\'e}e, ceci implique que $b ^+ - b^{+_j}$ tend vers $0$ dans $\HLL^2(M)$, d'o{\`u}
le corollaire.
\end{demo}

\section{{\'E}quations de Seiberg-Witten}
\label{secdemo}
Soit $M$ une surface complexe {\`a} bouts paraboliques
munie d'une m{\'e}trique $g$ asymptotique au mod{\`e}le local $\gmod$.
Nous commen{\c c}ons par expliquer comment g{\'e}n{\'e}raliser la
d{\'e}monstration de Le\,Brun~\cite{L} {\`a} ce cadre de volume fini au
\S\ref{secswm}, puis  comment obtenir une solution des
{\'e}quations de Seiberg-Witten pour $g$ 
par un proc{\'e}d{\'e} de convergence au \S\ref{secconv} dans le cas o{\`u}
les poids des bouts paraboliques v{\'e}rifient la condition de
rationalit{\'e}.

\subsection{{\'E}quations de Seiberg-Witten et bouts paraboliques}
\label{secswm}
La surface  complexe $M$ est munie d'une structure
$spin^c$ \emph{canonique} de  fibr{\'e} d{\'e}terminant $L=K_M^{-1}$
avec les  fibr{\'e}s de spineurs
$$W^+= \Lambda ^{0,0} M \oplus \Lambda
^{0,2} M \quad \mbox{ et }\quad W^-=\Lambda ^{0,1} M$$
munis d'une \emph{action de Clifford} de $TM$.

 Soit $g$ une m{\'e}trique \aparab sur $M$ et $A$ une connexion unitaire sur $M$. On en d{\'e}duit un op{\'e}rateur de Dirac $\Dirac_A$ et les 
{\'e}quations de Seiberg-Witten s'{\'e}crivent~:
 \begin{eqnarray}
  \label{eq:SW}
  \Dirac_A \psi & = & 0 \nonumber \\
  F_A ^{+_g}  &= &q(\psi),
\end{eqnarray}
o{\`u} $\psi\in \Gamma(W^+)$ et $q(\psi)$ est la partie dans trace de
l'endomorphisme $\psi\otimes\psi ^*$. La courbure $F_A^+$ est une
forme imaginaire autoduale et agit par multiplication de Clifford
comme un endomorphisme hermitien sans trace, ce qui donne un sens {\`a}
la deuxi{\`e}me {\'e}quation.

\subsubsection{Connexions $L^2_1$}
Sur les bouts paraboliques de $M$, le mod{\`e}le local de m{\'e}trique
k{\"a}hl{\'e}rienne $\hat g $ induit une connexion $\hat A$ sur le fibr{\'e}
anti-canonique $K_M^{-1}$ que l'on {\'e}tend par partition de l'unit{\'e} au
fibr{\'e} tout entier.
  On d{\'e}finit alors \emph{la classe de Chern $L^2$} du fibr{\'e} anti-canonique de $M$
comme la classe de cohomologie de $\HLL^ 2(M)$
\begin{equation}
\LLc = \frac i{2\pi} \left [ F_{\hat A}\right ] .
\end{equation}
Si $A$ est une connexion sur $K_M^ {-1}$ telle que $A = \hat A + a$
avec $a$ et $da\in
L^2$, alors  $\LLc = \frac i{2\pi} [F_{A }]$.
En particulier si $g$ est une m{\'e}trique de K{\"a}hler \aparab sur les bouts paraboliques
de $M$, elle induit alors  une
connexion $A ^g$   sur $K_M^{-1}$~ telle que  $A ^g - \hat A \in
L^2_1$ d'o{\`u} $\LLc = \frac i{2\pi} [F_{A ^ g}]$.

\subsubsection{Un lemme clef}
\label{seclemmeclef}
Le lemme utilis{\'e} par Le\,Brun dans
\cite{L} pour d{\'e}montrer le th{\'e}or{\`e}me~\ref{theoA} a un analogue 
dans notre contexte de volume fini {\`a} 
condition de faire des hypoth{\`e}ses de r{\'e}gularit{\'e} {\`a} l'infini sur
la solution des {\'e}quations de Seiberg-Witten~(\ref{eq:SW}).
\begin{lemme}
\label{proplebrun}
  Supposons que les {\'e}quations de Seiberg-Witten non
  perturb{\'e}es (\ref{eq:SW}) sur $M$,  associ{\'e}es {\`a} une m{\'e}trique
  \aparab $g$, admettent une solution $(A,\psi)$ telle que
 $\psi \not \equiv 0 $,
 $A = \Amod + a $, avec $a\in L^2_1$,
    et  $\psi \in L^2_1$.
Alors 
$$ \frac 1{32\pi^2}\int_M s^2_g \vol^g \geq \left ( \LLc^{+_g} \right )^2 . $$
En outre, le cas d'{\'e}galit{\'e} est r{\'e}alis{\'e} si et seulement si 
$ \nabla_A\phi=0$ et $s=cste < 0$.
Dans ce cas, $g$ est une m{\'e}trique de K{\"a}hler
relativement {\`a} une structure complexe $J$  induite par la
m{\'e}trique $g$ et la forme de K{\"a}hler $F_A^+/\sqrt 2 |F_A^+|$.
\end{lemme}
\begin{demo}
Identique {\`a} celle {\`a} celle des surfaces compactes, bas{\'e}e sur la
formule de Lichnerowicz (cf. \cite{L}),
car les hypoth{\`e}ses de r{\'e}gularit{\'e} sur
$(A,\psi)$ nous permettent de faire les int{\'e}grations par
parties.
\end{demo}

\subsection{D{\'e}monstration du th{\'e}or{\`e}me \ref{theoA}}
\label{subdemotheoa}
Dans tout le \S\ref{subdemotheoa}, $M=\PP(\mathcal E)_\Sigma$ est une surface r{\'e}gl{\'e}e
 obtenue {\`a} partir d'un fibr{\'e} parabolique $\mathcal E$. Pour
 simplifier, on supposera que $M$ poss{\`e}de au moins un bout
 parabolique et on se r{\'e}f{\'e}rera {\`a} \cite{BB} et \cite{L} pour
 les d{\'e}monstrations dans le cas compact.
 
\subsubsection{Une antipodie sur $M$}
Soit $h$ la m{\'e}trique hermitienne sur $\mathcal E$ adapt{\'e}e {\`a} la
 structure parabolique qui nous a permis de construire le mod{\`e}le
 local $\gmod$. On en d{\'e}duit une m{\'e}trique de Fubini-Study sur
 chaque fibre de $M$ et une \emph{antipodie} $\xi : M\rightarrow \Sigma$,
qui {\`a} chaque
 point d'une fibre de $M$ associe le point diam{\'e}tralement
 oppos{\'e}. L'antipodie
 $\xi$ est clairement une involution du fibr{\'e}
 $M\rightarrow\Sigma$ renversant  l'orientation et
 une isom{\'e}trie pour le
mod{\`e}le local de m{\'e}trique $\gmod$~;
on d{\'e}duit de cette derni{\`e}re propri{\'e}t{\'e} que si $g$ est
\aparabp, alors $\xi ^ *g$ est {\'e}galement \aparab  et que $\xi$
 agit donc sur la cohomologie $\HLL^*(M)$.
\begin{lemme}
\label{lemmeorthog}
  L'antipodie $\xi$ agit sur $\HLL^2(M)$ en {\'e}changeant les facteurs de la
  d{\'e}composition
 $$\HLL^2(M)=\HLL^+(M)\oplus \HLL^-(M).$$
\end{lemme}
\begin{demo}
Comme $\xi$ est un diff{\'e}omorphisme de $M$ qui renverse
l'orientation, on en d{\'e}duit que pour des $2$-formes $L^ 2$,
$\gamma_1$ et 
$\gamma_2$,
$$ \int_M \xi ^ * \gamma_1 \wedge  \xi ^ * \gamma_2 =  - \int_M \gamma_1 \wedge  \gamma_2.
$$
Donc pour 
$a,b\in\HLL^2(M)$,
on a $ \xi ^* a \cdot \xi ^* b = -a\cdot b $~;
autrement dit $\xi$ ne fait que changer le signe de la forme
d'intersection.
En utilisant le fait que $\xi$ est une involution,  on en d{\'e}duit que 
$  a\cdot \xi^*a = 0$.
Puisque nous sommes dans la situation $b_2^+=b_2^-=1$, cette derni{\`e}re
identit{\'e} implique le lemme.
\end{demo}
\begin{cor}
\label{corpol}
  Les m{\'e}trique $g$ et $\xi^*g$ d{\'e}finissent la m{\^e}me d{\'e}composition 
$$
  \HLL^2(M)= \HLL^{+}(M)\oplus \HLL^{-}(M).
$$
\end{cor}
\begin{demo}
Soit une classe de cohomologie $a\in
\HLL^{+_g}(M)$. D'apr{\`e}s le lemme  \ref{lemmeorthog}, $\xi^*a\in\HLL^{-_g}(M)$~;
soit $\gamma$ le repr{\'e}sentant $L^2$, $g$-harmonique anti-autodual de $\xi^*a$. Alors 
$\xi^*\gamma$ repr{\'e}sente $a$ et c'est une forme $L^2$,
$\xi^*g$-harmonique autoduale,  d'o{\`u} $a\in \HLL^{+_{\xi^*g}}(M)$. 
Les autres inclusions se d{\'e}montrent de la m{\^e}me fa{\c c}on.
\end{demo}

\subsubsection{Fin de la preuve} Voici maintenant une proposition
essentielle dont on d{\'e}duit automatiquement  que
$(M,\gkahl) \simeq
\Sigma\times_\rho \CP$ (cf. \cite{L}).  Par le th{\'e}or{\`e}me de Mehta-Seshadri,
il en r{\'e}sulte que $\mathcal E$ est polystable ce qui d{\'e}montre le
th{\'e}or{\`e}me \ref{theoA}. 
\begin{prop}
\label{proprevuniv}
  Sous les hypoth{\`e}ses du th{\'e}or{\`e}me \ref{theoA}, le
  rev{\^e}tement 
  universel holomorphe riemannien de $(M,\gkahl)$ donn{\'e} par
  $\HH\times \CProj^1$, muni de sa m{\'e}trique k{\"a}hl{\'e}rienne standard
  {\`a} courbure scalaire constante $s=2(c-1)$.
\end{prop}
\begin{demode}{dans le cas o{\`u} $s  =0$}
On d{\'e}finit le nombre caract{\'e}ristique
$$\sigma = \int_M \gamma
\quad \mbox{ o{\`u} }\quad   \gamma= -\frac 1{24\pi ^2}\left [ \trace (R^\gkahl\wedge
R^\gkahl) \right ]. $$
En dimension $4$, on en d{\'e}duit que 
\begin{equation}
\label{exprgaussbonnet} 
\sigma = \frac1{12\pi ^2}\int_M |W^+|^2 - |W^-|^2 \; \vol^\gkahl,
\end{equation}
o{\`u} $W^\pm$ sont les parties autoduales et anti-autoduales du
tenseur de Weyl. Une surface k{\"a}hl{\'e}rienne est anti-autoduale
si et seulement si
$s=0$~;  Maintenant, il est facile de calculer $\sigma$~: pour
cela, on peut quitte {\`a} changer la structure holomorphe du fibr{\'e}
parabolique $\mathcal E\rightarrow \overline \Sigma$ supposer que
c'est un fibr{\'e} paraboliquement stable. On en d{\'e}duit une
m{\'e}trique de K{\"a}hler $\hat g$
obtenue par construction standard {\`a}
courbure scalaire $s=0$. Par la formule de transgression, $\sigma$ ne
d{\'e}pend que de la classe de 
cohomologie $L^2$ repr{\'e}sent{\'e}e par $\gamma$ donc
$$\sigma=\int_M \gamma =\int_M \hat\gamma 
\quad \mbox{ o{\`u} } \quad 
 \hat \gamma=-\frac 1{24\pi ^2}\left [ \trace (R^{\hat g}\wedge
R^{\hat g}) \right ].
$$
  
Puisque $\hat g$ est un produit
local de deux m{\'e}triques de K{\"a}hler de courbure oppos{\'e}es sur
des surface de Riemann, 
elle est donc conform{\'e}ment plate et d'apr{\`e}s (\ref{exprgaussbonnet}) il
en r{\'e}sulte que $\sigma = 0$. Finalement, on en d{\'e}duit que
$W^-=0$  pour
la m{\'e}trique $\gkahl$ et qu'elle est elle aussi conform{\'e}ment plate.
Or une surface  k{\"a}hl{\'e}rienne conform{\'e}ment plate 
est n\'cessairement un produit local de
deux m{\'e}triques k{\"a}hl{\'e}riennes {\`a} courbures scalaires constantes
oppos{\'e}es $(X_1,g_1)\times(X_2,g_2)$. En raison du comportement
asymptotique de la m{\'e}trique l'un des facteurs doit {\^e}tre $\CP$ et
l'autre $\HH$.
\end{demode}

\begin{demode}{dans le cas o{\`u} $s<0$}
On commence par le lemme
\begin{lemme}
\label{lemmerevuniv}
  Sous les hypoth{\`e}ses du th{\'e}or{\`e}me \ref{theoA} et en supposant
  que les {\'e}quations de Seiberg-Witten non perturb{\'e}es associ{\'e}es
  {\`a}  la m{\'e}trique
  $g=\xi ^*\gkahl$ ont une solution $(A,\psi)$ avec la r{\'e}gularit{\'e}
  sp{\'e}cifi{\'e}e  dans le lemme \ref{proplebrun}, on en d{\'e}duit que
le rev{\^e}tement holomorphe riemannien de $M$ est {\'e}gal {\`a}
$\HH\times \CProj^1$ muni de sa m{\'e}trique standard.
\end{lemme}
\remarque
la fin de cet article a pour but de la d{\'e}montrer l'existence de la
solution $(A,\phi)$ (cf.
 \S\ref{secconv}  et \S\ref{secirrat}).

La m{\'e}trique de K{\"a}hler {\`a} courbure scalaire constante
$\gkahl$ induit une connexion $A ^{\mathrm{K}}$ sur le fibr{\'e} anti-canonique
dont la courbure est harmonique et v{\'e}rifie 
$$ F_{A ^{\mathrm K}}^{+_\gkahl}= -\frac {is \omega }4,
$$
avec  $\omega$  la forme de K{\"a}hler de $\gkahl$~;
on en d{\'e}duit que $\LLc^{+_\gkahl} = \frac s{8\pi } [\omega ]$ d'o{\`u} 
$$   (\LLc^{+_\gkahl})^2= \frac 1{32\pi ^2}\int_M s^2 \vol^\gkahl. $$
Par d{\'e}finition $g= \xi ^*\gkahl$, donc
 par changement de variable 
$\int_M s_g^2 \vol^g = \int_M
s^2 \vol^\gkahl $~; en utilisant le 
corollaire~\ref{corpol}
$$ ( \LLc^{+_{\gkahl}})^2 = ( \LLc ^{+_{g}} )^2 =  \frac 1{32\pi
^2}\int_M s^2 \vol^g. $$
En  appliquant le lemme \ref{proplebrun} dans le cas
d'{\'e}galit{\'e}, on en d{\'e}duit que $g $ est de K{\"a}hler
relativement {\`a} une structure complexe $J$ compatible avec l'orientation de $M$.
 En prenant l'image r{\'e}ciproque de ces structures par  $\xi$, il
vient que
$\gkahl=\xi ^* g$ est
k{\"a}hl{\'e}rienne relativement {\`a} une structure complexe $J_1=\xi^* J$  compatible
avec l'orientation inverse de $M$. 
En notant $J_0$ la structure complexe originelle sur $M$ qui commute
avec $J_1$ on d{\'e}finit un automorphisme  de carr{\'e}
$1$  de $TM$ donn{\'e} par $J_0J_1$. 
On obtient une   d{\'e}composition parall{\`e}le de $TM$ en $L_0\oplus
L_1$ suivant les espaces propres de $J_0J_1$.
Finalement le rev{\^e}tement universel holomorphe riemannien de
$(M,g)$ est un produit $(X_1,g_1)\times (X_2,g_2)$.  Puisque la
courbure scalaire de $\gkahl$ est constante et que 
$s(x_1,x_2)=s_1(x_1)+s_2(x_2)$,
on en d{\'e}duit que $s_1$ et
$s_2$ sont constantes. Maintenant, en raison du comportement
asymptotique de $\gkahl$ 
l'un des deux facteurs
doit {\^e}tre {\`a} courbure sectionnelle  $-1$ et l'autre {\`a} courbure
sectionnelle $c$. 
\end{demode}

\subsection{Convergence des solutions des {\'e}quations de
Seiberg-Witten}
\label{secconv}
Notre objectif est  maintenant d'obtenir une solution non r{\'e}ductible
aux {\'e}quations de Seiberg-Witten (\ref{eq:SW}) pour une surface
complexe r{\'e}gl{\'e}e {\`a} bouts paraboliques, car c'est un outil
central dans la d{\'e}monstration 
du th{\'e}or{\`e}me \ref{theoA}. Dans cette optique, nous 
montrons que l'on sait faire converger une suite de
solutions  
des  {\'e}quations de 
Seiberg-Witten perturb{\'e}es~(\ref{eq:SWp}) pour les m{\'e}triques
$g_j$, vers une 
solution des {\'e}quations non perturb{\'e}es~(\ref{eq:SW}) pour la m{\'e}trique
$g=\xi ^*\gkahl$.  
Plus g{\'e}n{\'e}ralement, le th{\'e}or{\`e}me de
convergence~\ref{theoconvergence} s'applique {\`a} toute m{\'e}trique
$g$ \aparab sur une surface complexe.

Nous supposerons pour donc pour l'instant que $M$ est une surface
complexe {\`a} bouts paraboliques v{\'e}rifiant la condition
$\alpha_2-\alpha_1\in \Rat$. Nous expliquerons au \S\ref{secirrat}
comment se passer de cette hypoth{\`e}se. Pour simplifier 
nous supposons en outre que $M$ poss{\`e}de exactement un bout
parabolique, les d{\'e}monstrations s'{\'e}tendant trivialement au cas de
plusieurs bouts.
Soit $g$ une m{\'e}trique sur $M$ \aparab $\gmod$ et  $g_j$ sont ses
approximations 
lisses sur la compactification orbifold $\overline M =
M\cup D$, o{\`u} $D=\CP/\ZZ_q$.

\subsubsection{Compactification du fibr{\'e} $K_M^{-1}$}
Rappelons que pr{\`e}s de $D$, on a rev{\^e}tement local 
{\`a} $q$ feuillets de $\overline M$ donn{\'e} par 
\begin{eqnarray*}
\pq : I_a/\langle \tau ^ q \rangle \times \CProj ^1& \rightarrow &
I_a/\langle \tau \rangle \times \CProj ^1 \\
(x+iy ,[\tilde u : \tilde v]) &\mapsto& (x+ iy ,[\tilde u e
^{-i\alpha_1 x }: \tilde v e ^{-i\alpha_2 x}]),
\end{eqnarray*}
que l'on prolonge  en  un rev{\^e}tement ramifi{\'e} $\pq :
\Delta_{a/q}\times \CProj^1\rightarrow \overline \bparab $ tel que
$\pq (\{0\}\times \CProj^1)= D$.

Le fibr{\'e}  en droites complexes $[D]$ associ{\'e} au diviseur $D$ tel que   $c_1([D])$ est le
dual de Poincar{\'e} de $D$, est d{\'e}fini classiquement 
{\`a} l'aide de deux ouverts de trivialisation, l'un au dessus de $M$,
l'autre au dessus de $\Delta_{a/q} \times
\CProj ^1$
(dans le cas
des orbifolds, il faut se placer sur le rev{\^e}tement local
ramifi{\'e} pour que la construction ait un sens).
En utilisant les coordonn{\'e}es $(\xi=x+iy, [\tilde u : \tilde v])$, on
recolle ces deux trivialisations par la fonction de transition $ \rho
= e ^ {ix}$ qui est invariante sous l'action de $\ZZ_q$ et d{\'e}finit
ainsi un orbifibr{\'e} sur $\overline M$. Notons  que le fibr{\'e}
$[D]$ restreint {\`a} $M$ devient {\'e}videmment trivial.

On construit  sur $[D]$
des  connexions d{\'e}finies 
explicitement, 
en utilisant les notations de
 \S\ref{seccomp}, par
\begin{equation}
\label{exprBj}
B_j = d - i\chi_j(\dt\phi_j) d\theta.
\end{equation}
Ceci d{\'e}finit  bien une connexion lisse sur $[D]$ puisque $\dt\phi_j=-1$ pr{\`e}s
de $[D]$ d'o{\`u}
$\rho
 B_j \rho ^{-1}= d$ pr{\`e}s de $D$. 
Alors la $2$-forme  $\varpi_j= \frac i{2\pi} F_{B_j}$
repr{\'e}sente $c_1([D])$ et par construction,  la suite de connexions $B_j$
converge sur tout compact de $M$ vers la
connexion triviale.

\begin{lemme}
\label{lemmeconnl2}
La classe de Chern $L^2$ du fibr{\'e} $K^{-1}_M$ se d{\'e}compose en
  $$\LLc = c_1(K^{-1}_{\overline M}) +  c_1([D]^{-1}).$$
Plus pr{\'e}cis{\'e}ment, la connexion $\hat A$ sur $K_M^ {-1}$ 
s'{\'e}crit sous la forme $\hat A = C + a$,
o{\`u} 
 $a$ est une $1$-forme $ L_1^2(g)$ sur
$M$
et $C$ est une connexion lisse sur le fibr{\'e} $L = K_{\overline M}^
{-1}\otimes [D]^ {-1}$ dont la  courbure
v{\'e}rifie pr{\`e}s de $D$ 
$$F_C= -ic\volss,$$
avec  $c$ la courbure de la m{\'e}trique de Fubini-Study.
\end{lemme}
\begin{demo}
Dans le rev{\^e}tement  $\pq : \Delta ^*\times \CProj^ 1\rightarrow
\bparab $, 
la connexion
de Chern sur le fibr{\'e} anti-canonique  se d{\'e}compose  en  
 $ A ^{\Delta ^*} \otimes  A ^\mathrm{FS}$ avec  $A
 ^{\Delta ^*}$ et  $ A ^\mathrm{FS}$ 
 les connexions de Chern sur les fibr{\'e}s
 tangents de $\Delta    ^*$ et
 $\CProj^1$. Mais
$A^{\Delta^*}= d - \partial\ln (|dz|^2 ) = d - \partial\ln
( |z| ^2) - \partial \ln (\ln ^ 2 |z|  ) 
$
et le dernier terme est $L^2_1$~; on en d{\'e}duit le lemme en utilisant la
formule de Poincar{\'e}--Lelong.
\end{demo}

\subsubsection{{\'E}quations perturb{\'e}es pour les m{\'e}triques $g_j$}
La compactification $\overline M$ de $M$ est une surface complexe~;
elle est donc munie d'une structure $spin^c$ canonique de
fibr{\'e} d{\'e}terminant $L_0=K_{\overline M}^{-1}$ avec pour fibr{\'e}s
de spineurs 
$$W^+= \Lambda ^{0,0}\overline M \oplus \Lambda
^{0,2}\overline M \quad \mbox{ et }\quad W^-=\Lambda ^{0,1}\overline
M.$$
 Notons que
cette structure $\spinc$ co{\"\i}{}ncide avec celle des {\'e}quations non
perturb{\'e}e si on la restreint {\`a} $M$.

Chaque m{\'e}trique $g_j$ sur $\overline M$ induit un op{\'e}rateur de
Dirac modulo le choix d'une connexion unitaire $A$ sur $L_0$~; on
d{\'e}finit  alors les {\'e}quations de Seiberg-Witten
perturb{\'e}es par 
\begin{eqnarray}
  \label{eq:SWp}
  \Dirac_A \psi & = & 0 \nonumber \\
  (F_A  - F_{B_j} )^{+_{g_j}}  &= &q(\psi),
\end{eqnarray}
o{\`u} $\psi\in \Gamma (W^+)$. Le groupe de
jauge $\jauge =Map(\overline M,S^1)$ agit par 
$$ f\cdot(A,\psi) = (A-2\frac {df}f, f\psi)
$$ sur l'espace des solutions des
{\'e}quations de Seiberg-Witten. 
D'apr{\`e}s le lemme \ref{lemmeconnl2},
$[F_A-F_{B_j}] =  -2i\pi
c_1(L_0\otimes [D]^ {-1}) = -2i\pi \LLc$~;
c'est pr{\'e}cis{\'e}ment pour cette raison que nous avons perturb{\'e} les
{\'e}quations
et peut alors extraire une solution des {\'e}quations de Seiberg-Witten
non perturb{\'e}es associ{\'e}es {\`a} la m{\'e}trique $g$ gr{\^a}ce au
th{\'e}or{\`e}me suivant.
\begin{theo}
\label{theoconvergence}
  Soit $g$, une m{\'e}trique \aparab sur $M$, et soit $g_j$ les
  m{\'e}triques sur $\overline M$ qui approximent $g$. Soit $C$, la connexion  de r{\'e}f{\'e}rence
   sur $L=L_0\otimes [D]^{-1}$ d{\'e}duite de
  la connexion $\hat A$ {\`a} l'aide du lemme \ref{lemmeconnl2}.
Supposons donn{{\'e}e}  une suite $(A_j,\psi_j)$ de solutions des  {\'e}quations de Seiberg-Witten
  perturb{\'e}es (\ref{eq:SWp}) associ{\'e}es aux m{\'e}triques $g_j$. Alors, quitte {\`a} extraire une sous-suite et {\`a} faire des changements de jauge, les $(A_j,\psi_j)$
  convergent au sens $C^\infty$ sur tout compact de $M$ vers une
  solution $(A,\psi)$ des {\'e}quations de Seiberg-Witten non
  perturb{\'e}es (\ref{eq:SW}) sur  $(M,g)$ telle que~:
  \begin{itemize}
  \item $A = C + a$, o{\`u} $a$ est
    une $1$-forme  $L^2_1(g)$, v{\'e}rifiant $d^{*_g}a =0$.
  \item $\psi \in L^2_1(g)$.
  \end{itemize}
\end{theo}\medskip

\subsubsection{Contr{\^o}le $C^0$}
Le premier pas dans la d{\'e}monstration du th{\'e}or{\`e}me
\ref{theoconvergence} est un
r{\'e}sultat de contr{\^o}le $C^0$ uniforme sur $\psi_j$~:
\begin{lemme}
\label{lemmec0}
Il existe une constante $K$ telle que pour tout $j$ et pour toute  solution $(A_j,\psi_j)$ des {\'e}quations de
  Seiberg-Witten perturb{\'e}es associ{\'e}es {\`a} la 
  m{\'e}trique $g_j$, on ait  $|\psi_j|\leq K$.
\end{lemme}
Pour le voir il faut  analyser la courbure de $B_j$ qui est explicite.
\begin{lemme}
\label{courbureBj}
  La courbure de la connexion $B_j$ v{\'e}rifie 
$$F_{B_j}=-i \chi_j \frac{\dt^2 \phi_j}{\phi_j} dt\wedge \phi_j
d\theta +F^b_j, 
$$
avec $|F^b_j|$ born{\'e}e ind{\'e}pendamment de $j$.
\end{lemme}
\begin{demode}{du lemme \ref{lemmec0}}
  Par la formule de Lichnerowicz 
$$\Dirac_{A_j}^2\psi_j   =  \nabla_{A_j}^*\nabla_{A_j}\psi_j + \frac {s^{g_j}}4
\psi_j +\frac 12 F_{A_j}^+.\psi_j  $$
puis en utilisant les {\'e}quations de Seiberg-Witten, il vient
\begin{equation}
\label{eq:lichne1}
 \nabla_{A_j}^*\nabla_{A_j}\psi_j + \frac {s^{g_j}}4
\psi_j +\frac 12 F_{B_j}^+.\psi_j +\frac 12 q(\psi_j)\psi_j =0.
\end{equation}
Rappelons que d'apr{\`e}s les lemmes \ref{courburescalaire} et  \ref{courbureBj},
$$  s^{g_j}=s^{g_j}_b-2\chi_j\frac{\dt^2\phi_j}{\phi_j}, \quad   F_{B_j}=-i \chi_j \frac{\dt^2 \phi_j}{\phi_j} dt\wedge \phi_j
d\theta +F^b_j, $$
avec $s^{g_j}_b$ et  $|F^b_j|$  uniform{\'e}ment born{\'e}es. Par ailleurs $(dt \wedge
\phi_j d\theta)^{+_j}$  agit par produit de Clifford sur $W^+$ avec des
valeurs propres $\pm i$, ${\phi_j^{-1}} {\partial_t^2\phi_j}$ reste
born{\'e}e sup{\'e}rieurement mais pas inf{\'e}rieurement donc l'{\'e}quation
(\ref{eq:lichne1}) peut s'{\'e}crire  sous la forme 
$$0=\nabla_{A_j}^*\nabla_{A_j}\psi_j + P_j\psi_j+ P_j^b\psi_j
+\frac 12 q(\psi_j)\psi_j , 
$$
o{\`u} sur le bout parabolique
$$P_j=-\frac 12 \chi_j\frac{\dt^2\phi_j}{\phi_j}\psi_j-\frac i2  \chi_j \frac{\dt^2 \phi_j}{\phi_j} (dt\wedge \phi_j
d\theta)^{+_j}\psi_j $$
est un op{\'e}rateur dont les valeurs propres sont $0$ et
$-\phi_j^{-1}\chi_j\dt^2\phi_j$. L'op{\'e}rateur lin{\'e}aire $P_j^b$ est 
uniform{\'e}ment born{\'e}. Il en r{\'e}sulte que  la partie
n{\'e}gative de ses valeurs 
propres de $P_j+P_j^b$ reste born{\'e}e inf{\'e}rieurement.

On peut maintenant d{\'e}montrer le lemme en utilisant le principe du
maximum~: si $x_0\in \overline M$ un maximum local de $|\psi_j |^2$,
alors $\Delta ^{g_j} ( |\psi_j |^2 )(x_0) \geq 0$. Or d'apr{\`e}s l'identit{\'e}
$ \Delta ^{g_j} ( |\psi_j |^2 )= 2\Re \langle
\nabla_{A_j}^*\nabla_{A_j}\psi_j,\psi_j\rangle - 2|\nabla_{A_j}\psi_j |^2
$, on en d{\'e}duit qu'il existe une constante
$\kappa > 0$ telle que 
$$ 0 \geq - \kappa |\psi_j(x_0) |^2  +\frac 12 \langle
q(\psi_j)\psi_j(x_0),\psi_j(x_0)\rangle =  - \kappa |\psi_j(x_0) |^2  +\frac 14 |\psi_j(x_0) |^4,
$$
d'o{\`u} $|\psi_j(x_0) |^ 2 \leq  4\kappa$.
\end{demode}

\subsubsection{Convergence des connexions}
Notons   $C_j= A_j \otimes B_j^{-1}$ la suite de connexions
d{\'e}finies sur $L$.
On utilise
la jauge donn{\'e}e par la connexion $C$ de <<r{\'e}f{\'e}rence>>   du
lemme \ref{lemmeconnl2}  et on {\'e}crit 
$$C_j= C+ \beta_j + \mu_j,
$$
o{\`u} $\mu_j$ est une $1$-forme $g_j$-harmonique, et $\beta_j$, une
$1$-forme orthogonale aux formes $g_j$-harmoniques.
\medskip

\emph{Convergence de la partie harmonique.}
D'apr{\`e}s le th{\'e}or{\`e}me des coefficients universels,
valable pour les  $CW-$complexes et en particulier pour les orbifolds,
$H^1(\overline M,\ZZ)$ est un r{\'e}seau de $H^1(\overline M,
\R)$ et
finalement
 le quotient 
$$H^1(\overline M,\R)/H^1(\overline M, \ZZ)$$
 est un tore  compact.
Par  ailleurs, 
les composantes connexes du  groupe de jauge $Map(\overline M ,S^1)$
agissent comme le r{\'e}seau $i H^1(\overline M , \ZZ)$ sur  les classes
de cohomologie de $iH^1(\overline M, \R)$.
On en d{\'e}duit que quitte {\`a} faire agir le groupe de jauge sur les
solutions des {\'e}quations de Seiberg-Witten $(A_j,\psi_j)$, on peut
supposer que les classes de cohomologie d{\'e}finies par 
$\mu_j$ restent born{\'e}es~; quitte {\`a} extraire une sous-suite
de $(A_j,\psi_j)$, on peut donc supposer que  $[\mu_j]$ converge. D'apr{\`e}s
la proposition 
\ref{convcompact}, il en r{\'e}sulte que $\mu_j$ converge au sens
$C^\infty$ sur tout compact de  
$M$, vers une forme  $g$-harmonique $L^2$.
\medskip

\emph{Convergence de la partie orthogonale aux formes harmoniques.}
Pour commencer on peut quitte {\`a} faire des changements de jauge dans la composante connexe de
l'identit{\'e}, ce qui n'affecte pas l'argument pr{\'e}c{\'e}dent, on
peut se placer  dans une \emph{jauge de Hodge} telle que $d^{*_j}\beta_j = 0$.
 Puisque $(A_j,\psi_j)$ est solution des
{\'e}quations de Seiberg-Witten, 
$$ F^{+_j}_{C_j}= (F_C + d\beta_j)^{+_j}=q(\psi_j),
$$
d'o{\`u}  $ d^{+_j}\beta_j = q(\psi_j)- F_C^{+_j}$.
Le lemme~\ref{lemmec0} donne alors une borne $C^0$ sur
 $q(\psi_j)$ et
comme la connexion de r{\'e}f{\'e}rence $C$ v{\'e}rifie
$F_C=\lambda\volss$
pr{\`e}s de $D$, elle admet  {\'e}galement  une borne $C^0$~; en utilisant le fait que le volume des m{\'e}triques
est  
uniform{\'e}ment born{\'e} par le lemme~\ref{volume}, 
on en d{\'e}duit une borne uniforme sur $\| d^{+_j}\beta_j\NLLj$.
Sur une vari{\'e}t{\'e} compacte 
$$2\|d^{+_j}\beta_j\NLLj^2 = \|d\beta_j\NLLj^2~; $$ 
alors, d'apr{\`e}s
la proposition~\ref{laplacien2}
$$ c \|\beta_j\NLLj^2 \leq \| d\beta_j\NLLj^2. $$  
Finalement, Par cons\'equent,  $\|\beta_j\|_{L^2_1(g_j)}$ est
uniform\'ement born\'e et on peut extraire une limite faible sur tout
compact $\beta\in L^2_1(g)$ de
$\beta_j$ orthogonale aux formes
$g$-harmoniques $L^2$ et v{\'e}rifiant $ d^{*_g}\beta = 0$.

En conclusion, les connexions $C_j= C+\mu_j+\beta_j$  convergent au sens
$L^2_1$-faible sur tout 
compact de $M$  vers une connexion $A = C+a$ avec $a=\mu+\beta \in
L^2_1(g)$. Par construction, $B_j$ converge vers la connexion
triviale sur tout compact de $M$~; on en d{\'e}duit que $A_j$ converge
{\'e}galement vers $C$ sur tout compact.
\medskip

\emph{Convergence de la partie spineur.}
Comme dans la d{\'e}monstration du lemme \ref{lemmec0}, on {\'e}crit la
formule de Lichnerowicz sous la forme 
\begin{equation}
\label{exprintermlichn}
0 = D_{A_j}\psi_j = \nabla_{A_j} ^*\nabla_{A_j}
\psi_j + (P_j ^b + P_j)\psi_j + \frac 12 q(\psi_j)\psi_j,
\end{equation}
o{\`u} $P^b_j$ est un op{\'e}rateur uniform{\'e}ment born{\'e}, et $P^b_j$ un
op{\'e}rateur dont la partie n{\'e}gative des  valeurs propres est  uniform{\'e}ment
minor{\'e}e. Alors, il existe une constante $\kappa>0$ telle que 
$ - \langle (P_j ^b + P_j)\psi_j,\psi_j\rangle \leq \kappa 
|\psi_j|^2\; ; $  
en int{\'e}grant sur $M$ la formule (\ref{exprintermlichn}) contre
$\psi_j$, on obtient le contr{\^o}le
$$ \|\nabla_{A_j}\psi_j\NLLj^2 \leq \kappa \|\psi_j\NLLj^2 - \frac 14
\int_M
|\psi_j|^4\volj   \leq \kappa \|\psi_j\NLLj^2 .
$$
Des  bornes uniformes $C^0$ sur $\psi_j$ et sur le volume des
m{\'e}triques $g_j$, on d{\'e}duit alors une borne
uniforme sur $\|\nabla_{A_j}\psi_j\NLLj$. 
On peut donc extraire de $\psi_j$ une limite faible sur tout compact $\psi$ telle
que $\psi, \nabla_C\psi\in L^2(g)$.
Finalement, $(A ,\psi)$ est 
obtenue comme une  limite faible sur tout compact  de solutions des
{\'e}quations de 
Seiberg-Witten perturb{\'e}es $(A_j,\psi_j)$. Or $B_j$ tend vers la
connexion triviale sur tout compact de $M$
 il en r{\'e}sulte que $(A, \psi)$ est solution {\`a}
la limite  des {\'e}quations de Seiberg-Witten non
perturb{\'e}es~(\ref{eq:SW}). 

Une fois  ces contr{\^o}les  $L^2_1$  obtenus, 
un proc{\'e}d{\'e} classique de
bootstraping utilisant l'ellipticit{\'e} des {\'e}quations de
Seiberg-Witten modulo l'action du groupe de jauge montre que la convergence est en r{\'e}alit{\'e}  $C^\infty$
sur tout compact de $M$ (cf. \cite{Mg}).

\subsection{Cas irrationnel}
\label{secirrat}
Nous allons maintenant expliquer comment g{\'e}n{\'e}raliser le
th{\'e}or{\`e}me \ref{theoA} dans le cas d'une m{\'e}trique
asymptotique {\`a} un 
mod{\`e}le $\bparabarg$ dont les poids sont tels que
v{\'e}rifient $\alpha_2-\alpha_1 \not \in \Rat$. 
Dans ce cas, le rev{\^e}tement
$p : I_a\times \CProj \rightarrow \bparab$ ne se factorise pas
(cf. lemme \ref{lemmemodrat} et corollaire \ref{corcomprat}),
la 
compactification naturelle du bout parabolique 
 devrait \^etre un quotient par $\ZZ$ sur lequel nous ne disposons
pas d'un th\'eor\`eme de l\'indice  appropri{\'e} pour d{\'e}velopper une 
th{\'e}orie de Seiberg-Witten (cf. \cite{Ka} pour les orbifolds).

Pour {\'e}viter ce probl{\`e}me, nous allons faire <<bouger>> les
poids, afin d'obtenir des mod{\`e}les $\gmodjj$ 
v{\'e}rifiant la condition de rationalit{\'e} et qui approximent
$\gmod$ sur tout compact de $M$. En recollant $g$ aux diff{\'e}rents mod{\`e}les $\gmodjj$,
nous en d{\'e}duirons des approximations 
par des m{\'e}triques $g^j$~; par ce qui pr{\'e}c{\`e}de nous pourrons
extraire  des solutions des
{\'e}quations de Seiberg-Witten pour ces m{\'e}triques que nous
ferons converger {\`a} leur tour vers une solution des {\'e}quations pour la m{\'e}trique $g$.

\subsubsection{Approximation  dans le cas irrationnel} 
Rappelons que la construction du mod{\`e}le local se fait via le choix
d'une m{\'e}trique hermitienne
$$h= \left (\begin{array}{cc}
  |z| ^ {\alpha_1} & 0 \\
  0 & |z| ^ {\alpha_2 }
  \end{array} \right ) \; 
$$
et il appara{\^\i}{}t clairement que $\gmod$ d{\'e}pend de fa{\c c}on
$C^\infty$ du choix des poids sur un compact de $M$.
On choisit  donc une suite de poids $(\alpha_1 ^ j, \alpha
^j_2)$ qui tend vers
$(\alpha_1, \alpha_2)$, telle que $\alpha_2^j -\alpha_1^j\in \Rat$ et
on note 
$\gmodjj$ le mod\`ele local associ\'e. En reprenant les notations de
\ref{subapproxmetraparab}, on pose
$$g^j = (1- \chi_j )g + \chi_j \gmodjj.
$$
On fait en sorte que les poids convergent  suffisamment vite de sorte
 que les m{\'e}triques 
 $\gmodjj$ soient tr{\`e}s proches de $\gmod$
au sens $C^2$ sur l'anneau $j\leq t\leq j+1$.

\remarque la m{\'e}trique $g^j$ est {\'e}videmment \aparab $\gmodjj$ et
tout ce qui concerne le cas
rationnel s'applique {\`a} elle.
Par ailleurs, d'apr{\`e}s le corollaire~\ref{cohomologie}, on a les 
isomorphismes
$$\HLLg{g^j}^*(M)\simeq \HDR^*(\overline M^j) \simeq \HLLg{g}^*(M),$$
via lesquels $\LLc(M,g)= \LLc (M, g^j)= c_1(K^{-1}_{\overline M^j} \otimes [D_j]^{-1})$.

\subsubsection{Convergence des $1$-formes harmoniques}
On a une in{\'e}galit{\'e} de Poincar{\'e} uniforme pour les fonctions qui
est l'analogue du lemme \ref{laplacien}~:
\begin{lemme}
\label{laplacienb}
  Il existe une constante $c>0$,  telle que pour toute m{\'e}trique $g^j$
  et toute fonction $f\in L ^2_1(g^j)$ v{\'e}rifiant $\int_M f\voljj
  =0$ on ait 
$$
    \int_M|df|^2 \voljj\geq c\int_M|f|^2\voljj.
$$
\end{lemme}
\begin{demo}
  On refait la d{\'e}monstration du lemme~\ref{laplacien} {\`a}
  l'identique. Il n'y a 
  que deux points {\`a} v{\'e}rifier pour que les arguments s'appliquent
  bien~:
 le lemme \ref{volume} est valable pour les m{\'e}triques
    $g^j$ et les constantes qui interviennent dans le
  lemme~\ref{courburemoyenne} ne d{\'e}pendent que de  la constante
  $h_0$ qui est la m{\^e}me pour toutes les m{\'e}triques $\gmodjj$, 
ce qui n'est pas difficile.
\end{demo}

On d{\'e}duit du lemme~\ref{laplacienb} par une d{\'e}monstration
presque identique {\`a} celle de la proposition~\ref{convcompact} 
l'analogue suivant~:
\begin{prop}
\label{convcompactb}
Soit $b\in \HLL^1(M)$ et $\gamma_j$  la suite de ses repr{\'e}sentants
$g^j$-harmoniques $L^2$. Alors $\gamma_j$  converge au sens 
$C^\infty$ sur tout compact de $M$, vers le repr{\'e}sentant
$g$-harmonique $L^2$ de $b$. 
\end{prop}

\subsubsection{In{\'e}galit{\'e} de Poincar{\'e} pour les $1$-formes}
Afin de g{\'e}n{\'e}raliser les r{\'e}sultats de convergence obtenus pour
les $1$-formes dans le
cas rationnel, on commence par montrer que  la constante $T$
du lemme~\ref{courburemoyenne2} peut {\^e}tre choisie ind{\'e}pendamment
de la m{\'e}trique $\gmodjj$. Notons $\gmodjj_k$ les approximations
successives suivant la m{\'e}thode expos{\'e}e au \S\ref{seccomp} de
chaque m{\'e}trique $\gmodjj$, on  a alors le lemme suivant~:
\begin{lemme}
    Il existe une constante $T$ suffisamment grande telle que le lemme
  \ref{courburemoyenne2} soit v{\'e}rifi{\'e} pour toutes les  m{\'e}triques
  $g=\gmod, \gmodjj$ ou $\gmodjj_k$.
\end{lemme}
\begin{demo}
  La seule modification {\`a} apporter dans la d{\'e}monstration du
  lemme \ref{courburemoyenne2} concerne l'argument utilisant les
  s{\'e}ries de Fourier~: comme 
$$(\alpha_2^j- \alpha_1^j)=\frac{r_j}{q_j}\rightarrow (\alpha_2-
\alpha_1) \not\in \Rat,$$
ceci entra{\^\i}{}ne que $q_j$ explose~; plus g{\'e}om{\'e}triquement,  les orbites de
$\Xtheta$ pour la m{\'e}trique $g$ limite ne se referment pas et s'enroulent
sur un tore. Nous allons donc raffiner l'argument en utilisant les
s{\'e}ries de Fourier sur le tore.

En reprenant les notations de la d{\'e}monstration du
lemme~\ref{courburemoyenne2}
la m{\'e}trique $g$ s'{\'e}crit en coordonn{\'e}es locales sur le bout
parabolique 
$$g=dt ^2 + \phi ^ 2(t)d\theta ^ 2 + \frac{4/c}{(1+|u|^2)^2}|du -i\alpha
u d\theta | ^ 2,\mbox{ o{\`u} $\alpha=\alpha_2-\alpha_1$.}
$$
En utilisant les coordonn{\'e}es polaires $u=\rho e ^{i\theta_2}$, on calcule
$$g=dt ^2 + e ^{-2t}d\theta ^ 2 + \frac{4/c}{(1+\rho^2)^2}\left (|d\rho | ^2
+ \rho ^ 2 |d\theta_2 -\alpha
 d\theta | ^ 2 \right ) 
$$
et $\Xtheta = \dtheta + \alpha
\dthetad$ -- notons que  $(d\theta,d\rho,d\theta_2)$ est une base duale
de  $(\Xtheta,\partial_\rho,\dthetad)$.
 Les orbites de $\Xtheta$ s'enroulent sur les tores
$$T(t_0,\rho_0)= \{
(t,\rho,\theta,\theta_2) / t=t_0 \mbox{ et } \rho=\rho_0 \}.$$
suivant lesquels on {\'e}crit la d{\'e}composition de $f$ en s{\'e}rie de Fourier
$$f=\sum_{k,l} c_{k,l}(t,\rho) e ^ {i(k\theta + l \theta _2)}\;;$$
on en d{\'e}duit que
$$\Xtheta \cdot f=(\dtheta + \alpha
\dthetad) \cdot f =\sum_{k,l} ic_{k,l}(k+\alpha l) e ^ {i(k\theta + l
\theta _2)}.$$ 
La difficult{\'e} dans le cas o{\`u} $\alpha$ est irrationnel provient du
fait que l'ensemble des valeurs prises par $k+l\alpha$ s'accumule en
$0$. Notons tout d'abord que
$$  l\rightarrow \infty \mbox{  et  }
k\rightarrow \infty \mbox{ lorsque }  k+l\alpha \rightarrow 0.
$$
\remarque
dans le cas de
la m{\'e}trique $\gmodjj$, on a  
$\alpha_2^j-\alpha_1^j = \alpha ^j =\frac{r_j}{q_j}
\rightarrow \alpha \not\in \Rat$.
Alors, pour tout $K>0$, il existe un $\epsilon>0$ suffisamment petit, tel que
pour $j$ suffisamment grand on ait $|k+\alpha ^j l| ^ 2 \leq \epsilon
\Rightarrow k,l\geq K$. En effet, si cette assertion {\'e}tait fausse,
ceci impliquerait $\alpha\in\Rat$, ce qui est contradictoire. 
On v{\'e}rifiera en vertu de cette remarque que tout ce
que nous dirons dans la suite de cette d{\'e}monstration s'applique {\'e}galement aux m{\'e}triques
$\gmodjj$ et $\gmodjj_k$ en choisissant $j$ et $k$ suffisamment grands.
\medskip

Maintenant
$$ \dthetad  \cdot f = \sum_{k,l} il c_{k,l}  e ^ {i(k\theta + l
\theta _2)}$$
et puisque $d^{\CP} f =  \dthetad  \cdot f (d\theta_2-\alpha d\theta) +
\partial_\rho \cdot f d\rho $
est une somme de termes
orthogonaux, on en d{\'e}duit que  que
\begin{equation}
\label{cc1}
|d^{\CP} f|^2\geq  4/c(\rho+ 1/\rho)^2
|\dthetad \cdot f| ^ 2 \geq   4/c
|\dthetad \cdot f| ^ 2.
\end{equation}

On d{\'e}compose l'espace des fonctions sur le tore en
l'espace des fonctions presque invariantes suivant $\Xtheta$ et son
orthogonal 
$$\Inv=\langle e ^{i(k\theta_1 +l\theta_2)}  \; / \;|k+\alpha l|^ 2\leq
\epsilon \rangle ,\quad \Inv^\perp=\langle e ^{i(k\theta_1 +l\theta_2)}  \; / \;\epsilon < |k+\alpha
l|^ 2 
 \rangle .$$
Pour tout $K>0$ il existe un
$\epsilon$ suffisamment petit, tel que la d{\'e}composition de $f\in
\Inv$  en s{\'e}rie de Fourier s'{\'e}crive
$$f=\sum_{k,l\geq K} c_{k,l} e ^{i(k\theta+l\theta_2)}.$$
\begin{equation}
\label{cc2}
\mbox{d'o{\`u}}\quad\int_N | \dthetad \cdot  f |^2 \volgt \geq K ^
2\int_N |f|^2\volgt.
\end{equation}
Revenons  {\`a} l'in{\'e}galit{\'e} (\ref{exprintermpoinc1}) de la
d{\'e}monstration du lemme \ref{courburemoyenne2}~:
\begin{multline*}
\int_N |\nabla \beta |^2 + \Ric^g(\beta,\beta) \geq \int_N |\nabla_\dt \beta
|^ 2 +
|\phi ^{-1} \Xtheta  \cdot  f|^2 +  |\phi ^{-1}\Xtheta \cdot f_2|^2 \\ 
 - 4 \langle f , \phi ^{-1} \Xtheta
\cdot f_2 \rangle + |d^{\CProj ^1}f|^2+ |d^{\CProj ^1}f_2|^2 \; ;
\end{multline*}
on voit donc que pour $f\in \Inv$, le terme mixte $- 4 \langle f , \phi
^{-1} \Xtheta \cdot f_2 \rangle$ est contr{\^o}l{\'e}  {\`a} 
l'aide de  $|\phi
^{-1}\Xtheta \cdot f_2|^2$ et $|d^{\CProj^1 }f|^2$ d'apr{\`e}s (\ref{cc1})
et (\ref{cc2}).

Supposons  $f\in \Inv^\perp$. Alors
$$\int_N |\Xtheta  \cdot  f| \volgt\geq \epsilon\int_N |f|^ 2 \volgt ,
$$ 

$$\mbox{d'o{\`u}}\quad \int_N |\phi ^{-1} \Xtheta \cdot f| \volgt\geq \epsilon \phi ^{-2}\int_N |f|^ 2 \volgt .
$$ 
Puisque $\phi ^ {-2}(t)\geq e ^{2t}\rightarrow +\infty$, on en d{\'e}duit
que pour $t\geq T$, suffisamment grand, le terme mixte est contr{\^o}l{\'e} par $|\phi
^{-1}\Xtheta \cdot f|^2$ et  $|\phi
^{-1}\Xtheta \cdot f_2|^2$. 
\end{demo}

On en  d{\'e}duit l'analogue de la proposition
 \ref{laplacien2}~:
\begin{prop}
\label{laplacien2b}
  Il existe une constante $c>0$ telle que pour toutes les m{\'e}triques
  $g^j$, et toute $1$-forme $\beta\in L^2_1(g^j)$, orthogonale aux
  formes $g^j$-har\-mo\-ni\-ques $L^2$,  on ait~:
$$\int_M |d\beta|^2+|d^{*_j}\beta|^2 \; \voljj \geq c\int_M|\beta|^2\voljj.
$$
\end{prop}\medskip

\subsubsection{Convergence des $2$-formes harmoniques}
Soit une classe de cohomologie $b\in \HLL^2(M)$. {\`A} l'aide du lemme
de Poincar{\'e}~\ref{lemmePoincare}, on commence par  choisir un 
repr{\'e}sentant $\gamma\in L^2(g)$ {\'e}gal {\`a}  
$\lambda\volss$ sur le bout parabolique.
On en d{\'e}duit  une suite  de repr{\'e}sentants $\gamma_j\in
L^2(g^j)$ via l'isomorphisme du du corollaire~\ref{cohomologie} tels que $\|\gamma_j\NLLjj$ soit born{\'e} et  $\gamma_j$ converge vers $\gamma$ au sens $C^\infty$ sur tout
compact de $M$.
 
On peut {\'e}crire le
repr{\'e}sentant $g^j-$harmonique $w_j \in L^2(g^j)$ de $b$ sous la
forme $w_j = \gamma_j + d\beta_j$,
o{\`u} $\beta_j\in L^2_1(g^j)$ est une $1$-forme orthogonale aux formes
$g^j$-harmoniques $L^2$ qui v{\'e}rifie $d^{*_j}\beta_j = 0$. 
Gr{\^a}ce {\`a} l'in{\'e}galit{\'e} de Poincar{\'e} de la proposition
\ref{laplacien2b},
on en d{\'e}duit par une d{\'e}monstration similaire {\`a} celle de la
proposition \ref{propconvharm}~:
\begin{prop} 
\label{propconvharmb}
Soit une classe de cohomologie $b\in \HLL(M)$. Soit $\Omega_j$ la
suite de ses repr{\'e}sentants $g^j$-harmoniques $L^2$. Alors
$\Omega_j$ converge au sens $C^\infty$ sur tout compact de
$M$ vers le repr{\'e}sentant $g$-harmonique $L^2(g)$ de $b$.
\end{prop}
Chaque m{\'e}trique $g^j$ induit une d{\'e}composition des
formes $g^j$-harmoniques $L^2$ en leur partie autoduale et 
anti-autoduale,
$$\HLL^2(M)= \HLL^{+^j}(M)\oplus  \HLL^{-^j}(M).$$
On d{\'e}montre {\'e}galement facilement
l'analogue du corollaire \ref{convpol}  
\begin{cor}
\label{convpolb}
  Soit une classe de cohomologie $b\in \HLL^2(M)$, $\pm ^j$ les
  projecteurs de $\HLL^2(M)$ associ{\'e}s {\`a} la m{\'e}trique $g^j$ et
  $\pm$ ceux associ{\'e}s 
  {\`a} $g$~; alors
$$ b^{+^j} \rightarrow b^+\quad \mbox{ et }  \quad b^{-^j} \rightarrow b^-,$$
lorsque $j$ tend vers l'infini. 
\end{cor}
\subsubsection{Convergence des
  solutions des {\'e}quations de Seiberg--Witten pour les m{\'e}triques $g^j$}
\label{subconvirrat}
Il faut commencer par fixer des connexions de r{\'e}f{\'e}rence $C^j$
convenables sur les fibr{\'e}s 
$$L_j = K^{-1}_{\overline M^j}\otimes
[D^j]^{-1}$$
 Soit $\Amod$ la connexion
induite par $\gmod$ sur $K_M^{-1}$ et $\Amod^j$ celles induites par
les m{\'e}triques $\gmodjj$. Alors les connexion $\Amod^j$  convergent
vers $\Amod$ au sens $C^1$ sur tout
compact de $M$ par construction de $\gmodjj$. 
D'apr{\`e}s  le  lemme \ref{lemmeconnl2}, on peut donc construire  alors une suite de
connexions $C^j$  qui converge sur tout compact  
$M$ vers la connexion $\Amod$, lisses sur  $L_j$ et telles que $\|C^j- \Amod^j \NLLj$
et $\|F_{C^j}- F_{\Amod^j} \NLLj$ soient born{\'e}s.
\medskip 

On d{\'e}montre alors a l'analogue du
th{\'e}or{\`e}me~\ref{theoconvergence} suivant~:
\begin{theo}
\label{theoconvergenceb}
  Soit $(A ^j, \psi ^j)$, une suite  de solutions des {\'e}quations de
  Seiberg-Witten non perturb{\'e}es  associ{\'e}es aux m{\'e}triques
  $g^j$  telles que
  \begin{itemize}
  \item $A ^j = C ^ j +  a ^j,$ o{\`u} $a ^j\in L^2_1(g^j)$ est en
    jauge de Hodge $d^{*_j} a ^j=0$ et la partie $g^j$-harmoniques
    $L^2$ de $a^j$ est uniform{\'e}ment born{\'e}e  en cohomologie.
   \item $\psi ^j\in L^2_1(g^j)$ relativement {\`a} la connexion $C^j$. 
   \item $\psi ^j$ est uniform{\'e}ment born{\'e}e en norme $C^0$.
  \end{itemize}
Alors,  quitte {\`a} extraire une sous suite, $(A ^j,\psi ^j)$ converge
au sens $C^\infty$ sur tout compact de $M$ vers une solution $(A
,\psi)$ des
{\'e}quations de Seiberg-Witten non perturb{\'e}es pour la m{\'e}trique $g$
telle que
\begin{itemize}
\item $ A= C+a$, avec $d^*a = 0$ et $a\in L^2_1(g)$.
\item $\psi\in L^2_1(g)$ relativement {\`a} $C$.
\end{itemize}
\end{theo}
\begin{demo}
 Notons que tout d'abord que la borne $C^0$ du
lemme~\ref{lemmec0}
ne d{\'e}pend pas non plus de $j$ pour les m{\'e}triques $g^j_k$ 
car la courbure scalaire des m{\'e}triques $g^j$ est
uniform{\'e}ment born{\'e}e. 
Les hypoth{\`e}ses de r{\'e}gularit{\'e} sur 
$(A ^j,\psi ^j)$ sont donc justement celles qu'on
peut faire si $(A ^j,\psi ^j)$ a {\'e}t{\'e} obtenue {\`a} l'aide du
th{\'e}or{\`e}me  \ref{theoconvergence}.

Par ailleurs on a une  borne uniforme sur $[\mu ^j]$, o{\`u} $\mu^j$ est
la partie $g^j$-harmonique de $a^j$
car on s'est ramen{\'e} {\`a} un domaine fondamental du tore
$H^1(\overline M^j,\R) / H^1(\overline M^j,\ZZ)$ afin de faire
converger $\mu ^j_k$.
On peut quitte {\`a}
  extraire une sous-suite supposer que $[\mu ^j]$ converge et en
  d{\'e}duire que $\mu ^j$ converge vers une forme $g$-harmonique
$L^2$ gr{\^a}ce {\`a} la proposition \ref{propconvharmb}. 

  Il est facile de faire converger $\beta ^j$ et $\psi ^j$
  (cf. th. \ref{theoconvergence}) maintenant
  que nous disposons de l'in{\'e}galit{\'e} de
  Poincar{\'e} de la proposition~\ref{laplacien2b}
 et de connexions de r{\'e}f{\'e}rence $C^j$
ad{\'e}quates.
\end{demo}

\subsection{Calcul de l'invariant de Seiberg--Witten}
\label{subcalcinv}
Soit $M$ une surface complexe r{\'e}gl{\'e}e obtenue {\`a} partir d'un fibr{\'e}
parabolique $\mathcal E$. Il nous reste prouver l'existence de
solutions des {\'e}quations de 
Seiberg-Witten sur les compactifications de $M$. Pour cela, on
commence par choisir une m{\'e}trique k{\"a}hl{\'e}rienne particuli{\`e}re quitte {\`a}
perturber la structure complexe sur  $\overline M$ pour laquelle on
sait que l'espace des module est constitu{\'e} d'un point sous la
condition~$\deg L<0$. Comme nous somme dans le cas o{\`u} $b_2^+=1$,
l'invariant de Seiberg-Witten ne  d{\'e}pend pas de la  m{\'e}trique sous
une <<condition de
chambre>>~(\ref{wallcross}) ce qui nous permet de trouver une solution
des {\'e}quations pour ces m{\'e}triques. 

\subsubsection{Cas rationnel}
\label{subsubcasrat}
 Si les poids de  la structure parabolique
v{\'e}rifient la condition de rationnalit{\'e}, on en d{\'e}duit une
compactification orbifold $\overline M$. Comme la condition de
stabilit{\'e} est g{\'e}n{\'e}rique, on peut quitte {\`a} perturber la
structure holomorphe de $\mathcal E$ le supposer stable. On en
d{\'e}duit une m{\'e}trique k{\"a}hl{\'e}rienne $\gmod$ {\`a} courbure scalaire
constante $s=2(c-1)<0$ d{\'e}finie sur $M$. Par construction, les
approximations $\gmodj$ de $\gmod$ sont k{\"a}hl{\'e}riennes sur
$\overline M$ et on la m{\^e}me classe de K{\"a}hler $[ \omega]$.
Comme $\gmod$ est k{\"a}hl{\'e}rienne {\`a} courbure scalaire constante,
\begin{equation}
\label{sscal}
\LLc\cdot[\omega] = \frac s{8\pi}\int \omega\wedge\omega,
\end{equation}
donc $\deg L = \LLc \cdot [\omega ] <0$ ce entra{\^\i}{}ne  que
l'espace des modules  des {\'e}quations de Seiberg--Witten 
perturb{\'e}es~(\ref{eq:SWp}) sur $(\overline M,\gmodj)$ s'identifie
{\`a} l'espace des modules de couples $(A^{0,1}, \alpha )$ o{\`u}
$A^{0,1}$ est une structure holomorphe sur le fibr{\'e}  $L_0$ et
$\alpha$ une section holomorphe du fibr{\'e} trivial $L_0\otimes
K_{\overline M}$ (voir \cite{W})~; la seule solution est le fibr{\'e} $K^{-1}_{\overline M}$
avec $\alpha$ la section holomorphe standard du fibr{\'e} trivial $\mathcal
O_{\overline M}$. 
On en d{\'e}duit que pour toute m{\'e}trique $g$ sur
$\overline M$ v{\'e}rifiant 
\begin{equation}
\label{wallcross}
\LLc ^{+_g}\cdot [\omega]<0,
\end{equation}
 les {\'e}quations
de Seiberg--Witten perturb{\'e}es on un invariant $\SW(g)= 1$.

Soit $\gkahl$ une m{\'e}trique k{\"a}hl{\'e}rienne {\`a} courbure scalaire
constante $s$ \aparab $\gmod$ sur $M$. En notant $\omega^K$ sa forme
de K{\"a}hler, 
$$\int_{\CProj^1} \omega ^{\mathrm K} = \int_{\CProj^1} \omega
$$
en raison du comportement asymptotique de $\gkahl$~; d'apr{\`e}s~(\ref{sscal}),
$$
 \frac s{8\pi}= \frac {\LLc\cdot [\omega]}{ [\omega]^2} = \frac
{\LLc \cdot [\omega ^
{ K}]}{ [\omega ^ K]^2},
$$
ce qui implique $[\omega] = [\omega ^{\mathrm K}]$.

Soit $\gkahl$ une m{\'e}trique  k{\"a}hl{\'e}rienne sur $M$ {\`a} courbure
scalaire constante $s<0$ \aparab (dont les poids v{\'e}rifient la
condition de rationnalit{\'e})
Notons $g=\xi ^*\gkahl$, o{\`u} $\xi$ est l'antipodie sur $M$
d{\'e}finie au \S\ref{subdemotheoa} et $g_j$ la suite d'approximations
de $g$. 
D'apr{\`e}s le corollaire \ref{convpol} et le lemme \ref{corpol}
$$\LLc ^{+_j} \rightarrow\LLc ^{+_{g}} =\LLc ^{+_{\gkahl}},
$$
donc  d'apr{\`e}s (\ref{sscal}) pour $j$ suffisamment grand,
$$\LLc ^{+_j}\cdot[\omega]<0 \; ;
$$
on en d{\'e}duit que $\SW(g_j)=1$ et que les {\'e}quations de
Seiberg--Witten perturb{\'e}es associ{\'e}es {\`a} la 
m{\'e}trique $g_j$ admettent une solution $(A_j,\psi_j)$.
En utilisant le th{\'e}or{\`e}me~\ref{theoconvergence}, on en d{\'e}duit la
solution recherch{\'e}e des {\'e}quation de Seiberg--Witten non
perturb{\'e}es associ{\'e}es {\`a} la m{\'e}trique $g=\xi ^*\gkahl$.

\subsubsection{Cas irrationnel}
Dans le cas o{\`u} le mod{\`e}le local $\gmod$ ne v{\'e}rifie
plus la condition $\alpha_2-\alpha_1\in \Rat$, on approxime la
m{\'e}trique $g=\xi ^*\gkahl$ par les m{\'e}triques $g^j$ d{\'e}finies au
\S\ref{secirrat} qu'on approxime {\`a} leur tour
par des m{\'e}triques lisses $g^j_k$ sur des compactifications orbifold
$\overline M^j$ suivant le \S\ref{seccomp}.
D'apr{\`e}s les corollaires \ref{convpol} et \ref{convpolb}, en notant
$\pm^j_k$ les projecteurs des d{\'e}compositions de $\HLL ^2(M)$ induites
par les m{\'e}triques $g_k^j$, on
d{\'e}duit que
$$\LLc ^{+^j_k} \rightarrow\LLc ^{+_{g}} =\LLc ^{+_{\gkahl}}
$$
pour $j$ et $k$ tendant vers l'infini. On en d{\'e}duit, en
utilisant le \S\ref{subsubcasrat} que  les {\'e}quations de
Seiberg--Witten non perturb{\'e}es 
associ{\'e}es {\`a} la m{\'e}trique $g^j$ admettent une solution via le
th{\'e}or{\`e}me~\ref{theoconvergence}. On fait converger {\`a} nouveau
cette suite vers une solution des {\'e}quations pour la m{\'e}trique
$g=\xi ^*\gkahl$ en utilisant le th{\'e}or{\`e}me~\ref{theoconvergenceb}.

\bigskip

\end{document}